\documentclass[A4j,11pt]{article}
\usepackage{amsfonts}
\usepackage{graphics}
\begin{document}

\title{{\bf The homogeneous slice theorem\\
for the complete complexification of\\
a proper complex equifocal submanifold}}
\author{{\bf Naoyuki Koike}}

\date{}
\maketitle

{\small\textit{Department of Mathematics, Faculty of Science, 
Tokyo University of Science,}}

{\small\textit{26 Wakamiya-cho Shinjuku-ku, Tokyo 162-8601 Japan}}

{\small\textit{E-mail address}: koike@ma.kagu.tus.ac.jp}

%
%
%

%
\begin{abstract}
The notion of a complex equifocal submanifold 
in a Riemannian symmetric space of non-compact type has been recently 
introduced as a generalization of isoparametric hypersurfaces in the 
hyperbolic space.  As its subclass, the notion of a proper complex equifocal 
submanifold has been introduced.  
Some results for a proper complex equifocal submanifold has been recently 
obtained by investigating the lift of its complete complexification to some 
path space.  
In this paper, we give a new construction of the complete complexification 
of a proper complex equifocal submanifold and, by using the construction, show 
that leaves of focal distributions of the complete complexification are 
the images by the normal exponential map of principal orbits of a certain kind 
of pseudo-orthogonal representation on the normal space of the corresponding 
focal submanifold.  
\end{abstract}

\vspace{0.3truecm}

{\small\textit{Keywords}: $\,$ proper complex equifocal submanifold, 
anti-Kaehlerian holonomy system,} 

\hspace{2truecm}{\small aks-representation}

\vspace{0.3truecm}

{\small\textit{MS classification}: 53C35, 53C40.}

\vspace{0.5truecm}

\section{Introduction}
C.L. Terng and G. Thorbergsson [TT1] introduced the notion of an equifocal 
submanifold in a Riemannian symmetric space , which is defined as a compact submanifold 
with globally flat and abelian normal bundle such that the focal radii for 
each parallel normal vector field are constant.  
This notion is a generalization of isoparametric submanifolds in the Euclidean 
space and isoparametric hypersurfaces in the sphere or the hyperbolic space.  
For (not necessarily compact) submanifolds in a Riemannian symmetric space  of non-compact type, the equifocality is a rather weak property.  So, we [K1,2] 
introduced the notion of a complex focal radius as 
a general notion of a focal radius and defined the notion of a complex 
equifocal submanifold as a submanifold with globally flat and abelian normal 
bundle such that the complex focal radii for each parallel normal vector field 
are constant and that they have constant multiplicties.  
E. Heintze, X. Liu and C. Olmos [HLO] defined the notion of an isoparametric 
submanifold with flat section as a submanifold with globally flat and abelian 
normal bundle such that the sufficiently close parallel submanifolds are of 
constant mean curvature with respect to the radial direction.  
The following fact is known (see Theorem 15 of [K2]):

\vspace{0.2truecm}

{\sl All isoparametric submanifolds with flat section are complex equifocal 
and, conversely, all curvature-adapted and complex equifocal submanifolds are 
isoparametric submanifolds with flat section.}

\vspace{0.2truecm}

\noindent
Furthermore, as its subclass, we [K1,2] introduced the notion of 
a proper complex equifocal submanifold.  
For a proper complex equifocal submanifold, the following fact is known ([K3]):

\vspace{0.2truecm}

{\sl Principal orbits of Hermann type actions on a Riemannian symmetric space 
of non-compact type are curvature-adapted and proper complex equifocal.}

\vspace{0.2truecm}

For a (general) submanifold in a Riemannian symmetric space  of non-compact 
type, the (non-real) complex focal raidii are defined algebraically.  We 
needed to find their geometrical essence.  For its purpose, we defined the 
complexification of the ambient Riemannian symmetric space  and defined the 
extrinsic complexification of the submanifold as a certain kind of submanifold 
in the complexified symmetric space, where the original submanifold needs to 
be assumed to be complete and real analytic.  
In the sequel, we assume that all submanifolds in the Riemannian symmetric 
space are complete and real analytic.  
We [K2] showed that the complex focal radii of the original submanifold 
indicate the positions of the focal points of the complexified submanifold.  
If the original submanifold is complex equifocal, then the extrinsic 
complexification is an anti-Kaehlerian equifocal submanifold in the sense of 
[K2].  Also, if the original submanifold is proper complex equifocal, 
then the complexified one is a proper anti-Kaehlerian equifocal submanifold 
in the sense of this paper.  Thus, the study of an anti-Kaehlerian equifocal 
(resp. proper anti-Kaehlerian equifocal) submanifold leads to that of a 
complex equifocal (resp. proper complex equifocal) submanifold.  
The complexified submanifold is not necessarily complete.  
In the global research, we need to extend the complexified submanifold 
to a complete one.  
In [K2], we obtained the complete extension of the complexified submanifold 
in the following method.  We first lifted the complexified submanifold to some 
path space (which is an infinite dimensional anti-Kaehlerian space) through 
some submersion, extended the lifted submanifold to the complete one 
by repeating some kind of extension infinite times and obtained the complete 
extension of the original complexified submanifold as the image of the 
complete one by the submersion.  
In this paper, we give a new construction of the complete extension of the 
complexified submanifold without repeating infinite times of processes 
(see the proof of Theorem B) and investigate the detailed structure of the 
complete extension in terms of the new construction.  
First we prove the following fact 
for an anti-Kaehlerian equifocal submanifold.  

\vspace{0.5truecm}

\noindent
{\bf Theorem A.} {\sl Let $M$ be an anti-Kaehlerian equifocal submanifold 
in a semi-simple anti-Kaehlerian symmetric space of non-positive 
(or non-negative) curvature having a focal submanifold 
$F$.  If the sections of $M$ are properly embedded, 
then $M$ is an open potion of a partial tube over $F$ whose each fibre is 
the image by the normal exponential map of a principal orbit of 
a pseudo-orthogonal representation on the normal space of $F$ which is 
equivalent to the direct sum representation of an aks-represenation and a 
trivial representation.}

\vspace{0.5truecm}

\noindent
{\it Remark 1.1.} 
(i) For a focal submanifold $F$ of $M$, we call 
$(\exp^{\perp}\vert_{T^{\perp}_xF})^{-1}(\exp^{\perp}(T^{\perp}_xF)\cap M)$ 
(rather than $\exp^{\perp}(T^{\perp}_xF)\cap M$) the slice of $M$.  This 
theorem states that slices of a complete anti-Kaehlerian equifocal submanifold 
are homogeneous.  

(ii) The dual action $H^{\ast}$ of a Hermann type action $H$ 
on a Riemannian symmetric space  $G/K$ of non-compact type is a Hermann action on 
the compact dual $G^{\ast}/K$ of $G/K$, where $G$ is assumed to be a 
connected semi-simple Lie group admitting a faithful real representation.  
Note that the existence of the dual action $H^{\ast}$ is assured by replacing 
$H$ by the conjugate group if necessary.  
Hence the sections of the $H^{\ast}$-action are flat tori.  From this fact, 
we see that the sections of the $H^{\bf c}$-action on 
$G^{\bf c}/K^{\bf c}$ are properly embedded, where $H^{\bf c}$ is the complexification of $H$ and $G^{\bf c}/K^{\bf c}$ is the anti-Kaehlerian symmetric 
space associated with $G/K$.  On the other hand, 
the principal orbits of the 
$H^{\bf c}$-action are proper anti-Kaehlerian equifocal.  Thus the principal 
orbits are submanifolds as in the statement of Theorem A.  

(iii) This result is an analogy of that of M. Br$\ddot u$ck 
[B] for an equifocal submanifold in a simply connected Riemannian symmetric 
space of compact type.  

\vspace{0.5truecm}

In [K4,5], we proved some global results for a proper complex equifocal 
submanifold by investigating the lift of the complete complexification 
of the submanifold to some path space through some submersion.  
Thus, in the global study of a proper 
complex equifocal submanifold, it is important to investigate the detailed 
structure of its complete complexification.  
By using Theorem A, we obtain a new construction of the complete 
complexification of a proper complex equifocal submanifold 
(see the proof of Theorem B).  
From the construction and Theorem A, we obtain the following homogeneous slice 
theorem for the complete complexification of a proper complex equifocal 
submanifold.  

\vspace{0.5truecm}

\noindent
{\bf Theorem B.} {\sl Assume that the sections of the complexification of 
a proper complex equifocal submanifold $M$ in a Riemannian 
symmetric space $G/K$ of non-compact type are properly embedded.  
Then the following statements $({\rm i})$ and $({\rm ii})$ hold:

$({\rm i})$ Each leaf of any focal distribution of the complete 
complexification 
$\widehat{M^{\bf c}}$ of $M$ is the image by the normal exponential map of 
a principal orbit of a pseudo-orthogonal representation 
on the normal space of a focal submanifold 
which is equivalent to the direct sum representation of an aks-representation 
and a trivial representation.  

$({\rm ii})$ Let $E_0$ be the disribution on $\widehat{M^{\bf c}}$ defined by $(E_0)_x:=\displaystyle{\mathop{\cap}_{v\in T^{\perp}_x\widehat{M^{\bf c}}}
\left({\rm Ker}\,R^{\bf c}(\cdot,v)v\cap{\rm Ker}\,A^{\bf c}_v\right)}$ 
($x\in\widehat{M^{\bf c}}$), where $R^{\bf c}$ is the curvature tensor of 
$G^{\bf c}/K^{\bf c}$ and $A^{\bf c}$ is the shape tensor of 
$\widehat{M^{\bf c}}$.  Then there exists a family 
$\{E_i\,\vert\,i=1,\cdots,k\}$ of focal distributions of $\widehat{M^{\bf c}}$ 
such that the leaves of $E_i$ ($i=1,\cdots,k$) are the images by the normal 
exponential map of complex spheres in the normal spaces of focal submanifolds 
and that $E_0\oplus\sum\limits_{i=1}^kE_i=T\widehat{M^{\bf c}}$ holds.}

\vspace{0.5truecm}

For a curvature-adapted and proper complex equifocal submanifold, we obtain 
the following fact in terms of Theorem B.  

\vspace{0.5truecm}

\noindent
{\bf Theorem C.} {\sl Let $M$ be a proper complex equifocal submanifold 
in a Riemannian symmetric space of non-compact type as in 
Theorem B and $\{E_0,\cdots,E_k\}$ be as in the statement $({\rm ii})$ of 
Theorem B.  Assume that $M$ is curvature-adapted.  Then 
$E_i^{\bf R}:=E_i\vert_M\cap TM$ ($i=0,\cdots,k$) are integrable distributions 
on $M$, leaves of $E_i^{\bf R}$ are half-dimensional totally real submanifolds 
of leaves of $E_i$ and $TM=E_0^{\bf R}\oplus\sum\limits_{i=1}^kE_i^{\bf R}$, 
where $E_i\vert_M$ is the restriction of $E_i$ to $M$.}

\vspace{0.5truecm}

\noindent
{\it Remark 1.2.} B. Wu ([W2]) showed that leaves of curvature distributions 
of a complete isoparametric submanifold in a hyperbolic space are totally 
umbilic spheres, totally umbilic hyperbolic spaces or horoshperes, where we 
note that the complexifications of a totally umbilic sphere and 
a totally umbilic hyperbolic space are totally anti-Kaehlerian umbilic 
complex spheres in the complexification (which is a complex sphere) of 
the ambient hyperbolic space.  See [K2] about the definition of the totally 
anti-Kaehlerian umbilicity.  Thus the statement of Theorem C is interpreted as 
an analogy of this result by B. Wu.  

\vspace{0.5truecm}

\noindent
{\bf Future plan of research.} {\sl By using Theorem B, we will investigate 
whether the complete complexifications of proper complex equifocal 
submanifolds are homogeneous.  Also, by using Theorems B and C, we will 
investigate whether curvature-adapted and proper complex equifocal 
submanifolds are homogeneous.}

\section{Basic notions}
In this section, we recall basic notions introduced in [K1$\sim$3].  
We first recall the notion of a complex equifocal submanifold introduced 
in [K1].  Let $M$ be an immersed submanifold with abelian normal bundle 
(i.e., the sectional curvature for each 2-plane in the normal space is equal 
to zero) of in a symmetric space $N=G/K$ of non-compact type.  
Denote by $A$ the shape 
tensor of $M$.  Let $v\in T^{\perp}_xM$ and $X\in T_xM$ ($x=gK$).  Denote 
by $\gamma_v$ the geodesic in $N$ with $\dot{\gamma}_v(0)=v$.  
The strongly $M$-Jacobi field $Y$ along $\gamma_v$ with $Y(0)=X$ (hence 
$Y'(0)=-A_vX$) is given by 
$$Y(s)=(P_{\gamma_v\vert_{[0,s]}}\circ(D^{co}_{sv}-sD^{si}_{sv}\circ A_v))
(X),$$
where $Y'(0)=\widetilde{\nabla}_vY,\,\,P_{\gamma_v\vert_{[0,s]}}$ is 
the parallel translation along $\gamma_v\vert_{[0,s]}$ and 
$D^{co}_{sv}$ (resp. $D^{si}_{sv}$) is given by 
$$\begin{array}{c}
\displaystyle{
D^{co}_{sv}=g_{\ast}\circ\cos(\sqrt{-1}{\rm ad}(sg_{\ast}^{-1}v))
\circ g_{\ast}^{-1}}\\
\displaystyle{\left({\rm resp.}\,\,\,\,
D^{si}_{sv}=g_{\ast}\circ
\frac{\sin(\sqrt{-1}{\rm ad}(sg_{\ast}^{-1}v))}
{\sqrt{-1}{\rm ad}(sg_{\ast}^{-1}v)}\circ g_{\ast}^{-1}\right).}
\end{array}$$ 
Here ${\rm ad}$ is the adjoint 
representation of the Lie algebra $\mathfrak g$ of $G$.  
All focal radii of $M$ along $\gamma_v$ are obtained as real numbers $s_0$ 
with ${\rm Ker}(D^{co}_{s_0v}-s_0D^{si}_{s_0v}\circ A_v)\not=\{0\}$.  So, we 
call a complex number $z_0$ with ${\rm Ker}(D^{co}_{z_0v}-
z_0D^{si}_{z_0v}\circ A_v^{{\bf c}})\not=\{0\}$ a {\it complex 
focal radius of} $M$ {\it along} $\gamma_v$ and call ${\rm dim}\,
{\rm Ker}(D^{co}_{z_0v}-z_0D^{si}_{z_0v}\circ A_v^{{\bf c}})$ the 
{\it multiplicity} of the complex focal radius $z_0$, 
where $D^{co}_{z_0v}$ (resp. $D^{si}_{z_0v}$) 
is a ${\bf C}$-linear transformation of $(T_xN)^{\bf c}$ defined by 
$$\begin{array}{c}
\displaystyle{
D^{co}_{z_0v}=g^{\bf c}_{\ast}\circ\cos(\sqrt{-1}{\rm ad}^{\bf c}
(z_0g_{\ast}^{-1}v))\circ (g^{\bf c}_{\ast})^{-1}}\\
\displaystyle{\left({\rm resp.}\,\,\,\,
D^{si}_{sv}=g^{\bf c}_{\ast}\circ
\frac{\sin(\sqrt{-1}{\rm ad}^{\bf c}(z_0g_{\ast}^{-1}v))}
{\sqrt{-1}{\rm ad}^{\bf c}(z_0g_{\ast}^{-1}v)}\circ(g^{\bf c}_{\ast})^{-1}
\right),}
\end{array}$$
where $g_{\ast}^{\bf c}$ (resp. ${\rm ad}^{\bf c}$) is the complexification 
of $g_{\ast}$ (resp. ${\rm ad}$).  
Here we note that, in the case where $M$ is of class $C^{\omega}$, 
complex focal radii along $\gamma_v$ 
indicate the positions of focal points of the extrinsic 
complexification $M^{\bf c}(\hookrightarrow G^{\bf c}/K^{\bf c})$ of $M$ 
along the complexified geodesic $\gamma_{\iota_{\ast}v}^{\bf c}$, where 
$G^{\bf c}/K^{\bf c}$ is the anti-Kaehlerian symmetric space associated with 
$G/K$ and $\iota$ is the natural immersion of $G/K$ into 
$G^{\bf c}/K^{\bf c}$.  
See the following paragraph about the definitions of 
$G^{\bf c}/K^{\bf c},\,M^{\bf c}(\hookrightarrow G^{\bf c}/K^{\bf c})$ and 
$\gamma_{\iota_{\ast}v}^{\bf c}$.  
Also, for a complex focal radius $z_0$ of $M$ along $\gamma_v$, we 
call $z_0v$ ($\in (T^{\perp}_xM)^{\bf c}$) a 
{\it complex focal normal vector of} $M$ {\it at} $x$.  
Furthermore, assume that $M$ has globally flat normal bundle (i.e., 
the normal holonomy group of $M$ is trivial).  
Let $\tilde v$ be a parallel unit normal vector field of $M$.  
Assume that the number (which may be $\infty$) of distinct complex 
focal radii along $\gamma_{\tilde v_x}$ is independent of the choice of 
$x\in M$.  
Let $\{r_{i,x}\,\vert\,i=1,2,\cdots\}$ 
be the set of all complex focal radii along $\gamma_{\tilde v_x}$, where 
$\vert r_{i,x}\vert\,<\,\vert r_{i+1,x}\vert$ or 
"$\vert r_{i,x}\vert=\vert r_{i+1,x}\vert\,\,\&\,\,{\rm Re}\,r_{i,x}
>{\rm Re}\,r_{i+1,x}$" or 
"$\vert r_{i,x}\vert=\vert r_{i+1,x}\vert\,\,\&\,\,
{\rm Re}\,r_{i,x}={\rm Re}\,r_{i+1,x}\,\,\&\,\,
{\rm Im}\,r_{i,x}=-{\rm Im}\,r_{i+1,x}<0$".  
Let $r_i$ ($i=1,2,\cdots$) be complex 
valued functions on $M$ defined by assigning $r_{i,x}$ to each $x\in M$.  
We call these functions $r_i$ ($i=1,2,\cdots$) {\it complex 
focal radius functions for} $\tilde v$.  
We call $r_i\tilde v$ a {\it complex focal normal vector field for} 
$\tilde v$.  If, for each parallel 
unit normal vector field $\tilde v$ of $M$, the number of distinct complex 
focal radii along $\gamma_{\tilde v_x}$ is independent of the choice of 
$x\in M$, each complex focal radius function for $\tilde v$ 
is constant on $M$ and it has constant multiplicity, then 
we call $M$ a {\it complex equifocal submanifold}.  
Let $\phi:H^0([0,1],\mathfrak g)\to G$ be the parallel transport map for $G$.  
See Section 4 of [K1] about the definition of 
the parallel transport map.  
This map $\phi$ is a pseudo-Riemannian submersion.  
Let $\pi:G\to G/K$ be the natural projection.  
It follows from 
Theorem 1 of [K2] that, 
$M$ is complex equifocal if and only if each component of 
$(\pi\circ\phi)^{-1}(M)$ is complex isoparametric.  See Section 2 of [K1] 
about the definition of a complex isoparametric submanifold.  
In particular, if each component of $(\pi\circ\phi)^{-1}(M)$ is 
proper complex isoparametric (i.e., complex isoparametric and, for each unit 
normal vector $v$, the complexified shape operator $A^{\bf c}_v$ is 
diagonalizable with respect to a pseudo-orthonormal base), then 
we call $M$ a {\it proper complex equifocal submanifold}.  
For a complex equifocal submanifold, the following fact holds:

\vspace{0.3truecm}

{\sl For a curvature-adapted and complex equifocal submanifold $M$, 
it is proper complex equifocal submanifold if and only if it has no focal 
point of non-Euclidean type on the ideal boundary of the ambient symmetric 
space.}

\vspace{0.3truecm}

Here the curvature-adaptedness means that, for each unit normal vector $v$, 
the Jacobi operator $R(\cdot,v)v$ ($R\,:\,$the curvature tensor of $G/K$) 
prerserves the tangent space invariantly and it commutes with the shape 
operator $A_v$.  See [K6] about the detail of the notion of a focal point 
of non-Euclidean type on the ideal boundary.  

Next we recall the notions of an anti-Kaehlerian symmetric space associated 
with a symmetric space of non-compact type which was introduced in [K2].  
Let $J$ be a parallel complex structure on an even dimensional 
pseudo-Riemannian manifold $(M,\langle\,\,,\,\,\rangle)$ of half index.  If 
$\langle JX,JY\rangle=-\langle X,Y\rangle$ holds for every $X,\,Y\in TM$, 
then $(M,\langle\,\,,\,\,\rangle,J)$ is called an 
{\it anti-Kaehlerian manifold}.  



\vspace{0.2truecm}

\centerline{
\unitlength 0.1in
%
\hspace{5.8truecm}
}

\vspace{0.1truecm}

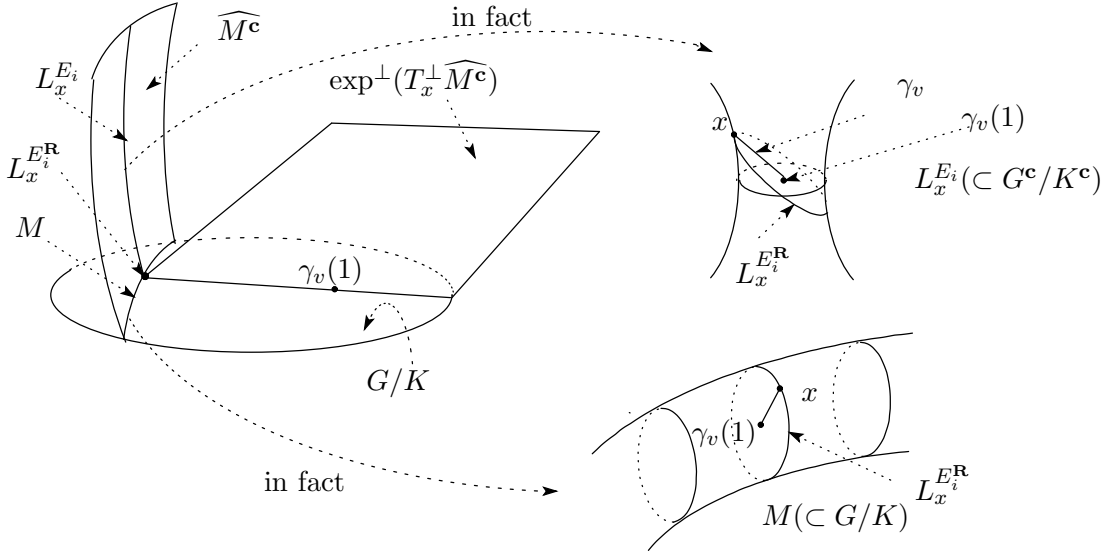
\centerline{{\bf Fig. 2.}}


\vspace{0.5truecm}

Let $N=G/K$ be a symmetric space of non-compact type and 
$({\mathfrak g},\sigma)$ be its orthogonal symmetric Lie algebra.  
Let ${\mathfrak g}={\mathfrak f}+{\mathfrak p}$ be the Cartan decomposition 
associated with a symmetric pair $(G,K)$.  
Note that $\mathfrak f$ is the Lie algebra of $K$ and $\mathfrak p$ is 
identified with the tangent space $T_{eK}N$, where $e$ is the identity 
element of $G$.  Let $\langle\,\,,\,\,\rangle$ be the 
${\rm Ad}(G)$-invariant non-degenerate inner product of $\mathfrak g$ 
inducing the Riemannian metric of $N$ 
and $\mathfrak a$ be a maximal abelian subspace of $\mathfrak p$ and 
$\mathfrak p=\mathfrak a+\sum\limits_{\alpha\in\triangle_+}
\mathfrak p_{\alpha}$ be 
the root space decomposition with respect to $\mathfrak a$, that is, 
${\mathfrak p}_{\alpha}=\{X\in{\mathfrak p}\,\vert\,{\rm ad}(a)^2(X)
=\alpha(a)^2X\,\,{\rm for}\,\,{\rm all}\,\,a\in\mathfrak a\}$.  
Let ${\mathfrak g}^{{\bf c}},
\,\,{\mathfrak f}^{{\bf c}},\,\,{\mathfrak p}^{{\bf c}},\,\,
{\mathfrak a}^{\bf c}\,\,{\mathfrak p}_{\alpha}^{\bf c}$ and 
$\langle\,\,,\,\,\rangle^{{\bf c}}$ be the complexifications of 
$\mathfrak g,\,\,\mathfrak f,\,\,\mathfrak p\,\,\mathfrak a\,\,
{\mathfrak p}_{\alpha}$ and $\langle\,\,,\,\,\rangle$, 
respectively.  If $\mathfrak g^{\bf c}$ and $\mathfrak f^{\bf c}$ are 
regarded as real Lie algebras, then 
$(\mathfrak g^{\bf c},\mathfrak f^{\bf c})$ is 
a semi-simple symmetric pair, $\mathfrak a$ is 
a maximal split abelian subspace of $\mathfrak p^{\bf c}$ and 
$\mathfrak p^{{\bf c}}=\mathfrak a^{{\bf c}}+
\sum\limits_{\alpha\in\triangle_+}\mathfrak p_{\alpha}^{\bf c}$ is 
the root space decomposition with respect to $\mathfrak a$.  
Here we note that ${\mathfrak a}^{\bf c}$ is the centralizer of $\mathfrak a$ 
in $\mathfrak p^{\bf c}$ and $\mathfrak p^{\bf c}_{\alpha}=
\{X\in \mathfrak p^{\bf c}\,\vert\,({\rm ad}(a)^{\bf c})^2(X)
=\alpha(a)^2X\,\,{\rm for}\,\,{\rm all}\,\,a\in\mathfrak a\}$.  
See [R] and [OS] about general theory of a semi-simple symmetric pair.  
Let $G^{\bf c}$ (resp. $K^{\bf c}$) be the complexification of $G$ 
(resp. $K$).  The 2-multiple of the real part 
${\rm Re}\langle\,\,,\,\,\rangle^{\bf c}$ of 
$\langle\,\,,\,\,\rangle^{\bf c}$ is the Killing form 
of $\mathfrak g^{\bf c}$ regarded as a real Lie algebra.  The restriction 
$2{\rm Re}\langle\,\,,\,\,\rangle^{\bf c}\vert_{{\mathfrak p}^{\bf c}\times
{\mathfrak p}^{\bf c}}$ is an ${\rm Ad}(K^{\bf c})$-invariant non-degenerate 
inner product of ${\mathfrak p}^{\bf c}$ ($=T_{eK^{\bf c}}
(G^{\bf c}/K^{\bf c})$).  
Denote by $\langle\,\,,\,\,\rangle'$ the $G^{\bf c}$-invariant 
pseudo-Riemannian metric on $G^{\bf c}/K^{\bf c}$ 
induced from $2{\rm Re}\langle\,\,,\,\,\rangle^{\bf c}
\vert_{{\mathfrak p}^{\bf c}\times{\mathfrak p}^{\bf c}}$.  
Define an almost complex structure $J_0$ of ${\mathfrak p}^{\bf c}$ 
by $J_0X=\sqrt{-1}X$ ($X\in{\mathfrak p}^{\bf c}$).  It is clear that 
$J_0$ is ${\rm Ad}(K^{\bf c})$-invariant.  Denote by 
$\widetilde J$ the $G^{\bf c}$-invariant almost complex structure 
on $G^{\bf c}/K^{\bf c}$ induced from $J_0$.  It is shown 
that $(G^{\bf c}/K^{\bf c},\langle\,\,,\,\,\rangle',\widetilde J)$ is 
an anti-Kaehlerian manifold and a (semi-simple) pseudo-Riemannian symmetric 
space.  We call this anti-Kaehlerian 
manifold an {\it anti-Kaehlerian symmetric space associated with} $G/K$ 
and simply denote it by $G^{\bf c}/K^{\bf c}$.  
Next we shall recall the notion of an anti-Kaehlerian equifocal submanifold 
which was introduced in [K2].  Let $f$ be an isometric immersion of 
an anti-Kaehlerian manifold $(M,\langle\,\,,\,\,\rangle,J)$ into 
$G^{\bf c}/K^{\bf c}$.  
If $\widetilde J\circ f_{\ast}=f_{\ast}\circ J$, then 
$M$ is called an {\it anti-Kaehlerian submanifold} immersed by $f$.  
If, for each $x\in M$, $\exp^{\perp}(T^{\perp}_xM)$ is totally geodesic, then 
$M$ is called a {\it a submanifold with section}.  
Denote by $\exp^{\perp}$ the normal exponential map of $M$.  
Let $v\in T^{\perp}_xM$.  
If $\exp^{\perp}(av_x+bJv_x)$ is a focal point of $(M,x)$, then we call 
the complex number $a+b\sqrt{-1}$ a {\it complex focal radius along the 
geodesic} $\gamma_{v_x}$.  
Assume that the normal bundle of $M$ is abelian and globally flat and that, 
for each unit normal vector field $v$, the number (which may be $\infty$) of 
distinct complex focal radii along the geodesic $\gamma_{v_x}$ 
is independent of the choice of $x\in M$.  
Then we can define the complex radius functions as above.  
If, for parallel unit normal vector field $v$, 
the number of distinct complex 
focal radii along $\gamma_{v_x}$ is independent of the choice of $x\in M$, 
complex focal radius functions for $v$ are constant on $M$ and they have 
constant multiplicity, 
then $M$ is called an {\it anti-Kaehlerian equifocal submanifold}.  
Let $\phi^{\bf c}:H^0([0,1],\mathfrak g^{\bf c})\to G^{\bf c}$ be the parallel 
transport map for $G^{\bf c}$.  See Section 6 of [K2] about the definition of 
the parallel transport map.  
This map $\phi^{\bf c}$ is an anti-Kaehlerian submersion.  
Let $\pi^{\bf c}:G^{\bf c}\to G^{\bf c}/K^{\bf c}$ be the natural projection.  
It is shown that 
$M$ is anti-Kaehlerian equifocal if and only if each component of 
$(\pi^{\bf c}\circ\phi^{\bf c})^{-1}(M)$ is anti-Kaehlerian isoparametric.  
See Section 5 of [K2] about the definition of an anti-Kaehlerian isoparametric 
submanifold.  In particular, if each component of 
$(\pi^{\bf c}\circ\phi^{\bf c})^{-1}(M)$ is 
proper anti-Kaehlerian isoparametric (i.e., anti-Kaehlerian isoparametric and, 
for each unit normal vector $v$, the shape operator $A_v$ is 
diagonalizable with respect to an orthonormal base of the tangent space 
regarded as a complex vector space), then 
we call $M$ a {\it proper anti-Kaehlerian equifocal submanifold}.  
Assume that $M$ is an anti-Kaehlerian equifocal.  
Let $r$ be a complex focal radius for a parallel unit normal vector field 
$v$.  Then $rv$ is called a {\it focal normal vector field} of $M$.  Then 
a focal map $f_{rv}:M\to G^{\bf c}/K^{\bf c}$ is defined by $f_{rv}(x)=
\exp^{\perp}(rv_x)$ ($x\in M$).  Set $F_{rv}:=f_{rv}(M)$.  We call $F_{rv}$ 
the {\it focal submanifold of} $M$ {\it for} $rv$.  Define a distribution 
$E_{rv}$ on $M$ by $(E_{rv})_x:={\rm Ker}(f_{rv})_{\ast x}$ ($x\in M$). 
We call $E_{rv}$ the {\it focal distribution on} $M$ {\it for} $rv$.  
It is clear that $E_{rv}$ is integrable.  It is shown that the focal set of 
$M$ at $x$ consists of the images by $\exp^{\perp}$ of infinitely many complex 
hyperplanes (which are called complex focal hyperplanes) 
in $T^{\perp}_xM$ (see [K2]).  Denote by $S$ the set of all complex focal 
hyperplanes of $M$ at $x$.  If $\sharp\{{\it l}\in S\,\vert\,rv_x\in{\it l}\}
=1$, then the leaves of $E_{rv}$ are the images by the normal exponential map 
of complex spheres in normal spaces of $F_{rv}$, where $\sharp(\cdot)$ is 
the cardinal number of $(\cdot)$.  
Let $r_1$ (resp. $r_2$) be a complex focal radius for a parallel unit normal 
vector field $v_1$ (resp. $v_2$).  
If $\{{\it l}\in S\,\vert\,r_1(v_1)_x\in{\it l}\}
=\{{\it l}\in S\,\vert\,r_2(v_2)_x\in{\it l}\}$, then we have 
$E_{r_1v_1}=E_{r_2v_2}$.  

\newpage


\centerline{
\unitlength 0.1in
\begin{picture}( 62.9000, 36.9400)(  1.6000,-43.4400)
%
\special{pn 8}%
\special{pa 1810 1150}%
\special{pa 3410 1150}%
\special{pa 3410 2550}%
\special{pa 1810 2550}%
\special{pa 1810 1150}%
\special{fp}%
%
\special{pn 8}%
\special{ar 2610 1750 400 200  6.2831853 6.2831853}%
\special{ar 2610 1750 400 200  0.0000000 3.1415927}%
%
\special{pn 8}%
\special{ar 2610 1750 400 200  3.1415927 3.1815927}%
\special{ar 2610 1750 400 200  3.3015927 3.3415927}%
\special{ar 2610 1750 400 200  3.4615927 3.5015927}%
\special{ar 2610 1750 400 200  3.6215927 3.6615927}%
\special{ar 2610 1750 400 200  3.7815927 3.8215927}%
\special{ar 2610 1750 400 200  3.9415927 3.9815927}%
\special{ar 2610 1750 400 200  4.1015927 4.1415927}%
\special{ar 2610 1750 400 200  4.2615927 4.3015927}%
\special{ar 2610 1750 400 200  4.4215927 4.4615927}%
\special{ar 2610 1750 400 200  4.5815927 4.6215927}%
\special{ar 2610 1750 400 200  4.7415927 4.7815927}%
\special{ar 2610 1750 400 200  4.9015927 4.9415927}%
\special{ar 2610 1750 400 200  5.0615927 5.1015927}%
\special{ar 2610 1750 400 200  5.2215927 5.2615927}%
\special{ar 2610 1750 400 200  5.3815927 5.4215927}%
\special{ar 2610 1750 400 200  5.5415927 5.5815927}%
\special{ar 2610 1750 400 200  5.7015927 5.7415927}%
\special{ar 2610 1750 400 200  5.8615927 5.9015927}%
\special{ar 2610 1750 400 200  6.0215927 6.0615927}%
\special{ar 2610 1750 400 200  6.1815927 6.2215927}%
%
\special{pn 8}%
\special{pa 2610 1150}%
\special{pa 2610 1880}%
\special{fp}%
\special{pa 2610 2030}%
\special{pa 2610 2550}%
\special{fp}%
%
\special{pn 8}%
\special{pa 1810 1950}%
\special{pa 3130 2550}%
\special{fp}%
%
\special{pn 8}%
\special{pa 2210 1750}%
\special{pa 2610 1750}%
\special{fp}%
\special{sh 1}%
\special{pa 2610 1750}%
\special{pa 2544 1730}%
\special{pa 2558 1750}%
\special{pa 2544 1770}%
\special{pa 2610 1750}%
\special{fp}%
%
\special{pn 20}%
\special{sh 1}%
\special{ar 2210 1750 10 10 0  6.28318530717959E+0000}%
\special{sh 1}%
\special{ar 2210 1750 10 10 0  6.28318530717959E+0000}%
%
\special{pn 8}%
\special{pa 2130 2700}%
\special{pa 2310 2170}%
\special{dt 0.045}%
\special{sh 1}%
\special{pa 2310 2170}%
\special{pa 2270 2228}%
\special{pa 2294 2222}%
\special{pa 2308 2240}%
\special{pa 2310 2170}%
\special{fp}%
\special{pa 3700 2690}%
\special{pa 3700 2690}%
\special{dt 0.045}%
%
\special{pn 8}%
\special{pa 2310 1510}%
\special{pa 2430 1750}%
\special{dt 0.045}%
\special{sh 1}%
\special{pa 2430 1750}%
\special{pa 2418 1682}%
\special{pa 2406 1702}%
\special{pa 2382 1700}%
\special{pa 2430 1750}%
\special{fp}%
%
\special{pn 8}%
\special{pa 2370 1010}%
\special{pa 2610 1290}%
\special{dt 0.045}%
\special{sh 1}%
\special{pa 2610 1290}%
\special{pa 2582 1226}%
\special{pa 2576 1250}%
\special{pa 2552 1252}%
\special{pa 2610 1290}%
\special{fp}%
%
\special{pn 8}%
\special{pa 2880 3030}%
\special{pa 2730 1940}%
\special{dt 0.045}%
\special{sh 1}%
\special{pa 2730 1940}%
\special{pa 2720 2010}%
\special{pa 2738 1994}%
\special{pa 2760 2004}%
\special{pa 2730 1940}%
\special{fp}%
%
\special{pn 8}%
\special{pa 3220 920}%
\special{pa 3120 1270}%
\special{dt 0.045}%
\special{sh 1}%
\special{pa 3120 1270}%
\special{pa 3158 1212}%
\special{pa 3136 1220}%
\special{pa 3120 1200}%
\special{pa 3120 1270}%
\special{fp}%
\put(22.0000,-16.6000){\makebox(0,0)[rt]{$x$}}%
\put(23.7000,-13.0000){\makebox(0,0)[rt]{$rv_x$}}%
\put(25.6000,-7.9000){\makebox(0,0)[rt]{$\exp^{\perp}({\it l}_1)$}}%
\put(35.6000,-6.5000){\makebox(0,0)[rt]{$\exp^{\perp}(T^{\perp}_xM)$}}%
\put(24.1000,-27.5000){\makebox(0,0)[rt]{$\exp^{\perp}({\it l}_2)$}}%
\put(35.3000,-30.8000){\makebox(0,0)[rt]{the leaf of $E_{rv}$}}%
\put(33.5000,-32.7000){\makebox(0,0)[rt]{through $x$}}%
%
\special{pn 8}%
\special{ar 2380 3880 290 550  2.1717133 4.0596559}%
%
\special{pn 8}%
\special{ar 3290 3890 290 550  5.3651221 6.2831853}%
\special{ar 3290 3890 290 550  0.0000000 0.9698793}%
%
\special{pn 8}%
\special{ar 2620 3950 260 480  5.3877216 6.2831853}%
\special{ar 2620 3950 260 480  0.0000000 0.8701619}%
%
\special{pn 8}%
\special{ar 3520 3940 260 480  2.2714308 4.0370563}%
%
\special{pn 8}%
\special{ar 3080 3950 190 90  0.0499584 3.1415927}%
%
\special{pn 8}%
\special{ar 3080 3950 200 100  3.1415927 3.2215927}%
\special{ar 3080 3950 200 100  3.4615927 3.5415927}%
\special{ar 3080 3950 200 100  3.7815927 3.8615927}%
\special{ar 3080 3950 200 100  4.1015927 4.1815927}%
\special{ar 3080 3950 200 100  4.4215927 4.5015927}%
\special{ar 3080 3950 200 100  4.7415927 4.8215927}%
\special{ar 3080 3950 200 100  5.0615927 5.1415927}%
\special{ar 3080 3950 200 100  5.3815927 5.4615927}%
\special{ar 3080 3950 200 100  5.7015927 5.7815927}%
\special{ar 3080 3950 200 100  6.0215927 6.1015927}%
\put(26.9000,-37.7000){\makebox(0,0)[rt]{in fact}}%
%
\special{pn 8}%
\special{pa 5462 2300}%
\special{pa 5468 2330}%
\special{pa 5464 2362}%
\special{pa 5454 2394}%
\special{pa 5438 2420}%
\special{pa 5420 2446}%
\special{pa 5398 2470}%
\special{pa 5376 2492}%
\special{pa 5352 2514}%
\special{pa 5324 2532}%
\special{pa 5298 2550}%
\special{pa 5270 2566}%
\special{pa 5242 2580}%
\special{pa 5214 2594}%
\special{pa 5184 2606}%
\special{pa 5154 2618}%
\special{pa 5124 2628}%
\special{pa 5092 2636}%
\special{pa 5062 2644}%
\special{pa 5030 2650}%
\special{pa 4998 2654}%
\special{pa 4966 2658}%
\special{pa 4934 2660}%
\special{pa 4902 2660}%
\special{pa 4872 2656}%
\special{pa 4840 2652}%
\special{pa 4808 2644}%
\special{pa 4778 2632}%
\special{pa 4750 2618}%
\special{pa 4728 2596}%
\special{pa 4710 2570}%
\special{pa 4706 2560}%
\special{sp}%
%
\special{pn 8}%
\special{pa 4706 2560}%
\special{pa 4700 2530}%
\special{pa 4704 2498}%
\special{pa 4714 2468}%
\special{pa 4730 2440}%
\special{pa 4748 2414}%
\special{pa 4770 2390}%
\special{pa 4794 2368}%
\special{pa 4818 2346}%
\special{pa 4844 2328}%
\special{pa 4870 2312}%
\special{pa 4898 2294}%
\special{pa 4926 2280}%
\special{pa 4956 2266}%
\special{pa 4986 2254}%
\special{pa 5014 2242}%
\special{pa 5046 2232}%
\special{pa 5076 2224}%
\special{pa 5108 2216}%
\special{pa 5138 2210}%
\special{pa 5170 2206}%
\special{pa 5202 2202}%
\special{pa 5234 2200}%
\special{pa 5266 2200}%
\special{pa 5298 2204}%
\special{pa 5330 2208}%
\special{pa 5360 2216}%
\special{pa 5390 2228}%
\special{pa 5418 2242}%
\special{pa 5442 2264}%
\special{pa 5460 2290}%
\special{pa 5462 2300}%
\special{sp -0.045}%
%
\special{pn 8}%
\special{ar 4820 1550 400 200  0.1721908 0.2001811}%
%
\special{pn 8}%
\special{pa 5138 1522}%
\special{pa 5150 1552}%
\special{pa 5152 1584}%
\special{pa 5146 1614}%
\special{pa 5136 1646}%
\special{pa 5122 1674}%
\special{pa 5106 1702}%
\special{pa 5086 1726}%
\special{pa 5064 1750}%
\special{pa 5042 1774}%
\special{pa 5018 1796}%
\special{pa 4994 1816}%
\special{pa 4970 1836}%
\special{pa 4944 1854}%
\special{pa 4916 1872}%
\special{pa 4888 1886}%
\special{pa 4860 1902}%
\special{pa 4832 1916}%
\special{pa 4802 1928}%
\special{pa 4772 1940}%
\special{pa 4742 1950}%
\special{pa 4710 1960}%
\special{pa 4680 1966}%
\special{pa 4648 1970}%
\special{pa 4616 1974}%
\special{pa 4584 1974}%
\special{pa 4552 1972}%
\special{pa 4520 1966}%
\special{pa 4490 1954}%
\special{pa 4464 1938}%
\special{pa 4442 1914}%
\special{pa 4428 1884}%
\special{pa 4424 1852}%
\special{pa 4428 1822}%
\special{pa 4438 1790}%
\special{pa 4450 1762}%
\special{pa 4468 1734}%
\special{pa 4486 1708}%
\special{pa 4506 1684}%
\special{pa 4528 1660}%
\special{pa 4552 1638}%
\special{pa 4576 1618}%
\special{pa 4602 1598}%
\special{pa 4626 1578}%
\special{pa 4654 1560}%
\special{pa 4682 1546}%
\special{pa 4710 1530}%
\special{pa 4738 1516}%
\special{pa 4768 1502}%
\special{pa 4798 1490}%
\special{pa 4828 1480}%
\special{pa 4858 1472}%
\special{pa 4890 1464}%
\special{pa 4920 1458}%
\special{pa 4952 1456}%
\special{pa 4984 1454}%
\special{pa 5016 1454}%
\special{pa 5048 1460}%
\special{pa 5078 1470}%
\special{pa 5106 1486}%
\special{pa 5130 1508}%
\special{pa 5138 1522}%
\special{sp}%
%
\special{pn 8}%
\special{pa 5618 3178}%
\special{pa 5616 3210}%
\special{pa 5604 3238}%
\special{pa 5586 3266}%
\special{pa 5564 3288}%
\special{pa 5538 3308}%
\special{pa 5512 3326}%
\special{pa 5482 3340}%
\special{pa 5454 3354}%
\special{pa 5424 3364}%
\special{pa 5394 3374}%
\special{pa 5362 3382}%
\special{pa 5330 3388}%
\special{pa 5300 3394}%
\special{pa 5268 3398}%
\special{pa 5236 3402}%
\special{pa 5204 3402}%
\special{pa 5172 3402}%
\special{pa 5140 3400}%
\special{pa 5108 3398}%
\special{pa 5076 3396}%
\special{pa 5044 3390}%
\special{pa 5014 3382}%
\special{pa 4982 3374}%
\special{pa 4952 3362}%
\special{pa 4924 3350}%
\special{pa 4896 3334}%
\special{pa 4870 3314}%
\special{pa 4848 3292}%
\special{pa 4830 3264}%
\special{pa 4822 3234}%
\special{pa 4820 3226}%
\special{sp}%
%
\special{pn 8}%
\special{pa 4820 3226}%
\special{pa 4824 3194}%
\special{pa 4836 3164}%
\special{pa 4854 3138}%
\special{pa 4876 3116}%
\special{pa 4902 3096}%
\special{pa 4928 3078}%
\special{pa 4956 3062}%
\special{pa 4986 3050}%
\special{pa 5016 3038}%
\special{pa 5046 3030}%
\special{pa 5076 3022}%
\special{pa 5108 3014}%
\special{pa 5140 3008}%
\special{pa 5172 3004}%
\special{pa 5204 3002}%
\special{pa 5236 3002}%
\special{pa 5268 3000}%
\special{pa 5300 3002}%
\special{pa 5332 3004}%
\special{pa 5364 3008}%
\special{pa 5394 3014}%
\special{pa 5426 3020}%
\special{pa 5456 3030}%
\special{pa 5486 3040}%
\special{pa 5516 3054}%
\special{pa 5544 3070}%
\special{pa 5568 3088}%
\special{pa 5592 3112}%
\special{pa 5608 3138}%
\special{pa 5618 3168}%
\special{pa 5618 3178}%
\special{sp -0.045}%
%
\special{pn 8}%
\special{ar 3920 3740 930 2210  5.0709435 6.1828940}%
%
\special{pn 8}%
\special{ar 4350 3940 1320 3060  5.1600065 6.1237801}%
%
\special{pn 8}%
\special{ar 4170 3530 1060 2290  5.0852484 6.2791495}%
%
\special{pn 20}%
\special{sh 1}%
\special{ar 5100 2440 10 10 0  6.28318530717959E+0000}%
\special{sh 1}%
\special{ar 4810 1710 10 10 0  6.28318530717959E+0000}%
\special{sh 1}%
\special{ar 5220 3210 10 10 0  6.28318530717959E+0000}%
\special{sh 1}%
\special{ar 5220 3210 10 10 0  6.28318530717959E+0000}%
%
\special{pn 8}%
\special{pa 4920 2660}%
\special{pa 5090 2440}%
\special{fp}%
\special{sh 1}%
\special{pa 5090 2440}%
\special{pa 5034 2482}%
\special{pa 5058 2482}%
\special{pa 5066 2506}%
\special{pa 5090 2440}%
\special{fp}%
%
\special{pn 8}%
\special{pa 5460 2300}%
\special{pa 5110 2440}%
\special{fp}%
\special{sh 1}%
\special{pa 5110 2440}%
\special{pa 5180 2434}%
\special{pa 5160 2420}%
\special{pa 5164 2398}%
\special{pa 5110 2440}%
\special{fp}%
%
\special{pn 8}%
\special{pa 4920 2300}%
\special{pa 5100 2440}%
\special{fp}%
\special{sh 1}%
\special{pa 5100 2440}%
\special{pa 5060 2384}%
\special{pa 5058 2408}%
\special{pa 5036 2416}%
\special{pa 5100 2440}%
\special{fp}%
%
\special{pn 8}%
\special{pa 4640 1970}%
\special{pa 4810 1720}%
\special{fp}%
\special{sh 1}%
\special{pa 4810 1720}%
\special{pa 4756 1764}%
\special{pa 4780 1764}%
\special{pa 4790 1786}%
\special{pa 4810 1720}%
\special{fp}%
%
\special{pn 8}%
\special{pa 5150 1520}%
\special{pa 4800 1720}%
\special{fp}%
\special{sh 1}%
\special{pa 4800 1720}%
\special{pa 4868 1704}%
\special{pa 4846 1694}%
\special{pa 4848 1670}%
\special{pa 4800 1720}%
\special{fp}%
%
\special{pn 8}%
\special{pa 4600 1600}%
\special{pa 4810 1710}%
\special{fp}%
\special{sh 1}%
\special{pa 4810 1710}%
\special{pa 4760 1662}%
\special{pa 4764 1686}%
\special{pa 4742 1698}%
\special{pa 4810 1710}%
\special{fp}%
%
\special{pn 8}%
\special{pa 5060 3390}%
\special{pa 5220 3210}%
\special{fp}%
\special{sh 1}%
\special{pa 5220 3210}%
\special{pa 5162 3248}%
\special{pa 5186 3250}%
\special{pa 5192 3274}%
\special{pa 5220 3210}%
\special{fp}%
%
\special{pn 8}%
\special{pa 5620 3170}%
\special{pa 5220 3210}%
\special{fp}%
\special{sh 1}%
\special{pa 5220 3210}%
\special{pa 5288 3224}%
\special{pa 5274 3206}%
\special{pa 5284 3184}%
\special{pa 5220 3210}%
\special{fp}%
%
\special{pn 8}%
\special{pa 5040 3030}%
\special{pa 5210 3200}%
\special{fp}%
\special{sh 1}%
\special{pa 5210 3200}%
\special{pa 5178 3140}%
\special{pa 5172 3162}%
\special{pa 5150 3168}%
\special{pa 5210 3200}%
\special{fp}%
%
\special{pn 8}%
\special{pa 5090 2430}%
\special{pa 4610 3070}%
\special{fp}%
%
\special{pn 8}%
\special{pa 5090 2430}%
\special{pa 5350 2080}%
\special{dt 0.045}%
%
\special{pn 8}%
\special{pa 5410 1990}%
\special{pa 5670 1600}%
\special{fp}%
%
\special{pn 20}%
\special{sh 1}%
\special{ar 4920 2660 10 10 0  6.28318530717959E+0000}%
\special{sh 1}%
\special{ar 4920 2660 10 10 0  6.28318530717959E+0000}%
%
\special{pn 8}%
\special{pa 4370 2470}%
\special{pa 4700 2160}%
\special{dt 0.045}%
\special{sh 1}%
\special{pa 4700 2160}%
\special{pa 4638 2192}%
\special{pa 4662 2198}%
\special{pa 4666 2220}%
\special{pa 4700 2160}%
\special{fp}%
%
\special{pn 8}%
\special{pa 4700 1080}%
\special{pa 4620 1480}%
\special{dt 0.045}%
\special{sh 1}%
\special{pa 4620 1480}%
\special{pa 4654 1420}%
\special{pa 4630 1428}%
\special{pa 4614 1412}%
\special{pa 4620 1480}%
\special{fp}%
\special{pa 5970 1910}%
\special{pa 5590 1710}%
\special{dt 0.045}%
\special{sh 1}%
\special{pa 5590 1710}%
\special{pa 5640 1760}%
\special{pa 5638 1736}%
\special{pa 5658 1724}%
\special{pa 5590 1710}%
\special{fp}%
\special{pa 6440 3770}%
\special{pa 5470 3360}%
\special{dt 0.045}%
\special{sh 1}%
\special{pa 5470 3360}%
\special{pa 5524 3404}%
\special{pa 5520 3382}%
\special{pa 5540 3368}%
\special{pa 5470 3360}%
\special{fp}%
%
\special{pn 8}%
\special{pa 6450 3770}%
\special{pa 5330 2530}%
\special{dt 0.045}%
\special{sh 1}%
\special{pa 5330 2530}%
\special{pa 5360 2594}%
\special{pa 5366 2570}%
\special{pa 5390 2566}%
\special{pa 5330 2530}%
\special{fp}%
%
\special{pn 8}%
\special{pa 6450 3770}%
\special{pa 5050 1770}%
\special{dt 0.045}%
\special{sh 1}%
\special{pa 5050 1770}%
\special{pa 5072 1836}%
\special{pa 5082 1814}%
\special{pa 5106 1814}%
\special{pa 5050 1770}%
\special{fp}%
\put(43.6000,-24.0000){\makebox(0,0)[rt]{$M$}}%
\put(50.3000,-26.9000){\makebox(0,0)[rt]{$x$}}%
\put(48.3000,-8.7000){\makebox(0,0)[rt]{$F_{rv}$}}%
\put(68.7000,-19.0000){\makebox(0,0)[rt]{$\exp^{\perp}(T^{\perp}_xM)$}}%
\put(70.5000,-38.1000){\makebox(0,0)[rt]{leaves of $E_{rv}$}}%
\end{picture}%
\hspace{6truecm}
}

\vspace{0.1truecm}

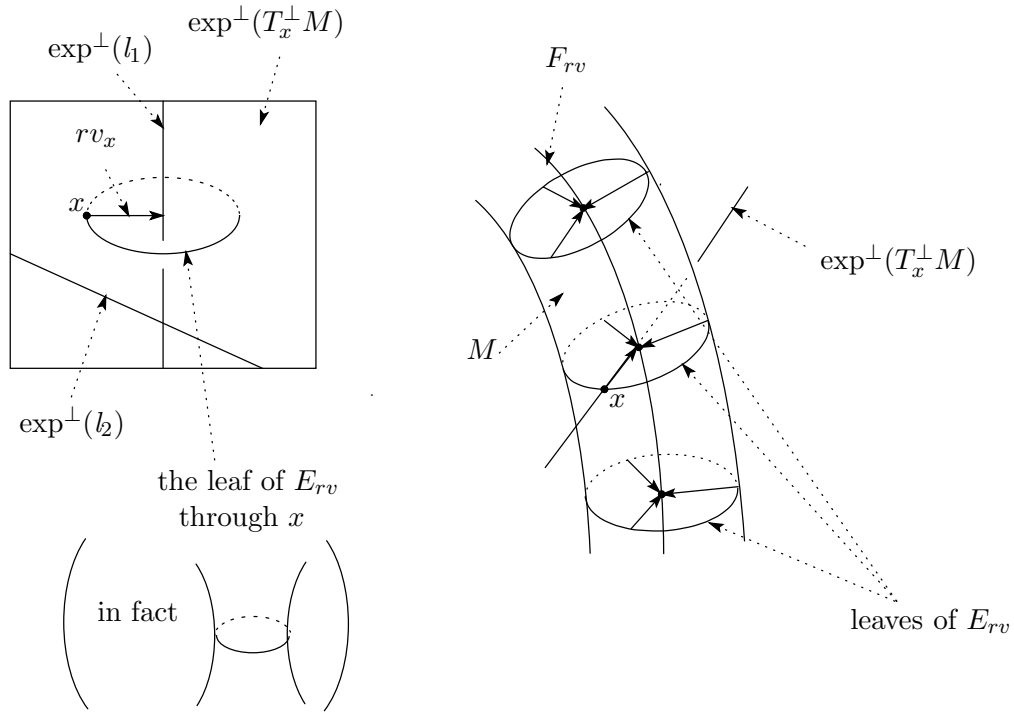
\centerline{{\bf Fig. 3.}}



\vspace{0.2truecm}

\centerline{
\unitlength 0.1in
\begin{picture}( 61.1000, 37.6100)(  7.3000,-44.6100)
%
\special{pn 8}%
\special{ar 3800 1800 400 200  6.2831853 6.2831853}%
\special{ar 3800 1800 400 200  0.0000000 3.1415927}%
%
\special{pn 8}%
\special{ar 3800 1800 400 200  3.1415927 3.1815927}%
\special{ar 3800 1800 400 200  3.3015927 3.3415927}%
\special{ar 3800 1800 400 200  3.4615927 3.5015927}%
\special{ar 3800 1800 400 200  3.6215927 3.6615927}%
\special{ar 3800 1800 400 200  3.7815927 3.8215927}%
\special{ar 3800 1800 400 200  3.9415927 3.9815927}%
\special{ar 3800 1800 400 200  4.1015927 4.1415927}%
\special{ar 3800 1800 400 200  4.2615927 4.3015927}%
\special{ar 3800 1800 400 200  4.4215927 4.4615927}%
\special{ar 3800 1800 400 200  4.5815927 4.6215927}%
\special{ar 3800 1800 400 200  4.7415927 4.7815927}%
\special{ar 3800 1800 400 200  4.9015927 4.9415927}%
\special{ar 3800 1800 400 200  5.0615927 5.1015927}%
\special{ar 3800 1800 400 200  5.2215927 5.2615927}%
\special{ar 3800 1800 400 200  5.3815927 5.4215927}%
\special{ar 3800 1800 400 200  5.5415927 5.5815927}%
\special{ar 3800 1800 400 200  5.7015927 5.7415927}%
\special{ar 3800 1800 400 200  5.8615927 5.9015927}%
\special{ar 3800 1800 400 200  6.0215927 6.0615927}%
\special{ar 3800 1800 400 200  6.1815927 6.2215927}%
%
\special{pn 8}%
\special{pa 3800 1200}%
\special{pa 3800 1930}%
\special{fp}%
\special{pa 3800 2080}%
\special{pa 3800 2600}%
\special{fp}%
%
\special{pn 8}%
\special{pa 3400 1800}%
\special{pa 3800 1800}%
\special{fp}%
\special{sh 1}%
\special{pa 3800 1800}%
\special{pa 3734 1780}%
\special{pa 3748 1800}%
\special{pa 3734 1820}%
\special{pa 3800 1800}%
\special{fp}%
%
\special{pn 20}%
\special{sh 1}%
\special{ar 3400 1800 10 10 0  6.28318530717959E+0000}%
\special{sh 1}%
\special{ar 3400 1800 10 10 0  6.28318530717959E+0000}%
%
\special{pn 8}%
\special{pa 3560 1060}%
\special{pa 3800 1340}%
\special{dt 0.045}%
\special{sh 1}%
\special{pa 3800 1340}%
\special{pa 3772 1276}%
\special{pa 3766 1300}%
\special{pa 3742 1302}%
\special{pa 3800 1340}%
\special{fp}%
%
\special{pn 8}%
\special{pa 4410 970}%
\special{pa 4310 1320}%
\special{dt 0.045}%
\special{sh 1}%
\special{pa 4310 1320}%
\special{pa 4348 1262}%
\special{pa 4326 1270}%
\special{pa 4310 1250}%
\special{pa 4310 1320}%
\special{fp}%
\put(33.9000,-17.1000){\makebox(0,0)[rt]{$x$}}%
\put(37.5000,-8.4000){\makebox(0,0)[rt]{$\exp^{\perp}({\it l}_1)$}}%
\put(47.5000,-7.0000){\makebox(0,0)[rt]{$\exp^{\perp}(T^{\perp}_xM)$}}%
\put(29.8000,-27.5000){\makebox(0,0)[rt]{$\exp^{\perp}({\it l}_2)$}}%
\put(52.0000,-32.0000){\makebox(0,0)[rt]{the leaf of $E_{r_1v_1}=E_{r_2v_2}$}}%
\put(45.5000,-33.9000){\makebox(0,0)[rt]{through $x$}}%
%
\special{pn 8}%
\special{ar 3580 4000 290 550  2.1717133 4.0596559}%
%
\special{pn 8}%
\special{ar 4260 4008 290 550  5.3651221 6.2831853}%
\special{ar 4260 4008 290 550  0.0000000 0.9698793}%
%
\special{pn 8}%
\special{ar 3590 4068 260 480  5.3877216 6.2831853}%
\special{ar 3590 4068 260 480  0.0000000 0.8701619}%
%
\special{pn 8}%
\special{ar 4490 4058 260 480  2.2714308 4.0370563}%
%
\special{pn 8}%
\special{ar 4050 4068 190 90  0.0499584 3.1415927}%
%
\special{pn 8}%
\special{ar 4050 4068 200 100  3.1415927 3.2215927}%
\special{ar 4050 4068 200 100  3.4615927 3.5415927}%
\special{ar 4050 4068 200 100  3.7815927 3.8615927}%
\special{ar 4050 4068 200 100  4.1015927 4.1815927}%
\special{ar 4050 4068 200 100  4.4215927 4.5015927}%
\special{ar 4050 4068 200 100  4.7415927 4.8215927}%
\special{ar 4050 4068 200 100  5.0615927 5.1415927}%
\special{ar 4050 4068 200 100  5.3815927 5.4615927}%
\special{ar 4050 4068 200 100  5.7015927 5.7815927}%
\special{ar 4050 4068 200 100  6.0215927 6.1015927}%
\put(37.8000,-38.8700){\makebox(0,0)[rt]{in fact}}%
%
\special{pn 8}%
\special{pa 2400 1200}%
\special{pa 4600 1200}%
\special{pa 4600 2600}%
\special{pa 2400 2600}%
\special{pa 2400 1200}%
\special{fp}%
%
\special{pn 8}%
\special{pa 2400 1880}%
\special{pa 4600 2540}%
\special{fp}%
%
\special{pn 8}%
\special{pa 3400 1820}%
\special{pa 3800 2300}%
\special{fp}%
\special{sh 1}%
\special{pa 3800 2300}%
\special{pa 3774 2236}%
\special{pa 3766 2260}%
\special{pa 3742 2262}%
\special{pa 3800 2300}%
\special{fp}%
%
\special{pn 8}%
\special{pa 3390 1800}%
\special{pa 3280 2140}%
\special{fp}%
\special{sh 1}%
\special{pa 3280 2140}%
\special{pa 3320 2084}%
\special{pa 3296 2090}%
\special{pa 3282 2070}%
\special{pa 3280 2140}%
\special{fp}%
%
\special{pn 8}%
\special{pa 3400 1800}%
\special{pa 3800 1400}%
\special{fp}%
\special{sh 1}%
\special{pa 3800 1400}%
\special{pa 3740 1434}%
\special{pa 3762 1438}%
\special{pa 3768 1462}%
\special{pa 3800 1400}%
\special{fp}%
%
\special{pn 8}%
\special{pa 4200 1800}%
\special{pa 3800 1800}%
\special{fp}%
\special{sh 1}%
\special{pa 3800 1800}%
\special{pa 3868 1820}%
\special{pa 3854 1800}%
\special{pa 3868 1780}%
\special{pa 3800 1800}%
\special{fp}%
%
\special{pn 8}%
\special{pa 4200 1800}%
\special{pa 3800 1400}%
\special{fp}%
\special{sh 1}%
\special{pa 3800 1400}%
\special{pa 3834 1462}%
\special{pa 3838 1438}%
\special{pa 3862 1434}%
\special{pa 3800 1400}%
\special{fp}%
%
\special{pn 8}%
\special{pa 3400 1800}%
\special{pa 3426 1820}%
\special{pa 3442 1848}%
\special{pa 3450 1878}%
\special{pa 3454 1910}%
\special{pa 3456 1942}%
\special{pa 3456 1974}%
\special{pa 3454 2006}%
\special{pa 3448 2038}%
\special{pa 3444 2068}%
\special{pa 3438 2100}%
\special{pa 3430 2132}%
\special{pa 3420 2162}%
\special{pa 3412 2192}%
\special{pa 3400 2222}%
\special{pa 3388 2252}%
\special{pa 3376 2282}%
\special{pa 3360 2310}%
\special{pa 3344 2338}%
\special{pa 3328 2364}%
\special{pa 3306 2388}%
\special{pa 3284 2412}%
\special{pa 3258 2430}%
\special{pa 3228 2436}%
\special{pa 3218 2436}%
\special{sp}%
%
\special{pn 8}%
\special{pa 3200 2430}%
\special{pa 3174 2412}%
\special{pa 3160 2384}%
\special{pa 3152 2354}%
\special{pa 3148 2322}%
\special{pa 3146 2290}%
\special{pa 3148 2258}%
\special{pa 3150 2226}%
\special{pa 3156 2194}%
\special{pa 3162 2162}%
\special{pa 3168 2132}%
\special{pa 3176 2100}%
\special{pa 3186 2070}%
\special{pa 3196 2040}%
\special{pa 3208 2010}%
\special{pa 3220 1980}%
\special{pa 3234 1952}%
\special{pa 3250 1924}%
\special{pa 3264 1896}%
\special{pa 3284 1870}%
\special{pa 3304 1846}%
\special{pa 3328 1822}%
\special{pa 3354 1804}%
\special{pa 3384 1798}%
\special{pa 3394 1800}%
\special{sp -0.045}%
%
\special{pn 8}%
\special{pa 2470 2710}%
\special{pa 2780 1990}%
\special{dt 0.045}%
\special{sh 1}%
\special{pa 2780 1990}%
\special{pa 2736 2044}%
\special{pa 2760 2040}%
\special{pa 2772 2060}%
\special{pa 2780 1990}%
\special{fp}%
%
\special{pn 20}%
\special{sh 1}%
\special{ar 3220 2440 10 10 0  6.28318530717959E+0000}%
\special{sh 1}%
\special{ar 4200 1800 10 10 0  6.28318530717959E+0000}%
\special{sh 1}%
\special{ar 4200 1800 10 10 0  6.28318530717959E+0000}%
%
\special{pn 8}%
\special{pa 3210 2440}%
\special{pa 3290 2140}%
\special{fp}%
\special{sh 1}%
\special{pa 3290 2140}%
\special{pa 3254 2200}%
\special{pa 3276 2192}%
\special{pa 3292 2210}%
\special{pa 3290 2140}%
\special{fp}%
%
\special{pn 8}%
\special{pa 3400 1800}%
\special{pa 3000 2050}%
\special{fp}%
\special{sh 1}%
\special{pa 3000 2050}%
\special{pa 3068 2032}%
\special{pa 3046 2022}%
\special{pa 3046 1998}%
\special{pa 3000 2050}%
\special{fp}%
%
\special{pn 8}%
\special{pa 3210 2440}%
\special{pa 3000 2060}%
\special{fp}%
\special{sh 1}%
\special{pa 3000 2060}%
\special{pa 3016 2128}%
\special{pa 3026 2108}%
\special{pa 3050 2110}%
\special{pa 3000 2060}%
\special{fp}%
%
\special{pn 20}%
\special{sh 1}%
\special{ar 3000 2060 10 10 0  6.28318530717959E+0000}%
\special{sh 1}%
\special{ar 3000 2060 10 10 0  6.28318530717959E+0000}%
%
\special{pn 8}%
\special{pa 3330 1620}%
\special{pa 3610 1800}%
\special{dt 0.045}%
\special{sh 1}%
\special{pa 3610 1800}%
\special{pa 3566 1748}%
\special{pa 3566 1772}%
\special{pa 3544 1782}%
\special{pa 3610 1800}%
\special{fp}%
%
\special{pn 8}%
\special{pa 3370 1390}%
\special{pa 3620 1570}%
\special{dt 0.045}%
\special{sh 1}%
\special{pa 3620 1570}%
\special{pa 3578 1516}%
\special{pa 3578 1540}%
\special{pa 3554 1548}%
\special{pa 3620 1570}%
\special{fp}%
%
\special{pn 8}%
\special{pa 2290 1630}%
\special{pa 3130 1960}%
\special{dt 0.045}%
\special{sh 1}%
\special{pa 3130 1960}%
\special{pa 3076 1918}%
\special{pa 3080 1940}%
\special{pa 3062 1954}%
\special{pa 3130 1960}%
\special{fp}%
%
\special{pn 8}%
\special{pa 2210 1400}%
\special{pa 3320 1990}%
\special{dt 0.045}%
\special{sh 1}%
\special{pa 3320 1990}%
\special{pa 3272 1942}%
\special{pa 3274 1966}%
\special{pa 3252 1976}%
\special{pa 3320 1990}%
\special{fp}%
\put(33.3000,-14.5000){\makebox(0,0)[rt]{$r_1(v_1)_x$}}%
\put(33.7000,-12.2000){\makebox(0,0)[rt]{$r_2(v_2)_x$}}%
\put(21.8000,-12.4000){\makebox(0,0)[rt]{$r_3(v_3)_x$}}%
\put(22.8000,-15.0000){\makebox(0,0)[rt]{$r_4(v_4)_x$}}%
%
\special{pn 8}%
\special{pa 4300 2720}%
\special{pa 3390 2270}%
\special{dt 0.045}%
\special{sh 1}%
\special{pa 3390 2270}%
\special{pa 3442 2318}%
\special{pa 3438 2294}%
\special{pa 3460 2282}%
\special{pa 3390 2270}%
\special{fp}%
\put(55.9000,-27.3000){\makebox(0,0)[rt]{the leaf of $E_{r_3v_3}=E_{r_4v_4}$}}%
\put(49.5000,-29.2000){\makebox(0,0)[rt]{through $x$}}%
%
\special{pn 8}%
\special{pa 3870 3160}%
\special{pa 3930 1990}%
\special{dt 0.045}%
\special{sh 1}%
\special{pa 3930 1990}%
\special{pa 3908 2056}%
\special{pa 3928 2044}%
\special{pa 3948 2058}%
\special{pa 3930 1990}%
\special{fp}%
%
\special{pn 8}%
\special{ar 5870 3138 290 550  2.1717133 4.0596559}%
%
\special{pn 8}%
\special{ar 6550 3158 290 550  5.3651221 6.2831853}%
\special{ar 6550 3158 290 550  0.0000000 0.9698793}%
%
\special{pn 8}%
\special{ar 5880 3218 260 480  5.3877216 6.2831853}%
\special{ar 5880 3218 260 480  0.0000000 0.8701619}%
%
\special{pn 8}%
\special{ar 6780 3208 260 480  2.2714308 4.0370563}%
%
\special{pn 8}%
\special{ar 6340 3218 190 90  0.0499584 3.1415927}%
%
\special{pn 8}%
\special{ar 6340 3218 200 100  3.1415927 3.2215927}%
\special{ar 6340 3218 200 100  3.4615927 3.5415927}%
\special{ar 6340 3218 200 100  3.7815927 3.8615927}%
\special{ar 6340 3218 200 100  4.1015927 4.1815927}%
\special{ar 6340 3218 200 100  4.4215927 4.5015927}%
\special{ar 6340 3218 200 100  4.7415927 4.8215927}%
\special{ar 6340 3218 200 100  5.0615927 5.1415927}%
\special{ar 6340 3218 200 100  5.3815927 5.4615927}%
\special{ar 6340 3218 200 100  5.7015927 5.7815927}%
\special{ar 6340 3218 200 100  6.0215927 6.1015927}%
\put(60.7000,-30.3700){\makebox(0,0)[rt]{in fact}}%
\put(67.4000,-38.9000){\makebox(0,0)[rt]{$E_{r_5v_5}=E_{r_1v_1}\oplus E_{r_3v_3}$}}%
%
\special{pn 8}%
\special{pa 4750 1840}%
\special{pa 3680 2150}%
\special{dt 0.045}%
\special{sh 1}%
\special{pa 3680 2150}%
\special{pa 3750 2152}%
\special{pa 3732 2136}%
\special{pa 3738 2112}%
\special{pa 3680 2150}%
\special{fp}%
\put(53.4000,-16.7000){\makebox(0,0)[rt]{$r_5(v_5)_x$}}%
\end{picture}%
\hspace{3truecm}}

\vspace{0.1truecm}

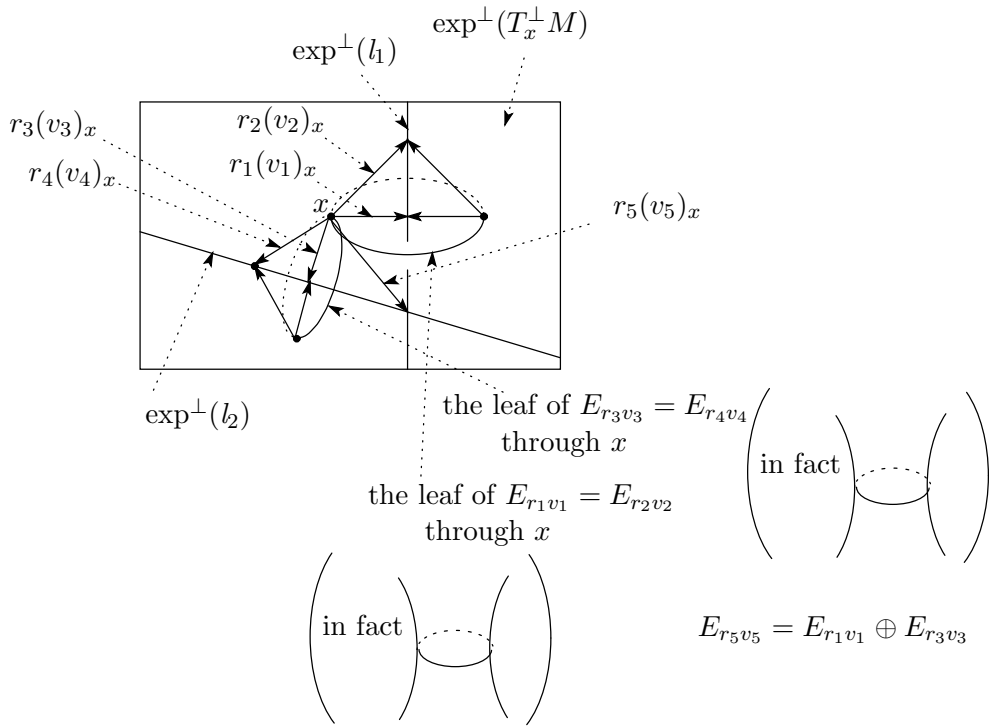
\centerline{{\bf Fig. 4.}}

\newpage


\noindent
Next we recall the notion of the extrinsic complexification of a complete 
$C^{\omega}$-submanifold in a symmetric space of non-compact type which 
was introduced in [K2].  First we recall the complexification 
of a complete $C^{\omega}$-Riemannian manifold.  
Let $M$ be a complete $C^{\omega}$-Riemannian manifold.  
The notion of the adapted complex structure on a neighborhood $U$ of the 
$0$-section of the tangent bundle $TM$ is defined as the complex structure 
(on $U$) such that, for each geodesic $\gamma:{\bf R}\to N$, 
the restriction of its differential $\gamma_{\ast}:T{\bf R}=
{\bf C}\to TM$ to $\gamma_{\ast}^{-1}(U)$ is holomorphic.  
We take $U$ as largely as possible under the condition that $U\cap T_xM$ 
is a star-shaped neighborhood of $0_x$ for each $x\in M$, where $0_x$ is 
the zero vector of $T_xM$.  If $N$ is of non-negative curvature, then 
we have $U=TM$.  Also, if all sectional curvatures of $M$ are bigger than 
or equal to $c$ ($c<0$), 
then $U$ contains the ball bundle $T^rM:=\{X\in TM\,\vert\,\vert\vert X
\vert\vert<r\}$ of radius $r:=\frac{\pi}{2\sqrt{-c}}$.  
In detail, see [Sz1$\sim$4].  
Denote by $J_A$ the adapted complex structure on $U$.  
The complex manifold $(U,J_A)$ is interpreted as the complexification of 
$N$.  We denote $(U,J_A)$ by $M^{\bf c}$ and call it 
the complexification of $M$, where we note that $M^{{\bf c}}$ 
is given no Riemannian metric.  
In particular, in case of $M={\bf R}^m$ (the Euclidean space), we 
have $(U,J_A)={\bf C}^m$.  
Also, in the case where $N$ is a symmetric space $G/K$ of 
non-compact type, there exists the holomorphic diffeomorphism $\delta$ of 
$(U,J_A)$ onto an open subset of $G^{{\bf c}}/K^{{\bf c}}$.  
Let $M$ be an immersed (complete) $C^{\omega}$-submanifold in $G/K$.  
Denote by $f$ its immersion.  Let $M^{{\bf c}}$ be the 
complexification of $M$ (defined as above).  
We shall define the 
complexification $f^{{\bf c}}:M^{{\bf c}}\rightarrow 
G^{{\bf c}}/K^{{\bf c}}$ of $f$, where we shrink 
$M^{{\bf c}}$ to a neighborhood of the $0$-section of $TM$ if 
necessary.  For its purpose, we first 
define the complexification of a $C^{\omega}$-curve 
$\alpha:{\bf R}\to G/K$.  Let $\mathfrak g=\mathfrak f+\mathfrak p$ be the 
Cartan decomposition associated with $G/K$ 
and $W:{\bf R}\to\mathfrak p$ be the curve in 
$\mathfrak p$ with $(\exp W(t))K=\alpha(t)$ ($t\in{\bf R}$), 
where we note that $W$ is uniquely determined because $G/K$ is of 
non-compact type.  
Since $\alpha$ is of class $C^{\omega}$, so is also $W$.  
Let $W^{{\bf c}}:D\to\mathfrak p^{{\bf c}}$ 
($D\,:\,$ a neighborhood of ${\bf R}$ in ${\bf C}$) be 
the holomorphic extension of $W$.  
We define the complexification 
$\alpha^{{\bf c}}:D\to G^{{\bf c}}/K^{{\bf c}}$ 
of $\alpha$ by 
$\alpha^{{\bf c}}(z)=(\exp W^{{\bf c}}(z))K^{{\bf c}}$.  
It is shown that this complexification of a $C^{\omega}$-curve in $G/K$ is 
a holomorphic curve in $G^{\bf c}/K^{\bf c}$.  
By using this complexification of a $C^{\omega}$-curve in $G/K$, 
we define the complexification 
$f^{{\bf c}}:M^{{\bf c}}\rightarrow G^{{\bf c}}/K^{{\bf c}}$ 
of $f$ by $f^{{\bf c}}(X):=(f\circ\gamma^M_X)^{{\bf c}}(\sqrt{-1})$ 
($X\in M^{{\bf c}}\,(\subset TM)$), 
where $\gamma^M_X$ is the geodesic in $M$ with $\dot{\gamma}^M_X(0)=X$.  
Here we shrink $M^{{\bf c}}$ to a neighborhood of the $0$-section 
of $TM$ if necessary in order to assure that $\sqrt{-1}$ belongs to 
the domain of $(f\circ\gamma^M_X)^{{\bf c}}$ 
for each $X\in M^{{\bf c}}$.  It is shown that the map 
$f^{{\bf c}}:M^{{\bf c}}\rightarrow G^{{\bf c}}/K^{{\bf c}}$ is holomorphic 
and that the restriction of $f^{{\bf c}}$ to a neighborhood $U'$ 
of the $0$-section of $TM$ is an immersion, where we take $U'$ 
as largely as possible.  Denote by $M^{{\bf c}}$ 
this neighborhood $U'$ newly.  Give $M^{{\bf c}}$ the Riemannian 
metric induced from that of $G^{{\bf c}}/K^{{\bf c}}$ by $f^{{\bf c}}$.  
Then $M^{{\bf c}}$ is an anti-Kaehlerian submanifold in $G^{{\bf c}}/
K^{{\bf c}}$ immersed by $f^{{\bf c}}$.  
We call this anti-Kaehlerian submanifold $M^{{\bf c}}$ immersed by 
$f^{{\bf c}}$ the {\it extrinsic complexification} of 
the submanifold $M$.  We consider the case 
where $M$ is (extrinsically) homogeneous.  
Concretely we consider the case where 
$M=H(g_0K)$ and $f$ is the inclusion map of $M$ into $G/K$, where $H$ is 
a closed subgroup of $G$.  
Let $\iota$ be a natural immersion of $G/K$ into $G^{{\bf c}}/
K^{{\bf c}}$, that is, $\iota(gK)=gK^{{\bf c}}$ 
$(X\in \mathfrak g)$.  It is shown that $\iota$ is totally geodesic.  
Let $\mathfrak g_H^{{\bf c}}$ be the complexification of the Lie 
algebra of $H$ and set $H^{{\bf c}}:=
\exp\mathfrak g_H^{{\bf c}}$.  
For a homogeneous submanifold $M=H(g_0K)$, 
the image $f^{{\bf c}}(M^{{\bf c}})$ is an open subset 
of the orbit $H^{{\bf c}}(g_0K^{{\bf c}})$.  
Hence this orbit is the complete extension of $M^{\bf c}$.  
It is shown that $M$ is complex equifocal if and only if $M^{\bf c}$ is 
anti-Kaehlerian equifocal (see Theorem 5 of [K2]).  Also, 
it is shown that $M$ is proper complex equifocal if and only if 
$M^{\bf c}$ is proper anti-Kaehlerian equifocal.  

\section{Aks-representations} 
In this section, we shall first introduce the notions of an anti-Kaehlerian 
symmetric pair and an anti-Kaehlerian symmetric Lie algebra, and investigate 
the correspondence relations of those notions with an anti-Kaehlerian 
symmetric space.  
Let $(M,J,\langle\,\,,\,\,\rangle)$ be an anti-Kaehlerian manifold (i.e., 
$J^2=-{\rm id}, \nabla J=0$ ($\nabla\,:\,$ the Levi-Civita connection of 
$\langle\,\,,\,\,\rangle$)) and $\langle JX,JY\rangle=-\langle X,Y\rangle$ 
($X,Y\in TM$)).  
In the sequel, denote by the same symbol ${\rm id}$ 
the identity transformations of various sets.  
If there exists an involutive holomorphic isometry $s_p$ of 
$M$ having $p$ as an isolated fixed point for each $p\in M$, then we call 
$(M,J,\langle\,\,,\,\,\rangle)$ an {\it anti-Kaehlerian symmetric space}.  
Also, if there exists a local involutive holomorphic isometry defined on 
a neighborhood of $p$ having $p$ as an isolated fixed point for each $p\in M$, 
then we call $(M,J,\langle\,\,,\,\,\rangle)$ a {\it locally anti-Kaehlerian 
symmetric space}.  In this section, we introduce the notions of an 
anti-Kaehlerian symmetric pair and an anti-Kaehlerian symmetric Lie algebra 
in relation with an anti-Kaehlerian symmetric spaces.  Let $G$ be a connected 
complex Lie group and $K$ be a closed complex subgroup of $G$.  If there 
exists an involutive (complex) automorphism $\rho$ of $G$ such that 
$G_{\rho}^0\subset K\subset G_{\rho}$ ($G_{\rho}\,:\,$ the group of all fixed 
points of $\rho$, $G_{\rho}^0\,:\,$ the identity component of $G_{\rho}$)
then we call the pair $(G,K)$ an {\it anti-Kaehlerian symmetric pair}.  
If $\mathfrak g$ be a complex Lie algebra and $\tau$ be a complex involution 
of $\mathfrak g$, then we call such a pair $(\mathfrak g,\tau)$ an 
{\it anti-Kaehlerian symmertric Lie algebra}.  
Let $\mathfrak f:={\rm Ker}(\tau-{\rm id})$ and $\mathfrak p:={\rm Ker}
(\tau+{\rm id})$.  Denote by ${\rm Ad}_G$ and ${\rm ad}_{\mathfrak g}$ 
the adjoint representations of $G$ and $\mathfrak g$, repsectively.  
Also, denote by $j$ the complex structure of $\mathfrak g$.  
Let ${\mathfrak p}_{\bf R}$ be the totally real subspace of $\mathfrak p$ such 
that $\langle\,\,,\,\,\rangle
\vert_{\mathfrak p_{\bf R}\times j\mathfrak p_{\bf R}}=0$ and that 
$\langle\,\,,\,\,\rangle\vert_{\mathfrak p_{\bf R}\times\mathfrak p_{\bf R}}$ 
is positive definite.  Here we note that such a totally real subspace is 
determined uniquely.  Set 
${\rm ad}_{\mathfrak g}\vert_{\mathfrak p}(\mathfrak f)
:=\{{\rm ad}_{\mathfrak g}(X)\vert_{\mathfrak p}\,\vert\,X\in \mathfrak f\}$, 
${\rm Ad}_G\vert_{\mathfrak p}(K):=\{{\rm Ad}_G(k)\vert_{\mathfrak p}\,
\vert\,k\in K\}$, 
${\rm ad}_{\mathfrak g}\vert_{\mathfrak p_{\bf R}}(\mathfrak f):=\{
{\rm pr}_{\mathfrak p_{\bf R}}\circ{\rm ad}_{\mathfrak g}(X)
\vert_{\mathfrak p_{\bf R}}\,\vert\,X\in\mathfrak f\}$ and 
${\rm Ad}_G\vert_{\mathfrak p_{\bf R}}(K):=\exp_{\rm GL(\mathfrak p_{\bf R})}
({\rm ad}_{\mathfrak g}\vert_{\mathfrak p_{\bf R}}(\mathfrak f))$, where 
$\exp_{{\rm GL}(\mathfrak p_{\bf R})}$ is the exponential map of 
${\rm GL}(\mathfrak p_{\bf R})$.  
Let $SO_{AK}(\mathfrak p)$ be the identity component of the group 
$\{A\in {\rm GL}(\mathfrak p)\,\vert\,A^{\ast}\langle\,\,,\,\,\rangle=
\langle\,\,,\,\,\rangle,\,\,A\circ j=j\circ A\}$ and set 
$\mathfrak{so}_{AK}(\mathfrak p):=\{A\in\mathfrak{gl}(\mathfrak p)\,\vert\,
A\circ j=j\circ A,\,\,\langle AX,Y\rangle=-\langle X,AY\rangle\,
(\forall\,X,Y\in\mathfrak p)\}$, which is the Lie algebra of 
$SO_{AK}(\mathfrak p)$.  Then we have the following fact.  

\vspace{0.5truecm}

\noindent
{\bf Lemma 3.1.} {\sl The complexification 
$\mathfrak{so}(\mathfrak p_{\bf R})^{\bf c}$ of 
$\mathfrak{so}(\mathfrak p_{\bf R})$ coincides with $\mathfrak{so}_{AK}
(\mathfrak p)$ and hence $SO(\mathfrak p_{\bf R})$ is a half-dimensional 
totally real compact subgroup of $SO_{AK}(\mathfrak p)$.  Also, the 
complexification $({\rm ad}_{\mathfrak g}\vert_{\mathfrak p_{\bf R}}
(\mathfrak f))^{\bf c}$ of 
${\rm ad}_{\mathfrak g}\vert_{\mathfrak p_{\bf R}}(\mathfrak f)$ coincides 
with ${\rm ad}_{\mathfrak g}\vert_{\mathfrak p}(\mathfrak f)$ and 
${\rm ad}_{\mathfrak g}\vert_{\mathfrak p_{\bf R}}(\mathfrak f)$ is contained 
in $\mathfrak{so}(\mathfrak p_{\bf R})$.  Hence 
${\rm Ad}_G\vert_{\mathfrak p_{\bf R}}(K)$ is a half-dimensional totally real 
compact subgroup of ${\rm Ad}_G\vert_{\mathfrak p}(K)$ contained in 
$SO(\mathfrak p_{\bf R})$.}

\vspace{0.5truecm}

\noindent
{\it Proof.} For $A\in\mathfrak{gl}(\mathfrak p_{\bf R})$, denote by 
$\widetilde A$ the element of $\mathfrak{gl}(\mathfrak p,j)
:=\{B\in\mathfrak{gl}(\mathfrak p)\,\vert\,B\circ j=j\circ B\}$ whose 
restriction to $\mathfrak p_{\bf R}$ is equal to $A$.  Let 
$C\in\mathfrak{so}_{AK}(\mathfrak p)$.  Set 
$A:={\rm pr}_{\mathfrak p_{\bf R}}\circ C\vert_{\mathfrak p_{\bf R}}$ and 
$B:=-j\circ{\rm pr}_{j\mathfrak p_{\bf R}}\circ C\vert_{\mathfrak p_{\bf R}}$.  Then we have $C=\widetilde A+j\widetilde B$.  Take 
$X,Y\in\mathfrak p_{\bf R}$.  Then it follows from 
$\langle\mathfrak p_{\bf R},j\mathfrak p_{\bf R}\rangle=0$ that 
$\langle CX,jY\rangle=-\langle BX,Y\rangle$ and $\langle X,C(jY)\rangle
=-\langle BY,X\rangle$.  
Hence it follows from $\langle CX,jY\rangle=-\langle X,C(jY)\rangle$ that 
$\langle BX,Y\rangle=-\langle X,BY\rangle$.  Thus we have 
$B\in\mathfrak{so}(\mathfrak p_{\bf R})$.  Also we have 
$\langle CX,Y\rangle=\langle AX,Y\rangle$ and 
$\langle X,CY\rangle=-\langle X,AY\rangle$.  
Hence we have $\langle AX,Y\rangle=-\langle X,AY\rangle$.  
Thus we have $A\in\mathfrak{so}(\mathfrak p_{\bf R})$.  
Therefore we have $C\in\mathfrak{so}(\mathfrak p_{\bf R})^{\bf c}$.  
Thus we have $\mathfrak{so}_{AK}(\mathfrak p)\subset\mathfrak{so}
(\mathfrak p_{\bf R})^{\bf c}$.  Since $\mathfrak{so}_{AK}(\mathfrak p)$ and 
$\mathfrak{so}(\mathfrak p_{\bf R})^{\bf c}$ are of the same dimension, 
we have $\mathfrak{so}_{AK}(\mathfrak p)=\mathfrak{so}
(\mathfrak p_{\bf R})^{\bf c}$.  Therefore the first-half statement of this 
lemma is shown.  Let $C\in{\rm ad}\vert_{\mathfrak p}(\mathfrak f)$.  
Set $A:={\rm pr}_{\mathfrak p_{\bf R}}\circ C\vert_{\mathfrak p_{\bf R}}$ 
and $B:=-j\circ{\rm pr}_{j\mathfrak p_{\bf R}}\circ C
\vert_{\mathfrak p_{\bf R}}$.  From the definition of 
${\rm ad}_{\mathfrak p_{\bf R}}(\mathfrak f)$, we have 
$A\in{\rm ad}\vert_{\mathfrak p_{\bf R}}(\mathfrak f)$.  
Also, it follows from $-j\circ C\in{\rm ad}\vert_{\mathfrak p}
(\mathfrak f)$ that $-({\rm pr}_{\mathfrak p_{\bf R}}\circ j\circ C)
\vert_{\mathfrak p_{\bf R}}\in{\rm ad}\vert_{\mathfrak p_{\bf R}}
(\mathfrak f)$.  Clearly we have 
$-({\rm pr}_{\mathfrak p_{\bf R}}\circ j\circ C)
\vert_{\mathfrak p_{\bf R}}=B$.  Thus we have 
$B\in\mathfrak{ad}\vert_{\mathfrak p_{\bf R}}(\mathfrak f)$.  
Therefore we have $C(=\widetilde A+j\circ\widetilde B)\in
({\rm ad}\vert_{\mathfrak p_{\bf R}}(\mathfrak f))^{\bf c}$.  
Thus ${\rm ad}\vert_{\mathfrak p}(\mathfrak f)\subset
({\rm ad}\vert_{\mathfrak p_{\bf R}}(\mathfrak f))^{\bf c}$ is obtained.  
From ${\rm dim}_{\bf R}{\rm ad}\vert_{\mathfrak p}(\mathfrak f)
={\rm dim}_{\bf R}({\rm ad}\vert_{\mathfrak p_{\bf R}}(\mathfrak f))^{\bf c}$, 
it follows that ${\rm ad}\vert_{\mathfrak p}(\mathfrak f)=
({\rm ad}\vert_{\mathfrak p_{\bf R}}(\mathfrak f))^{\bf c}$.  
Also, since $C\in{\rm ad}\vert_{\mathfrak p}(\mathfrak f)\subset
\mathfrak{so}_{AK}(\mathfrak p)$, we can show $A\in\mathfrak{so}
(\mathfrak p_{\bf R})$ as above.  Therefore we obtain 
${\rm ad}\vert_{\mathfrak p_{\bf R}}(\mathfrak f)\subset
\mathfrak{so}(\mathfrak p_{\bf R})$.  
Hence ${\rm Ad}\vert_{\mathfrak p_{\bf R}}(K)\subset SO(\mathfrak p_{\bf R})$ 
is obtained.  Furthermore, since ${\rm Ad}\vert_{\mathfrak p_{\bf R}}(K)$ is 
closed in $SO(\mathfrak p_{\bf R})$, it is compact.  Thus the second-half 
statement of this lemma follows.
\begin{flushright}q.e.d.\end{flushright}

\vspace{0.5truecm}

Now we show that an anti-Kaehlerian symmetric pair arises from an anti-Kaehlerian symmetric space.  

\vspace{0.5truecm}

\noindent
{\bf Proposition 3.2.} {\sl Let $(M,J,\langle\,\,,\,\,\rangle)$ be an 
anti-Kaehlerian symmetric space, $G$ be the identity component of the isometry 
group of $(M,J,\langle\,\,,\,\,\rangle)$ and $K$ be the isotropy group of $G$ 
at some $p_0\in M$.  Then the pair $(G,K)$ is an anti-Kaehlerian symmetric 
pair.}

\vspace{0.5truecm}

\noindent
{\it Proof.} Identify $M$ with $G/K$ under the correspondence 
$g(p_0)\leftrightarrow gK$ ($g\in G$).  Define a map $\rho:G\to G$ by 
$\rho(g)=s_{p_0}\circ g\circ s_{p_0}$ ($g\in G$), which is an involutive 
automorphism of $G$.  
Easily we can show that $G_{\rho}^0\subset K\subset G_{\rho}$ 
(see the proof of (ii) of Theorem 3.3 of Chapter IV in [H]).  
Let $\mathfrak f:={\rm Ker}(\rho_{\ast e}-{\rm id})$ and $\mathfrak p:={\rm Ker}(\rho_{\ast e}+{\rm id})$, 
where $e$ is the identity element of $G$.  
The space $\mathfrak p$ is identified with $T_{p_0}M$.  Define the $\sqrt{-1}$
-multiple in $\mathfrak g$ by $\sqrt{-1}X=J_{p_0}X$ ($X\in \mathfrak p
=T_{p_0}M$) and $[\sqrt{-1}Y,Z]=[Y,J_{p_0}Z]$ ($Y\in\mathfrak f,\,Z\in
\mathfrak p$), where  $[\,\,,\,\,]$ is the Lie bracket product of 
$\mathfrak g$.  
Note that this $\sqrt{-1}$-multiple in $\mathfrak g$ is well-defined because 
$\mathfrak f$ acts on $\mathfrak p$ effectively.  
Since ${\rm ad}(X)\circ J_{p_0}=J_{p_0}\circ{\rm ad}(X)$ on $\mathfrak p$ 
($X\in\mathfrak f$), $[Y,Z]=-R_{p_0}(Y,Z)$ ($Y,Z\in\mathfrak p$) and 
$R_{p_0}(J_{p_0}Y,Z)=J_{p_0}R_{p_0}(Y,Z)$ ($Y,Z\in\mathfrak p$) by 
anti-Kaehlerity of $M$, we see that $(\mathfrak g,[\,\,,\,\,])$ is a 
complex Lie algebra under this $\sqrt{-1}$-multiple.  Also, 
it is easy to show that $\mathfrak f$ is a complex Lie subalgebra and 
$\rho_{\ast e}$ is the complex involution.  Hence $G,\,K$ and $\rho$ are 
regarded as a complex Lie group, a complex Lie subgroup of $G$ and an 
involutive complex automorphism of $G$, respectively.  \hspace{10.7truecm}q.e.d.

\vspace{0.5truecm}

By using Lemma 3.1, we show that an anti-Kaehlerian symmetric space arises 
from an anti-Kaehlerian symmetric pair.  

\vspace{0.5truecm}

\noindent
{\bf Proposition 3.3.} {\sl Let $(G,K)$ be an anti-Kaehlerian symmetric pair.  
Then there exists an anti-Kaehlerian structure $(J,\langle\,\,,\,\,\rangle)$ 
of $G/K$ such that $(G/K,J,\langle\,\,,\,\,\rangle)$ is an anti-Kaehlerian 
symmetric space.}

\vspace{0.5truecm}

\noindent
{\it Proof.} Since $(G,K)$ is an anti-Kaehlerian symmetric pair, there 
exists an involutive (complex) automorphism $\rho$ of $G$ with 
$G_{\rho}^0\subset K\subset G_{\rho}$.  
Let $\mathfrak g:={\rm Lie}\,G$, $\mathfrak f:={\rm Lie}\,K$ and 
$\mathfrak p:={\rm Ker}(\rho_{\ast e}+{\rm id})$.  
Then we can show ${\rm Ad}_G(K)(\mathfrak p)\subset\mathfrak p$ 
(see the first part of the proof of Proposition 3.4 of Chapter IV in [H]).  
Define an almost complex structure $j$ of $\mathfrak p$ by $j(X)=\sqrt{-1}X$ 
($X\in\mathfrak p$).  It is clear that $j$ is ${\rm Ad}_G(K)$-invariant.  
Denote by $J$ the $G$-invariant almost complex structure on $G/K$ arising 
from $j$.  
Let ${\rm GL}((\mathfrak p,j)):=\{A\in {\rm GL}(\mathfrak p)\,\vert\, A\circ j
=j\circ A\}$, where ${\rm GL}(\mathfrak p)$ is the group of all (real) linear 
isomorphisms of $\mathfrak p$.  
Take a half-dimensional subspace $\mathfrak p_{\bf R}$ of $\mathfrak p$ 
with $\mathfrak p_{\bf R}\oplus j\mathfrak p_{\bf R}=\mathfrak p$.  
The group ${\rm GL}(\mathfrak p_{\bf R})$ of all linear 
isomorphisms of $\mathfrak p_{\bf R}$ is regarded as 
a half-dimensional totally real subgroup of ${\rm GL}((\mathfrak p,j))$ by 
identifying each $A\in {\rm GL}(\mathfrak p_{\bf R})$ with 
$\widetilde A\in {\rm GL}((\mathfrak p,j))$ 
defined by $\widetilde A(X+jY)=AX+jAY$ ($X,Y\in \mathfrak p_{\bf R}$).  
Let ${\rm ad}_{\mathfrak g}\vert_{\mathfrak p_{\bf R}}(\mathfrak f)$ be 
as in the proof of Proposition 3.2 and 
${\rm Ad}_G\vert_{\mathfrak p_{\bf R}}(K):=\exp_{{\rm GL}(\mathfrak p_{\bf R})}
({\rm ad}_{\mathfrak g}\vert_{\mathfrak p_{\bf R}}(\mathfrak f))$.  
It is clear that the group ${\rm Ad}_G\vert_{\mathfrak p_{\bf R}}(K)$ 
is regarded as a half-dimensional totally real subgroup of 
${\rm Ad}_G\vert_{\mathfrak p}(K)$.  
By taking an anti-Kaehlerian inner product $\beta$ of $(\mathfrak p,j)$ such 
that $\beta\vert_{\mathfrak p_{\bf R}\times j\mathfrak p_{\bf R}}=0$ and that 
$\beta\vert_{\mathfrak p_{\bf R}\times\mathfrak p_{\bf R}}$ is positive 
definite hojotekini and using Lemma 3.1, 
${\rm Ad}_G\vert_{\mathfrak p_{\bf R}}(K)$ is a half-dimensional totally real 
compact subgroup of ${\rm Ad}_G\vert_{\mathfrak p}(K)$.  
Define a real bilinear form $\beta_0$ on $\mathfrak p$ by 
$$\beta_0(X,Y)=\int_{a\in{\rm Ad}_G\vert_{\mathfrak p_{\bf R}}(K)}
\beta(aX,aY)\omega\quad(X,Y\in \mathfrak p),$$
where $\omega$ is the Haar measure of ${\rm Ad}_G
\vert_{\mathfrak p_{\bf R}}(K)$ and each $a\in{\rm Ad}_G
\vert_{\mathfrak p_{\bf R}}(K)$ is extended to the linear transformation 
of $\mathfrak p$ in the natural manner.  
We shall show that $\beta_0$ is an anti-Kaehlerian inner product of 
$(\mathfrak p,j)$.  Let $X\in\mathfrak p_{\bf R}$.  Since 
$\beta(aX,aX)\geq0$ for any 
$a\in{\rm Ad}_G\vert_{\mathfrak p_{\bf R}}(K)$, we have $\beta_0(X,X)\geq0$.  
If $\beta_0(X,X)=0$, then we have $\beta(aX,aX)=0$ for any $a\in{\rm Ad}_G
\vert_{\mathfrak p_{\bf R}}(K)$.  In particular, we have $\beta(X,X)=0$, 
that is, $X=0$.  Thus 
$\beta_0\vert_{\mathfrak p_{\bf R}\times\mathfrak p_{\bf R}}$ is positive 
definite.  Let $Y$ be another vector of $\mathfrak p_{\bf R}$.  Since 
$\beta(aX,ajY)=0$ ($a\in{\rm Ad}_G\vert_{\mathfrak p_{\bf R}}(K)$), we have 
$\beta_0(X,jY)=0$.  Thus it follows from the arbitrarinesses of $X$ and 
$Y$ that $\beta_0\vert_{\mathfrak p_{\bf R}\times j\mathfrak p_{\bf R}}=0$.  
On the other hand, it is clear that $\beta_0(jZ,jW)=-\beta_0(Z,W)$ 
($Z,W\in\mathfrak p$).  These facts imply that $\beta_0$ is an anti-Kaehlerian 
inner product of $(\mathfrak p,j)$.  Next we shall show that $\beta_0$ is 
${\rm Ad}_G\vert_{\mathfrak p}(K)$-invariant.  It is clear that $\beta_0$ is 
${\rm Ad}_G\vert_{\mathfrak p_{\bf R}}(K)$-invariant.  Fix 
$X,Y\in\mathfrak p$.  
Define a complex-valued function $f$ on ${\rm Ad}_G\vert_{\mathfrak p}(K)$ 
by $f(a)=\beta_0(aX,aY)-\sqrt{-1}\beta_0(aX,ajY)$ 
($a\in {\rm Ad}_G\vert_{\mathfrak p}(K)$).  
Since $f\equiv\beta_0(X,Y)-\sqrt{-1}\beta_0(X,jY)$ on 
${\rm Ad}_G\vert_{\mathfrak p_{\bf R}}(K)$, 
$f$ is holomorphic and ${\rm Ad}_G\vert_{\mathfrak p_{\bf R}}(K)$ is a 
half-dimensional totally real subgroup of ${\rm Ad}_G\vert_{\mathfrak p}(K)$, 
we see that $f\equiv\beta_0(X,Y)-\sqrt{-1}\beta_0(X,jY)$ on ${\rm Ad}_G
\vert_{\mathfrak p}(K)$, which implies that $\beta_0$ 
is ${\rm Ad}_G\vert_{\mathfrak p}(K)$-invariant.  
Denote by $\langle\,\,,\,\,\rangle$ the 
$G$-invariant pseudo-Riemannian metric on $G/K$ arising from $\beta$.  
It is clear that $(G/K,J,\langle\,\,,\,\,\rangle)$ is an anti-Kaehlerian 
manifold.  Next we shall show that $(G/K,J,\langle\,\,\,\,\,\rangle)$ is 
an anti-Kaehlerian symmetric space.  Let $\pi:G\to G/K$ be the natural 
projection.  Define a map $s_o:G/K\to G/K$ by $s_o(\pi(g))=\pi(\rho(g))$ 
($g\in G$).  It is clear that $s_o$ is well-defined and that 
$s_o^2={\rm id}$.  Also, it is shown that $s_o$ is an isometry of 
$(G/K,\langle\,\,,\,\,\rangle)$ (see the proof of Proposition 3.4 of 
Chapter IV in [H]).  
Furthermore, it is shown that $s_o$ is holomorphic.  Also, we have 
$s_{o\ast\pi(e)}\circ\pi_{\ast e}=\pi_{\ast e}\circ\rho_{\ast e}
=-\pi_{\ast e}$ on $\mathfrak p$, that is, $s_{o\ast\pi(e)}=-{\rm id}$, 
which implies that $\pi(e)$ is an isolated fixed point of $s_o$.  For each 
$g\in G$, define a map $s_{\pi(g)}:G/K\to G/K$ by $s_{\pi(g)}=g\circ s_o
\circ g^{-1}$.  Easily we can show that $s_{\pi(g)}$ is an involutive 
holomorphic isometry of $(G/K,J,\langle\,\,,\,\,\rangle)$ having $\pi(g)$ 
as an isolated fixed point.  Thus $(G/K,J,\langle\,\,,\,\,\rangle)$ is 
an anti-Kaehlerian symmetric space.  \hspace{2.9truecm}q.e.d.

\vspace{0.5truecm}

Let $(\mathfrak g,\tau)$ be an anti-Kaehlerian symmetric Lie algebra and 
$\mathfrak f:={\rm Ker}(\tau-{\rm id})$.  Let $G$ be a connected complex 
Lie group with ${\rm Lie}\,G=\mathfrak g$ and $K$ be a complex Lie subgroup 
of $G$ with ${\rm Lie}\,K=\mathfrak f$.  We call such a pair $(G,K)$ a 
{\it pair associated with} $(\mathfrak g,\tau)$.  

\vspace{0.5truecm}

\noindent
{\bf Proposition 3.4.} {\sl Let $(\mathfrak g,\tau)$ be an anti-Kaehlerian 
symmetric Lie algebra, $(G,K)$ be a pair associated with $(\mathfrak g,\tau)$ 
such that $K$ is connected and $(\widetilde G,\widetilde K)$ be a pair 
associated with $(\mathfrak g,\tau)$ such that $\widetilde G$ is simply 
connected and that $\widetilde K$ is connected.  Then the following 
statements {\rm (i)} and {\rm (ii)} hold:

{\rm (i)} $(\widetilde G,\widetilde K)$ is an anti-Kaehlerian symmetric pair.  

{\rm (ii)} Assume that $K$ is closed.  Let $(J,\langle\,\,,\,\,\rangle)$ 
be a $G$-invariant anti-Kaehlerian structure on $G/K$ defined as in the 
proof of Proposition 3.3.  
Then $(G/K,J,\langle\,\,,\,\,\rangle)$ is a locally 
anti-Kaehlerian symmetric space and the universal anti-Kaehlerian covering 
of $(G/K,J,\langle\,\,,\,\,\rangle)$ is isometric to an anti-Kaehlerian 
symmetric space $\widetilde G/\widetilde K$ equipped with a suitable 
anti-Kaehlerian structure defined as in the 
proof of Proposition 3.3.}

\vspace{0.5truecm}

\noindent
{\it Proof.} First we shall show the statement (i).  Since $\widetilde G$ is 
simply connected, there uniquely exists an involutive (complex) automorphism 
$\rho$ of $\widetilde G$ with $\rho_{\ast e}=\tau$.  
In a standard method, we can show that $\widetilde K$ is equal to the identity 
component $\widetilde G_{\rho}^0$ of the group of all fixed points of 
$\rho$ because $\widetilde K$ is connected.  
Thus $(\widetilde G, \widetilde K)$ is an anti-Kaehlerian symmetric pair.  

Next we shall show the statement (ii).  The groups ${\rm Ad}_G(K)$ and 
${\rm Ad}_{\widetilde G}(\widetilde K)$ coincide with each other because 
they are connected complex Lie subgroups of the adjoint group 
${\rm int}\,\mathfrak g$ and have the same Lie algebra.  
Let $(J,\langle\,\,,\,\,\rangle)$ (resp. $(\widetilde J,\langle\,\,,\,\,
\rangle\,\,\widetilde{}\,\,)$) be a $G$ (resp. $\widetilde G$)-invariant 
anti-Kaehlerian structure on $G/K$ (resp. $\widetilde G/\widetilde K$) 
as in the proof of Proposition 3.3.  
Let $\psi$ be the homomorphism of $\widetilde G$ onto 
$G$ with $\psi_{\ast e}={\rm id}$.  It is clear that 
$\widetilde K$ is the identity component of $\psi^{-1}(K)$.  Hence a map 
$\bar{\psi}:\widetilde G/\widetilde K\to G/K$ is well-defined by 
$\bar{\psi}(\widetilde g\widetilde K)=\psi(\widetilde g)K$ ($\widetilde g
\in\widetilde G$).  It is shown that this map $\bar{\psi}$ is a covering 
map (see Lemma 13.4 of Chapter I in [H]).  It is easy 
to show that $\bar{\psi}$ is an anti-Kaehlerian covering map of 
$(\widetilde G/\widetilde K,\widetilde J,\langle\,\,,\,\,\rangle\,\,
\widetilde{}\,\,)$ onto $(G/K,J,\langle\,\,,\,\,\rangle)$.  Hence 
$(G/K,J,\langle\,\,,\,\,\rangle)$ is a locally anti-Kaehlerian symmetric 
space.  Since $\widetilde G/\widetilde K$ is simply connected (see the 
proof of Proposition 3.6 of Chapter IV in [H]), 
$(\widetilde G/\widetilde K,\widetilde J,\langle\,\,,\,\,\rangle
\,\,\widetilde{}\,\,)$ is the universal anti-Kaehlerian covering of 
$(G/K,J,\langle\,\,,\,\,\rangle)$.  
\hspace{1.4truecm}q.e.d.

\vspace{0.5truecm}

Let $(M,J,\langle\,\,,\,\,\rangle)$ be an irreducible 
anti-Kaehlerian symmetric space, 
$G$ be the identity component of the isometry group of 
$(M,J,\langle\,\,,\,\,\rangle)$ and $K$ be the isotropy group of $G$ at 
some point $p_0\in M$, where the irreducibility implies that $M$ is 
not decomposed into the non-trivial product of two anti-Kaehlerian symmetric 
spaces.  Assume that $(M,J,\langle\,\,,\,\,\rangle)$ 
does not have the pseudo-Euclidean part in its de Rham decomposition.  
Note that an anti-Kaehlerian symmetric space without pseudo-Euclidean part 
is not necessarily semi-simple (see [CP],[W1]).  
Also, let $(\mathfrak g,\tau)$ be the anti-Kaehlerian 
symmetric Lie algebra associated with the anti-Kaehlerian symmetric pair 
$(G,K)$ and $\mathfrak p:={\rm Ker}(\tau+{\rm id})$.  
The space ${\rm Ker}(\tau-{\rm id})$ is equal to the Lie 
algebra $\mathfrak f$ of $K$ and $\mathfrak p$ is identified with 
$T_{p_0}M\,(=T_{eK}(G/K))$.  We call the linear isotropy representation 
${\rm Ad}_G\vert_{\mathfrak p}:K\to{\rm GL}(\mathfrak p)$ an 
{\it aks}-{\it representation}, where $\mathfrak p$ is regarded as an 
anti-Kaehlerian space under the identification $\mathfrak p=T_{p_0}M$.  Let 
$\mathfrak a_s$ be a maximal split abelian subspace of $\mathfrak p$ 
(see [R] or [OS] about the definition of a maximal split abelian subspace) and 
$\mathfrak p=\mathfrak p_0+\sum\limits_{\alpha\in \triangle_+}
\mathfrak p_{\alpha}$ be the root space decomposition with respect to 
$\mathfrak a_s$ (i.e., the simultaneously eigenspace decomposition of 
${\rm ad}(a)^2$'s ($a\in\mathfrak a_s$)), where 
the space $\mathfrak p_{\alpha}$ is defined by 
$\mathfrak p_{\alpha}:=\{X\in\mathfrak p\,\vert\,{\rm ad}(a)^2(X)
=\alpha(a)^2X\,\,{\rm for}\,\,{\rm all}\,\,a\in \mathfrak a_s\}$ and 
$\triangle_+$ is the positive root system with respect to $\mathfrak a_s$ 
under some lexicographic ordering of $\mathfrak a_s^{\ast}$.  
Set $\mathfrak a:=
\mathfrak p_0\,(\supset\mathfrak a_s),\,j:=J_{p_0}$ and 
$\langle\,\,,\,\,\rangle_0:=\langle\,\,,\,\,\rangle_{p_0}$.  It is shown 
that $\langle\,\,,\,\,\rangle_0\vert_{\mathfrak a_s\times\mathfrak a_s}$ is 
positive (or negative) definite, $\mathfrak a=\mathfrak a_s\oplus 
j\mathfrak a_s$ and $\langle\,\,,\,\,\rangle_0\vert_{\mathfrak a_s\times 
j\mathfrak a_s}=0$.  Note that $\mathfrak p_{\alpha}=\{X\in
\mathfrak p\,\vert\,{\rm ad}(a)^2(X)=\alpha^{\bf c}(a)^2X\,\,{\rm for}\,\,
{\rm all}\,\,a\in\mathfrak a\}$ for each $\alpha\in\triangle_+$, where 
$\alpha^{\bf c}$ is the complexification of $\alpha:\mathfrak a_s\to{\bf R}$, 
$\mathfrak a$ is regarded as the complexification 
$\mathfrak a_s^{\bf c}$ of $\mathfrak a_s$ and 
$\alpha^{\bf c}(a)^2X$ means ${\rm Re}(\alpha^{\bf c}(a)^2)X
+{\rm Im}(\alpha^{\bf c}(a)^2)jX$.  Let 
${\it l}_{\alpha}:=(\alpha^{\bf c})^{-1}(0)$ ($\alpha\in\triangle$) and 
$D:=\mathfrak a\setminus\displaystyle{
\mathop{\cup}_{\alpha\in\triangle_+}{\it l}_{\alpha}}$.  Take $u\in D$ and let 
$M$ be the orbit through $u$ of the $K$-action by the linear isotropy 
representation $({\rm Ad}_G\vert_{\mathfrak p})\vert_K$.  
Since $u\in D$, $M$ is a principal orbit.  Denote by $A$ the shape tensor 
of $M$.  Take $v\in T^{\perp}_uM(=\mathfrak a)$.  Then we have 
$T_uM=\sum\limits_{\alpha\in\triangle_+}\mathfrak p_{\alpha}$ and 
$A_v\vert_{\mathfrak p_{\alpha}}=-\frac{\alpha^{\bf c}(v)}{\alpha^{\bf c}(u)}
{\rm id}$ ($\alpha\in\triangle_+$).  
It is easy to show that the $K$-action by 
$({\rm Ad}_G\vert_{\mathfrak p})\vert_K$ is an anti-Kaehlerian 
polar action having $\mathfrak a$ as a section, where an anti-Kaehlerian 
polar action means the finite dimensional version of an anti-Kaehlerian polar 
action on an infinite dimensional anti-Kaehlerian space defined in [K2].  
Furthermore, from $A_v\vert_{\mathfrak p_{\alpha}}=-\frac{\alpha^{\bf c}(v)}
{\alpha^{\bf c}(u)}{\rm id}_{\mathfrak p_{\alpha}}$ and the arbitrarinesses 
of $v$ and $u$, we see that each principal orbit of the $K$-action is proper 
anti-Kaehlerian isoparametric in the sense of [K4].  

In the $2$-dimensional anti-Kaehlerian space $V=({\bf R}^2,J_0,
\langle\,\,,\,\,\rangle_0)$, there uniquely exists a $1$-dimensional 
totally real subspace $W$ of $V$ such that $\langle W,J_0W\rangle_0=0$ and 
that $\langle\,,\,\,\rangle_0\vert_{W\times W}$ is positive definite.  
Let $w\in W\cup J_0W$.  
The quotient manifold $V/{\bf Z}w$ is a flat anti-Kaehlerian manifold 
whose universal anti-Kaehlerian covering is $V$.  We call such 
an anti-Kaehlerian manifold an {\it anti-Kaehlerian cylinder}.  
Let $(G/K,J,\langle\,\,,\,\,\rangle)$ be a semi-simple anti-Kaehlerian 
symmetric space and $\mathfrak a$ be a maximal abelian subspace of 
$\mathfrak p=T_{eK}(G/K)$.  
It is easy to show that $\exp\,\mathfrak a$ is a flat 
totally geodesic submanifold in $G/K$ and that it is 
holomorphic and isometric to the product of some anti-Kaehlerian cylinders.  
We call $\exp\,\mathfrak a$ a {\it maximal anti-Kaehlerian cylindrical 
product}.  Here we note that, if $(M,J,\langle\,\,,\,\,\rangle)$ is not 
semi-simple, then $\exp\,\mathfrak a$ is 
holomorphic and isometric to the product of some anti-Kaehlerian 
cylinders and an anti-Kaehlerian space.  

At the end of this section, 
we shall recall the notion of the anti-Kaehlerian symmetric space associated 
with a Riemannian symmetric space  of non-compact type which was introduced 
in [K2].  Let $G/K$ be a Riemannian symmetric space  of non-compact type 
and $\rho$ be the Cartan involution, where $G$ 
is assumed to be a connected semi-simple Lie 
group admitting a faithful real representation and $K$ can be assumed to be 
a maximal compact subgroup of $G$.  Let $\mathfrak g:={\rm Lie}\,G,\,\,
\mathfrak f:={\rm Lie}\,K$ and $\mathfrak p:={\rm Ker}(\rho_{\ast e}+
{\rm id})$, where $\mathfrak p$ is identified with $T_{eK}(G/K)$.  
Also, let $\mathfrak g^{\bf c}$ (resp. $\rho_{\ast e}^{\bf c}$) be the 
complexification of $\mathfrak g$ (resp. $\rho_{\ast e}$).  
Since $G$ admits a faithful real representation, 
we can define the complexification $G^{\bf c}$ (resp. $K^{\bf c}$) of $G$ 
(resp. $K$) and the compact dual $G^{\ast}(\subset G^{\bf c})$ of $G$.  
It is shown that $(G^{\bf c},K^{\bf c})$ is an anti-Kaehlerian symmetric 
pair.  Let $\beta$ be the ${\rm Ad}_G(K)$-invariant (positive definite) inner 
product of $\mathfrak p$ arising the Riemannian metric of $G/K$.  
Let $\langle\,\,,\,\,\rangle$ be the pseudo-Riemannian 
metric of $G^{\bf c}/K^{\bf c}$ arising from 
${\rm Re}\,\beta^{\bf c}\,(\mathfrak p^{\bf c}\times\mathfrak p^{\bf c}
\to{\bf R})$ and $J$ be the natural almost complex structure of 
$G^{\bf c}/K^{\bf c}$, where $\mathfrak p^{\bf c}$ is identified with 
$T_{eK^{\bf c}}(G^{\bf c}/K^{\bf c})$.  Then 
$(G^{\bf c}/K^{\bf c},J,\langle\,\,,\,\,\rangle)$ is an 
anti-Kaehlerian symmetric space.  
We call this anti-Kaehlerian symmetric space 
the {\it anti-Kaehlerian symmetric space associated with} $G/K$, where we 
note that $G^{\bf c}/K^{\bf c}$ is a semi-simple 
anti-Kaehlerian symmetric space.  

\vspace{0.5truecm}

\noindent
{\it Remark 3.1.}  If $\beta$ is the Killing-Cartan form of $\mathfrak g$, 
$2{\rm Re}\,\beta^{\bf c}$ is that of $\mathfrak g^{\bf c}$ regarded as 
a real Lie algebra.  

\section{Anti-Kaehlerian holonomy systems} 
Let $(V,R,G)$ be a triple consisting of a Euclidean space $V$, a 
curvature-like tensor $R\,(\in V^{\ast}\otimes V^{\ast}\otimes V^{\ast}
\otimes V)$ and a compact connected Lie subgroup $G$ of the linear isometry 
group $O(V)$ of $V$.  J. Simons [Si] called $(V,R,G)$ a {\it holonomy system} 
if $R(v_1,v_2)\in{\rm Lie}\,G$ for all $v_1,v_2\in V$.  In this section, we 
introduce the notion of an anti-Kaehlerian holonomy system and show some 
facts for the system.  
Let $(V,J,\langle\,\,,\,\,\rangle)$ be a (finite dimensional) 
anti-Kaehlerian space and $R\,(\in V^{\ast}\otimes V^{\ast}\otimes 
V^{\ast}\otimes V)$ be a curvature-like tensor.  Let $SO_{AK}(V)$ be 
the identity component of the group 
$\{A\in GL(V)\,\vert\,A^{\ast}\langle\,\,,\,\,\rangle=\langle\,\,,\,\,
\rangle,\,\,[A,J]=0\}$ and $G$ be a connected 
complex Lie subgroup of $SO_{AK}(V)$.  We call the triple 
$((V,J,\langle\,\,,\,\,\rangle),R,G)$ an {\it anti-Kaehlerian holonomy 
system} if the following two conditions hold: 

\vspace{0.2truecm}

(AH-i) $J\circ R(v_1,v_2)=R(Jv_1,v_2)=R(v_1,v_2)\circ J$ for all $v_1,v_2
\in V$,

(AH-ii) $R(v_1,v_2)\in{\rm Lie}\,G$ for all $v_1,v_2\in V$.  

\vspace{0.2truecm}

\noindent
Furthermore, if the following condition (S) holds, then we call that the 
triple is {\it symmetric}:

\vspace{0.2truecm}

(S) $R(gv_1,gv_2)gv_3=gR(v_1,v_2)v_3$ for all $v_i\in V$ ($i=1,2,3$) and 
all $g\in G$.  

\vspace{0.2truecm}

\noindent
Also, if $G$ is weakly irreducible, then we call that the triple is 
{\it weakly irreducible}, where the weakly irreduciblity of $G$ implies that 
there exists no $G$-invariant non-degenerate subspace $W$ of $V$ with 
$W\not=\{0\}$ and $W\not=V$.  
Here we give examples of an anti-Kaehlerian holonomy system.  

\vspace{0.5truecm}

\noindent
{\it Example 1.}  Let $(M,J,\langle\,\,,\,\,\rangle)$ be an anti-Kaehlerian 
manifold.  Let $\nabla$ be the Levi-Civita connection of 
$\langle\,\,,\,\,\rangle$, $R$ be the curvature tensor of $\nabla$ and 
$\Phi_x$ be the restricted holonomy group of $\nabla$ at $x\,(\in M)$.  
Then the triple $((T_xM,J_x,\langle\,\,,\,\,\rangle_x),R_x,\Phi_x)$ is an 
anti-Kaehlerian holonomy system.  In particular, if $(M,J,\langle\,\,,\,\,
\rangle)$ is locally symmetric (resp. irreducible), then this 

anti-Kaehlerian holonomy system is symmetric (resp. weakly irreducible).  

\vspace{0.5truecm}

\noindent
{\it Example 2.}  Let $(M,J,\langle\,\,,\,\,\rangle)$ be a complex 
$n$-dimensional anti-Kaehlerian 
submanifold in an anti-Kaehlerian manifold $(\widetilde M,\widetilde J,
\langle\,\,,\,\,\rangle\,\,\widetilde{}\,\,)$, 
$T^{\perp}M$ be the normal bundle, $A$ be 
the shape tensor, $\nabla^{\perp}$ be the normal connection, $R^{\perp}$ be 
the curvature tensor of $\nabla^{\perp}$ and $\Phi^{\perp}_x$ be the 
restricted holonomy group of $\nabla^{\perp}$ at $x\,(\in M)$.  Define 
$\bar R^{\perp}_x\in T^{\perp}_xM^{\ast}\otimes T^{\perp}_xM^{\ast}\otimes 
T^{\perp}_xM^{\ast}\otimes T^{\perp}_xM$ by 
$$\bar R^{\perp}_x(v_1,v_2)v_3:=\sum_{i=1}^{2n}\langle e_i,e_i\rangle 
R^{\perp}_x(A_{v_1}e_i,A_{v_2}e_i)v_3,$$
where $(e_1,\cdots,e_{2n})$ is an orthonormal base of $T_xM$.  Then the 
triple $((T^{\perp}_xM,\widetilde J_x\vert_{T^{\perp}_xM},$\newline
$\langle\,\,,\,\,\rangle_x\widetilde{}\,\,
\vert_{T^{\perp}_xM\times T^{\perp}_xM}),
\bar R^{\perp}_x,\Phi^{\perp}_x)$ is an anti-Kaehlerian holonomy system.  

\vspace{0.5truecm}

We have the following fact for a weakly irreducible symmetric anti-Kaehlerian 
holonomy system.  

\vspace{0.5truecm}

\noindent
{\bf Lemma 4.1.}  {\sl Let $S=((V,J,\langle\,\,,\,\,\rangle),R,G)$ be a weakly 
irreducible symmetric anti-Kaehlerian holonomy system with $R\not=0$.  
Then the $G$-action on $V$ is equivalent to an aks-representation.}

\vspace{0.5truecm}

\noindent
{\it Proof.} Let $\mathfrak g^R$ be the Lie algebra generated by the set 
$\{R(v_1,v_2)\,\vert\,v_1,v_2\in V\}\,(\subset \mathfrak{so}_{AK}(V):=
{\rm Lie}(SO_{AK}(V)))$ and $G^R:=\exp\,\mathfrak g^R$, where $\exp$ is 
the exponential map of $SO_{AK}(V)$.  
Set $\mathfrak L:=\mathfrak g^R\oplus V$.  Define 
the $\sqrt{-1}$-multiples of elements of $\mathfrak L$ by $\sqrt{-1}v:=Jv$ 
($v\in V$) and $\sqrt{-1}R(v_1,v_2):=J\circ R(v_1,v_2)$ ($v_1,v_2\in V$).  
Also, define $[\,\,,\,\,]\,(:\mathfrak L\times\mathfrak L\to\mathfrak L)$ by 
$[A_1,A_2]:=A_1\circ A_2-A_2\circ A_1\,(A_1,A_2\in \mathfrak g^R),$ 
$[v_1,v_2]:=R(v_1,v_2)\,(v_1,v_2\in V)$ and $[A,v]:=A(v)\,(A\in 
\mathfrak g^R,\,v\in V)$.  Then it follows from the symmetricness of $S$ that 
$(\mathfrak L,[\,\,,\,\,])$ is a complex Lie algebra.  Define an (complex) 
involution $\rho$ of $(\mathfrak L,[\,\,,\,\,])$ by 
$\rho\vert_{\mathfrak g^R}={\rm id}$ and 
$\rho\vert_V=-{\rm id}$.  
Take a totally real subspace $W$ of $V$ such that $\langle\,\,,\,\,\rangle
\vert_{W\times JW}=0$ and that $\langle\,\,,\,\,\rangle\vert_{W\times W}$ 
is positive definite.  Let $(\mathfrak g^R)_W:=\{{\rm pr}_W\circ R(v_1,v_2)
\vert_W\,\vert\,v_1,v_2\in V\}$.  
By imitating the proof of Lemma 3.1, we can show that $(\mathfrak g^R)_W$ 
is a Lie subalgebra of $\mathfrak{so}(W)$ and $((\mathfrak g^R)_W)^{\bf c}
=\mathfrak g^R$.  
Thus $((\mathfrak L,[\,\,,\,\,]),\rho)$ is an anti-Kaehlerian 
symmetric Lie algebra.  
Let $(\widetilde L,\widetilde G)$ be a pair associated with 
$((\mathfrak L,[\,\,,\,\,]),\rho)$ such that $\widetilde L$ is simply 
connected and that $\widetilde G$ is connected.  According to 
Proposition 3.4, $(\widetilde L,\widetilde G)$ is an anti-Kaehlerian 
symmetric pair.  Hence, it follows from Proposition 3.3 that there exists 
an anti-Kaehlerian structure $(J,\langle\,\,,\,\,\rangle)$ such that 
$(\widetilde L/\widetilde G,J,\langle\,\,,\,\,\rangle)$ is an 
anti-Kaehlerian symmetric space.  On the other hand, 
we can show that the $G$-action on $V$ is equivalent to both the restricted 
holonomy group action $G^R$ of $\widetilde L/\widetilde G$ at $e\widetilde G$ 
and the linear isotropy group action ${\rm Ad}_{\widetilde L}
\vert_{T_{e\widetilde G}(\widetilde L/\widetilde G)}(\widetilde G)$ (see 
P359$\sim$360 of [W1]).  Since the $G$-action is weakly irreducible by the 
assumption, $\widetilde L/\widetilde G$ is irreducible.  Hence, 
${\rm Ad}_{\widetilde L}\vert_{T_{e\widetilde G}(\widetilde L/\widetilde G)}
(\widetilde G)$-action is an aks-representation.  Therefore, we obtain 
the statement of this lemma.  \hspace{9truecm}q.e.d.

\vspace{0.5truecm}

Now we shall define the notion of the complexification of a holonomy system.  
Let $S=((V,\langle\,\,,\,\,\rangle),R,G)$ be a holonomy system.  Then the 
triple $S^{\bf c}:=((V^{\bf c},{\rm Re}\langle\,\,,\,\,\rangle^{\bf c}),
R^{\bf c},G^{\bf c})$ gives an anti-Kaehlerian holonomy system, where 
$V^{\bf c},\,\langle\,\,,\,\,\rangle^{\bf c},\,R^{\bf c}$ and $G^{\bf c}$ 
are the complexifications of $V,\,\langle\,\,,\,\,\rangle,\,R$ and $G$, 
respectively.  We call this system $S^{\bf c}$ the {\it complexification of} 
$S$.  Next we shall define the notion of a totally real holonomy subsystem of 
an anti-Kaehlerian holonomy system.  Let 
$S=((V,J,\langle\,\,,\,\,\rangle),R,G)$ be an anti-Kaehlerian holonomy 
system.  Take a totally real subspace $W$ of $V$ such that 
$\langle\,\,,\,\,\rangle\vert_{W\times JW}=0$ and that 
$\langle\,\,,\,\,\rangle\vert_{W\times W}$ is positive definite.  
Set $R_W:={\rm pr}_W\circ R\vert_{W\times W\times W}$.  
Let $\mathfrak g_W$ be the Lie subalgebra of $\mathfrak{so}(W)$ spanned by 
$\{{\rm pr}_W\circ A\vert_W\,\vert\,A\in\mathfrak g\}$ and 
$G_W:=\exp_{SO(W)}(\mathfrak g_W)$.  It is shown that $G_W$ 
is compact and connected.  
Hence the triple $S_W:=((W,\langle\,\,,\,\,\rangle\vert_{W\times W}),R_W,G_W)$ 
is a holonomy system.  If $G_W^{\bf c}=G$, then we have $S_W^{\bf c}=S$.  
Then we call $S_W$ a {\it totally real holonomy subsystem of} $S$.  
Note that, if $S$ is symmetric and $R\not=0$, then $G_W^{\bf c}=G$ 
automatically holds.  In fact, the Lie algebra $\mathfrak g$ of $G$ is then 
generated by $\{R(v_1,v_2)\,\vert\,v_1,v_2\in V\}$ and the Lie algebra 
$\mathfrak g_W$ of $G_W$ includes $\{R_W(w_1,w_2)\,\vert\,w_1,w_2\in W\}$.  
Hence we have $\mathfrak g\subset\mathfrak g_W^{\bf c}$, that is, $G\subset 
G_W^{\bf c}$.  On the other hand, it is clear that $G_W^{\bf c}\subset G$.  
After all we have $G_W^{\bf c}=G$.  

Now we show the following fact for a weakly irreducible anti-Kaehlerian 
holonomy system.  

\vspace{0.5truecm}

\noindent
{\bf Lemma 4.2.} {\sl Let $S=((V,J,\langle\,\,,\,\,\rangle),R,G)$ be a 
weakly irreducible anti-Kaehlerian holonomy system.  Assume that there exists 
a totally real holonomy subsystem of $S$ having non-zero scalar curvature.  
Then the $G$-action on $V$ is equivalent to an aks-representation.}

\vspace{0.5truecm}

\noindent
{\it Proof.} Let $S':=((W,\langle\,\,,\,\,\rangle\vert_{W\times W}),
R_W,G_W)$ be a totally real holonomy subsystem of 
$S$ having non-zero scalar curvature, which is irreducible because $S$ 
is weakly irreducible.  According to the proof of Theorem 5 
of [Si], we can construct a non-zero curvature-like tensor $R'\,(:W\times W
\times W\to W)$ such that $((W,\langle\,\,,\,\,\rangle\vert_{W\times W}),
R',G_W)$ is a symmetric holonomy system.  Define 
$\psi:G\times V^3\to V$ by $\psi(g,v_1,v_2,v_3)=gR'^{\bf c}(g^{-1}v_1,
g^{-1}v_2)g^{-1}v_3-R'^{\bf c}(v_1,v_2)v_3$ ($(g,v_1,v_2,v_3)\in G\times 
V^3$), where $R'^{\bf c}$ is the complexification of $R'$.  
Since $\psi$ is holomorphic and $\psi=0$ over a totally real 
submanifold $G_W\times W^3$ of $G\times V^3$, we have $\psi\equiv0$ by 
the theorem of identity.  Then the triple $((V,J,\langle\,\,,\,\,\rangle),
R'^{\bf c},G)$ is a weakly irreducible symmetric anti-Kaehlerian holonomy 
system.  Hence we obtain the statement of this lemma by Lemma 4.1.  
\hspace{9truecm} q.e.d.

\section{Partial tubes with flat and abelian normal bundle} 
For a submanifold in a Riemannian symmetric space  of non-positive 
(or non-negative) curvature, 
M. Br$\ddot u$ck [B] defined a certain kind of partial tube 
with abelian normal bundle including the normal holonomy tube, where 
the submanifold is assumed to admit the $\varepsilon$-tube for a sufficiently 
small positive number $\varepsilon$.  In this section, we shall define 
the similar partial tube for an anti-Kaehlerian submanifold in a 
non-flat anti-Kaehlerian symmetric space of non-positive (or non-negative) 
curvature.  Let $M$ be an anti-Kaehlerian submanifold in such an 
anti-Kaehlerian symmetric space $N=G/K$.  
Let $\varepsilon_{\gamma}:={\rm inf}\{\vert r\vert\,\vert\,r:{\rm focal}\,\,
{\rm radius}\,\,{\rm of}\,\,M\,\,{\rm along}\,\,\gamma\}$, where $\gamma$ 
is a unit speed normal geodesic of $M$.  Denote by $\varepsilon^M_+$ (resp. 
$\varepsilon^M_-$) ${\rm inf}\{\varepsilon_{\gamma}\,\vert\,\gamma:{\rm unit}
\,\,{\rm speed}\,\,{\rm spacelike}\,\,{\rm (resp.}\,\,{\rm timelike)}\,\,
{\rm normal}\,\,{\rm geodesic}\}$.  Assume that $\varepsilon^M_+>0$ (resp. 
$\varepsilon^M_->0$).  
Denote the metric, the curvature tensor and 
the complex structure of $N$ by $\langle\,\,,\,\,\rangle,\,\widetilde R$ 
and $\widetilde J$, respectively.  
Fix $x_0\in M$.  
Let $\mathfrak C_{x_0}:=\{c:[0,1]\to M\,:\,{\rm a}\,\,{\rm piecewise}\,\,
{\rm smooth}\,\,{\rm path}\,\,{\rm with}\,\,c(0)=x_0\}$, 
$\Phi_{x_0}^0$ be the restricted normal holonomy group 
of $M$ at $x_0$ and $\mathfrak L_{x_0}$ be the Lie subalgebra of 
$\mathfrak{so}_{AK}(T^{\perp}_{x_0}M)$ generated by $\{P_c^{-1}\circ
{\rm pr}_{T^{\perp}_{c(1)}M}\circ\widetilde R_{c(1)}(P_cv_1,P_cv_2)\circ 
P_c\,\vert\,v_1,v_2\in T^{\perp}_{x_0}M,\,\,c\in\mathfrak C_{x_0}\}$, 
where $\mathfrak{so}_{AK}(T^{\perp}_{x_0}M):=\{A\in\mathfrak{gl}
(T^{\perp}_{x_0}M)\,\vert\,\langle Av_1,v_2\rangle_{x_0}+
\langle v_1,Av_2\rangle_{x_0}=0\,\,(\forall v_1,v_2\in T^{\perp}_{x_0}M),
\,\,[A,\widetilde J_{x_0}\vert_{T^{\perp}_{x_0}M}]=0\}$, $P_c$ is the 
parallel transport along $c$ with respect to the normal connection 
$\nabla^{\perp}$ of $M$ and ${\rm pr}_{T^{\perp}_{c(1)}M}$ is the orthogonal 
projection onto $T^{\perp}_{c(1)}M$.  Also, let 
$\widehat{\mathfrak L}_{x_0}$ be the Lie algebra generated by 
$\mathfrak L_{x_0}$ and ${\rm Lie}\,\Phi_{x_0}^0$.  Let $L_{x_0}:=\exp\,
\mathfrak L_{x_0}$ and $\widehat L_{x_0}:=\exp\,\widehat{\mathfrak L}_{x_0}$, 
where $\exp$ is the exponential map of $GL(T^{\perp}_{x_0}M)$.  Note that 
$L_{x_0}$ and $\widehat L_{x_0}$ are Lie subgroups of $SO_{AK}
(T^{\perp}_{x_0}M):=\{A\in GL(T^{\perp}_{x_0}M)\,\vert\,\langle Av_1,Av_2
\rangle_{x_0}=\langle v_1,v_2\rangle_{x_0}\,(\forall v_1,v_2\in 
T^{\perp}_{x_0}M),\,\,[A,\widetilde J_{x_0}\vert_{T^{\perp}_{x_0}M}]=0\}$.  
Set $\widetilde R_c:=P_c^{-1}\circ{\rm pr}_{T^{\perp}_{c(1)}M}\circ
\widetilde R_{c(1)}(P_c(\cdot),P_c(\cdot))\circ P_c$ for each 
$c\in \mathfrak C_{x_0}$.  
For each $c\in\mathfrak C_{x_0}$, it is clear that $S_c:=
(T^{\perp}_{x_0}M,\widetilde R_c,L_{x_0})$ is an anti-Kaehlerian holonomy 
system.  Fix $c_0\in\mathfrak C_{x_0}$ and a totally real subspace 
$W$ of $T^{\perp}_{x_0}M$ such that $\langle\,\,,\,\,\rangle_{x_0}
\vert_{W\times W}$ is positive definite.  
Let $\mathfrak L^W_{x_0}$ be the Lie subalgebra of $\mathfrak{so}(W)$ 
generated by $\{{\rm pr}_W\circ\widetilde R_c(v_1,v_2)\vert_W\,\vert\,
v_1,v_2\in V,\,\,c\in \mathfrak C_{x_0}\}$ and set $L^W_{x_0}:=\exp\,
\mathfrak L^W_{x_0}$, where $\exp$ is the exponential map of ${\rm GL}(W)$.  
The group $L^W_{x_0}$ is compact because it is a closed subgroup of 
the compact group $SO(W)$.  
Hence $S_{c_0}\vert_W:=((W,\langle\,\,,\,\,\rangle_{x_0}\vert_{W\times W}),
{\rm pr}_W\circ\widetilde R_{c_0}\vert_{W\times W\times W},L^W_{x_0})$ is a 
holonomy system.  Clearly we have $(\mathfrak L^W_{x_0})^{\bf c}
=\mathfrak L_{x_0}$, that is, $(L^W_{x_0})^{\bf c}=L_{x_0}$.  
It is shown that $S_{c_0}\vert_W$ is a totally real holonomy 
subsystem of $S_{c_0}$.  Let $W=W_0\oplus W_1\oplus\cdots\oplus W_k$ 
be the decomposition of $W$ such that $W_i$ ($i=0,1,\cdots,k$) are 
$L^W_{x_0}$-invariant, 
$L^{W_0}_{x_0}=\{{\rm id}_{W_0}\}$ and that $L^{W_i}_{x_0}$ ($i=1,\cdots,k$) 
are irreducible (non-trivial), where $L^{W_i}_{x_0}:=\{g\vert_{W_i}\,\vert\,
g\in L^W_{x_0}\}$ ($i=0,1,\cdots,k$).  
Let $V_i:=W_i\oplus JW_i(=W_i^{\bf c})$ ($i=0,1,\cdots,k$).  
Note that the Lie algebra of $L^{W_i}_{x_0}$ is equal to 
$\{{\rm pr}_{W_i}\circ\widetilde R_c(v_1,v_2)\vert_{W_i}\,\vert\,v_1,v_2\in V,
\,c\in{\mathfrak C}_{x_0}\}$.  Let $\mathfrak L^{V_i}_{x_0}$ ($i=0,1,\cdots,k$) be the Lie subalgebra of $\mathfrak{so}_{AK}(V_i)$ generated by 
$\{{\rm pr}_{V_i}\circ\widetilde R_c(v_1,v_2)\vert_{V_i}\,\vert\,v_1,v_2\in V,
\,c\in{\mathfrak C}_{x_0}\}$ and $L_{x_0}^{V_i}:=\exp\,
\mathfrak L_{x_0}^{V_i}$, where $\exp$ is the exponential map of 
${\rm GL}(V_i)$.  
Clearly we have $T^{\perp}_{x_0}M=V_0\oplus V_1\oplus\cdots\oplus V_k$ and 
$L^{V_i}_{x_0}=(L^{W_i}_{x_0})^{\bf c}$ ($i=0,1,\cdots,k$).  
Also, it is easy to show that $V_i$ ($i=0,1,\cdots,k$) are 
$L_{x_0}$-invariant, $L^{V_0}_{x_0}=\{{\rm id}_{V_0}\}$ and that 
$L_{x_0}^{V_i}$ ($i=1,\cdots,k$) are weakly irreducible (non-trivial).  

We have the following fact.  

\vspace{0.5truecm}

\noindent
{\bf Lemma 5.1.} {\sl 
The action of $L_{x_0}^{V_i}$ on $V_i$ is equivalent to 
an aks-representation.}

\vspace{0.5truecm}

\noindent
{\it Proof.} 
It is easy to show that $S_i:=(V_i,
{\rm pr}_{V_i}\circ\widetilde R_{c_0}\vert_{V_i\times V_i\times V_i},
L_{x_0}^{V_i})$ is a weakly irreducible anti-Kaehlerian holonomy system and 
that 
$(W_i,{\rm pr}_{W_i}\circ\widetilde R_{c_0}\vert_{W_i\times W_i\times W_i},
L^{W_i}_{x_0})$ is an irreducible totally real holonomy subsystem of $S_i$.  
Since $N$ is of non-positive (or non-negative) curvatures, we see that the 
scalar curvature of ${\rm pr}_{W_i}\circ\widetilde R_{c_0}
\vert_{W_i\times W_i\times W_i}$ does not vanish.  Hence, it follows from 
Lemma 4.2 that the $L_{x_0}^{V_i}$-action is equivalent to 
an aks-representation.  \hspace{0.1truecm} q.e.d.

\vspace{0.5truecm}

In similar to Lemma 3.3 of [B], we have the following statements.  

\vspace{0.5truecm}

\noindent
{\bf Lemma 5.2.} {\sl {\rm (i)} $V_i$ {\rm(}$i=0,1,\cdots,k${\rm)} are 
$\Phi^0_{x_0}$-invariant.

{\rm (ii)} $\Phi^0_{x_0}\vert_{V_i}\subset L_{x_0}^{V_i}$ {\rm(}
$i=1,\cdots,k${\rm)}, 
where $\Phi^0_{x_0}\vert_{V_i}:=\{g\vert_{V_i}\,\vert\,g\in\Phi^0_{x_0}\}$.  

{\rm (iii)} Let $W_0=W_{0,0}\oplus W_{0,1}\oplus\cdots\oplus W_{0,{\it l}}$ 
be the decomposition of $W_0$ such that $W_{0,j}$ \newline
{\rm(}$j=0,1,\cdots,{\it l}${\rm)} 
are $\Phi^0_{x_0}\vert_{W_0}$-invariant, $\Phi^0_{x_0}\vert_{W_{0,0}}
=\{{\rm id}_{W_{0,0}}\}$ and that $\Phi^0_{x_0}\vert_{W_{0,j}}$ 
{\rm(}$j=1,\cdots,{\it l}${\rm)} are irreducible, where 
$\Phi^0_{x_0}\vert_{W_{0,j}}
:=\{g\vert_{W_{0,j}}\,\vert\,g\in \Phi^0_{x_0}\}$ 
{\rm(}$j=0,1,\cdots,{\it l}${\rm)}.  
Set $V_{0,j}:=W_{0,j}^{\bf c}$ ($j=1,\cdots,{\it l}$).  
Then the $\Phi^0_{x_0}\vert_{V_{0,j}}$-action on $V_{0,j}$ is 
equivalent to an aks-representation {\rm(}$j=1,\cdots,{\it l}${\rm)}.}

\vspace{0.5truecm}

\noindent
{\it Proof.} From the definition of $L_{x_0}$, it follows that 
$\Phi^0_{x_0}$ is contained in the normalizer of $L_{x_0}$ in 
${\rm SO}_{\rm AK}(T^{\perp}_{x_0}M)$.  Hence $V_i$ ($i=0,1,\cdots,k$) are 
$\Phi_{x_0}^0$-invariant.  The group $\Phi_{x_0}^0\vert_{V_i}$ is contained 
in the normalizer $N(L_{x_0}^{V_i})$ of $L_{x_0}^{V_i}$ ($i\geq1$). 
On the other hand, according to Theorem 5 of [Si], the normalizer of 
$L^{W_i}_{x_0}$ coincides with oneself.  From this fact, 
$N(L^{V_i}_{x_0})=L^{V_i}_{x_0}$ follows.  Hence we have $\Phi_{x_0}^0
\vert_{V_i}\subset L^{V_i}_{x_0}$ ($\geq1$).  We define 
$\bar R^{\perp}_{x_0}\in T^{\perp}_{x_0}M^{\ast}\otimes T^{\perp}_{x_0}
M^{\ast}\otimes T^{\perp}_{x_0}M^{\ast}\otimes T^{\perp}_{x_0}M$ by 
$\bar R^{\perp}_{x_0}(v_1,v_2)v_3:=\sum\limits_{i=1}^{2n}
\langle e_i,e_i\rangle R^{\perp}_{x_0}(A_{v_1}e_i,A_{v_2}e_i)v_3$, where 
$(e_1,\cdots,e_{2n})$ is an orthonormal base of $T_{x_0}M$.  Let 
$(\bar R^{\perp}_{x_0})_{W_0}:={\rm pr}_{W_0}\circ\bar R^{\perp}_{x_0}
\vert_{W_0\times W_0\times W_0}$ and $(\Phi_{x_0}^0)_{W_0}$ be 
the image by the exponential map of the Lie subalgebra of 
$\mathfrak{so}(W_0)$ generated by 
$\{{\rm pr}_{W_0}\circ P_c^{-1}\circ R^{\perp}_{c(1)}(P_cX, P_cY)\circ P_c
\vert_{W_0}\,\vert\,X,Y\in T_{x_0}M,\,c\in\mathfrak C_{x_0}\}$.  
The triple $(W_0,(\bar R^{\perp}_{x_0})_{W_0},(\Phi_{x_0}^0)_{W_0})$ 
is a holonomy system.  Since $\widetilde R(w_1,w_2)=0$ for all 
$w_1,w_2\in W_0$, we have 
$$\langle R^{\perp}_{x_0}(X,Y)w_1,w_2\rangle=\langle[A_{w_2},A_{w_1}]X,Y
\rangle\,\,(X,Y\in T_{x_0}M,\,\,w_1,w_2\in W_0) \leqno{(5.1)}$$
by the Ricci equation.  By using this relation, we have 
$$\langle(\bar R^{\perp}_{x_0})_{W_0}(w_1,w_2)w_3,w_4\rangle
=\frac12{\rm Tr}([A_{w_1},A_{w_2}]\circ[A_{w_3},A_{w_4}])\,\,
(w_1,\cdots,w_4\in W_0). \leqno{(5.2)}$$
By imitating the proof of Theorem 3.1 of [O] (in terms of $(5.1)$ and 
$(5.2)$), we can show that the triples 
$S_{W_{0,j}}:=(W_{0,j},\,({\rm pr}_{W_{0,j}}\circ\bar R^{\perp}_{x_0})
\vert_{W_{0,j}\times W_{0,j}\times W_{0,j}},(\Phi_{x_0}^0)_{W_0}
\vert_{W_{0,j}})$ ($j=1,\cdots,{\it l}$) are holonomy systems having 
non-zero scalar curvature, where we use the fact that $N$ is of 
non-positive (or non-negative) curvature.  Also, it is clear that 
$S_{V_{0,j}}:=(V_{0,j},({\rm pr}_{V_{0,j}}\circ\bar R_{x_0}^{\perp})
\vert_{V_{0,j}\times V_{0,j}\times V_{0,j}},\Phi_{x_0}^0\vert_{V_{0,j}})$ 
($j=1,\cdots,{\it l}$) are weakly irreducible anti-Kaehlerian holonomy systems 
having $S_{W_{0,j}}$ as a totally real holonomy subsystem.  Hence it 
follows from Lemma 4.2 that the 
$\Phi_{x_0}^0\vert_{V_{0,j}}$-action ($j=1,\cdots,{\it l}$) is 
equivalent to an aks-representation.\hspace{0.6truecm} q.e.d.

\vspace{0.5truecm}

From these lemmas, we have the following fact directly.  

\vspace{0.5truecm}

\noindent
{\bf Theorem 5.3.}  {\sl There exists a decomposition $T^{\perp}_{x_0}M
=V_0\oplus V_1\oplus\cdots\oplus V_{{\it l}}\oplus V'_1\oplus\cdots\oplus 
V'_k$ of $T^{\perp}_{x_0}M$ such that $V_i$ {\rm(}$i=0,1,\cdots,{\it l}${\rm)} 
and $V'_i$ {\rm(}$i=1,\cdots,k${\rm)} are $\widehat L_{x_0}$-invariant, 
$\widehat L_{x_0}\vert_{V_0}=\{{\rm id}_{V_0}\}$, the 
$\widehat L_{x_0}\vert_{V_i}$-actions {\rm(}$i=1,\cdots,{\it l}${\rm)} 
and the $\widehat L_{x_0}\vert_{V'_i}$-actions {\rm(}$i=1,\cdots,k${\rm)} 
are equivalent to aks-representations, $\widehat L_{x_0}
\vert_{V_1\oplus\cdots\oplus V_{{\it l}}}=\Phi^0_{x_0}
\vert_{V_1\oplus\cdots\oplus V_{{\it l}}}$ and that $\widehat L_{x_0}
\vert_{V'_1\oplus\cdots\oplus V'_k}=L_{x_0}
\vert_{V'_1\oplus\cdots\oplus V'_k}$.}

\vspace{0.5truecm}

For $v_0\in T^{\perp}_{x_0}M$, define a subbundle $B_{v_0}(M)$ of 
$T^{\perp}M$ by 
$$B_{v_0}(M):=\{P_c(gv_0)\,\vert\,g\in\widehat L_{x_0},\,\,\,
c\in\mathfrak C_{x_0}\}$$
and $\widetilde B_{v_0}(M):=\exp^{\perp}(B_{v_0}(M))$, where $\exp^{\perp}$ 
is the normal exponential map of $M$.  
For each spacelike (resp. timelike) vector $v_0$ with 
$\vert\vert v_0\vert\vert\,<\,\varepsilon_+^M$ (resp. $\varepsilon^M_-$), 
$\widetilde B_{v_0}(M)$ is is an immersed submanifold, that is, a partial 
tube over $M$ whose fibre over $x_0$ is $\exp^{\perp}(\widehat L_{x_0}v_0)$.  
This partial tube $\widetilde B_{v_0}(M)$ is a notion similar to a partial 
tube defined for a submanifold in a Riemannian symmetric space of 
non-positive (or non-negative) curvature by M. Br$\ddot u$ck [B].  Denote by 
${\rm Hol}_{v_0}(M)$ the normal holonomy tube over $M$ through $v_0$.  
Clearly we have ${\rm Hol}_{v_0}(M)\subset\widetilde B_{v_0}(M)$.  Also, we 
have the following facts.  

\vspace{0.5truecm}

\noindent
{\bf Theorem 5.4.} {\sl Assume that $\widehat L_{x_0}v_0$ is a principal 
orbit of the $\widehat L_{x_0}$-action.  Then the following statements 
{\rm (i)}$\sim${\rm (iii)} hold.  

{\rm (i)} The normal connection of $\widetilde B_{v_0}(M)$ is flat, 

{\rm (ii)} $\widetilde B_{v_0}(M)$ has abelian normal bundle,

{\rm (iii)} Assume that $M$ is simply connected.  The $\widehat L_{x_0}$
-action and the normal parallel transport map of $M$ preserve the focal 
structure of $M$ if and only if $\widetilde B_{v_0}(M)$ is anti-Kaehlerian 
equifocal.  Then $M$ is a focal submanifold of $\widetilde B_{v_0}(M)$.}

\vspace{0.5truecm}

\noindent
{\it Proof.} These statements are shown by imitating the discussions in 
Sections $4.2\sim4.4$ of \newline
[B].\hspace{13.5truecm}q.e.d.

\section{Anti-Kaehlerian submanifolds with abelian normal bundle} 
Let $N=G/K$ be a semi-simple anti-Kaehlerian symmetric space.  Denote by 
$\langle\,\,,\,\,\rangle$ (resp. $J$) the metric (resp. the complex structure) 
of $N$.  Let $E$ be a vector bundle along a smooth curve $c:[0,1]\to N$ 
(i.e., $E\,:\,$ a subbundle of $c^{\ast}TN$) such that each fibre 
$E_t$ ($t\in[0,1]$) is an anti-Kaehlerian and abelian subspace of $T_{c(t)}N$ 
and that each $\exp_N(E_t)$ ($t\in[0,1]$) is properly embedded into $N$.  
Since $N$ is semi-simple, $\exp_N(E_t)$ is an anti-Kaehlerian cylindrical 
product.  There exists a totally real subspace $E_t^{\bf R}$ of $E_t$ 
such that $\exp_N(E_t^{\bf R})$ is a torus (with a flat pseudo-Riemannian 
metric).  Denote by $G$ the full holomorphical isometry 
group of $N$ newly.  
Also, denote by $K_t$ the isotropy group of $G$ at $c(t)$ and denote by 
$(K_t)_{v_0}$ the isotropy group of the linear isotropy action 
$K_t\times T_{c(t)}N\to T_{c(t)}N$ at $v_0\in T_{c(t)}N$.  Then we have the 
following fact.  

\vspace{0.5truecm}

\noindent
{\bf Lemma 6.1.} {\sl The set $E^{{\it l}}:=\displaystyle{
\mathop{\cup}_{t\in[0,1]}\{v_0\in E_t\,\vert\,{\rm dim}\,(K_t)_{v_0}
\leq{\it l}\}}$ is open in $E$ for each ${\it l}\in{\bf N}$.}

\vspace{0.5truecm}

\noindent
{\it Proof.} The statement of this lemma is shown by imitating the discussion 
in Page 81 of \newline
[PT].  \hspace{13.1truecm}q.e.d.

\vspace{0.5truecm}

Set ${\it l}_0:={\rm min}\{{\it l}\,\vert\,E^{\it l}\not=\emptyset\}$.  Fix 
$t_0\in[0,1]$ and $v_0\in E_{t_0}\cap E^{{\it l}_0}$.  By using some 
$J$-orthonormal frame field $(\widetilde v_1,J\widetilde v_1,\cdots,
\widetilde v_r,J\widetilde v_r)$ of $E$, we define maps $\psi_{t_0t}
:E_{t_0}\to E_t$ ($t\in[0,1]$) by $\psi_{t_0t}((\widetilde{v_i})_{t_0})=
(\widetilde{v_i})_t$ and $\psi_{t_0t}(J(\widetilde{v_i})_{t_0})=J
(\widetilde{v_i})_t$ ($i=1,\cdots,r$).  Let $v_t:=\psi_{t_0t}(v_0)$.  Let 
$I_0$ be the maximal sub-interval of $[0,1]$ containing $t_0$ such that 
$v_t\in E^{{\it l}_0}$ for all $t\in I_0$, which is open because 
$E^{{\it l}_0}$ is open in $E$.  Take a smooth curve $\hat c:I_0\to G$ 
satisfying $\hat c(t_0)=e$ ($e\,:\,$ the identity element of $G$) and 
$\hat c(t)(c(t_0))=c(t)$ for all $t\in I_0$.  Let $\widehat E_t:=\hat c
(t)_{\ast}^{-1}(E_t)$ and $h(t):=\hat c(t)_{\ast}^{-1}(v_t)$ ($t\in I_0$).  
Take a tubular neighborhood $T$ of the principal orbit $K_{t_0}v_0$ in 
$T_{c(t_0)}N$.  Let $I_1$ be the maximal sub-interval of $I_0$ containing 
$t_0$ satisfying $h(I_1)\subset T$ and define $\gamma:I_1\to K_{t_0}v_0$ by 
$h(t)\in S_{\gamma(t)}$ ($t\in I_1$), where $S_{\gamma(t)}$ is the slice of 
$K_{t_0}v_0$ through $\gamma(t)$.  Let $o:I_1\to K_{t_0}$ be a smooth curve 
such that $o(t_0)=e$ and $o(t)(v_0)=\gamma(t)$ for all $t\in I_1$.  
Then we can prove the following fact by imitating the proof of Lemma 5.2 
of [B].  

\vspace{0.5truecm}

\noindent
{\bf Lemma 6.2.} {\sl The set $\displaystyle{\mathop{\cup}_{t\in I_1}
o(t)^{-1}(\widehat E_t)}$ is contained in a maximal abelian anti-Kaehlerian 
subspace of $T_{c(t_0)}N$.}

\vspace{0.5truecm}

\noindent
{\it Proof.} Let $w\in S_{\gamma(t)}\cap\widehat E_t$.  From $w\in 
S_{\gamma(t)}$, we have $(K_{t_0})_w\subset(K_{t_0})_{\gamma(t)}$ (see Page 81 
of [PT]).  This together with ${\rm dim}(K_{t_0})_{v_0}={\it l}_0$ deduces 
that ${\rm dim}(K_{t_0})_w={\rm dim}(K_{t_0})_{\gamma(t)}$, which implies 
that $(K_{t_0})_w=(K_{t_0})_{\gamma(t)}$.  Let $\mathfrak a_t:=T^{\perp}
_{\gamma(t)}(K_{t_0}v_0)$ ($t\in I_1$), which is the maximal abelian 
anti-Kaehlerian subspace of $T_{c(t_0)}N$ containing $\gamma(t)$.  
Since $K_{t_0}w$ is parallel to 
$K_{t_0}v_0$ and $w\in S_{\gamma_t}$, we have $T^{\perp}_wK_{t_0}w=T^{\perp}
_{\gamma(t)}K_{t_0}v_0$.  Similarly we have $T^{\perp}_{h(t)}K_{t_0}h(t)
=T^{\perp}_{\gamma(t)}K_{t_0}v_0=\mathfrak a_t$, where we use 
$h(t)\in E^{{\it l}_0}$. 
Hence, since $\mathfrak a_t$ is the maximal abelian anti-Kaehlerian subspace 
containing $h(t)$, $h(t)\in \widehat E_t$ and $\widehat E_t$ is abelian, we 
have $\widehat E_t\subset \mathfrak a_t$, that is, $o(t)^{-1}\widehat E_t
\subset \mathfrak a_0$.  
Thus the statement of this lemma follows.  \hspace{7.5truecm}q.e.d.

\vspace{0.5truecm}

Furthermore we can show the following fact by imitating the proof of 
Lemma 5.3 of [B].  

\vspace{0.5truecm}

\noindent
{\bf Lemma 6.3.} {\sl The space $o(t)^{-1}(\widehat E_t)$ is independent of 
the choice of $t\in I_1$.}

\vspace{0.5truecm}

\noindent
{\it Proof.} According to Lemma 6.2, $\displaystyle{\mathop{\cup}_{t\in I_1}
o(t)^{-1}(\widehat E_t)}$ is contained in some maximal abelian anti-Kaehlerian 
subspace $\mathfrak a_0$ of $T_{c(t_0)}N$.  
Since $N$ is semi-simple, $\exp\,\mathfrak a_0$ is an anti-Kaehlerian 
cylindrical product.  There exists a totally real subspace $\mathfrak a_0
^{\bf R}$ of $\mathfrak a_0$ such that $\exp\,\mathfrak a_0^{\bf R}$ is 
a torus.  Denote $\exp\,\mathfrak a_0^{\bf R}$ by $T^k$ 
($k=\frac12{\rm rank}\,N$).  
Since $\exp\,E_t$ is an anti-Kaehlerian cylindrical product by the 
assumption, so is also $\exp(o(t)^{-1}(\widehat E_t))$.  Hence 
$\exp(o(t)^{-1}(\widehat E_t)\cap \mathfrak a_0^{\bf R})$ is a torus, which 
we denote by 
$T_t^r$ ($r=\frac12{\rm dim}\,E_t$).  Let $\{{\bf e}_1,\cdots,{\bf e}_k\}$ be 
the lattice of $T^k$.  Since $T_t^r$ is a sub-torus of $T^k$, the lattice 
of $T_t^r$ is expressed as $\{{\bf a}_i:=\sum\limits_{j=1}^ka_{ij}(t)
{\bf e}_j\,\vert\,i=1,\cdots,r\}$ ($a_{ij}(t)\in{\bf Z}$).  Furthermore, 
since $T_t^r$ variates continuously with respect to $t$, 
$a_{ij}$'s are continuous.  Hence, since $a_{ij}$'s are constant.  
Hence $T_t^r$ is 
independent of the choice of $t$.  This implies 
that $o(t)^{-1}(\widehat E_t)$ is independent of the choice of $t$. 
\hspace{3.7truecm}q.e.d.

\vspace{0.5truecm}

From this lemma, we have the following fact.  

\vspace{0.5truecm}

\noindent
{\bf Lemma 6.4.} {\sl There exists a smooth curve $w:I_1\to G$ with 
$w(t)_{\ast}E_{t_0}=E_t$ {\rm(}$t\in I_1${\rm)}.}

\vspace{0.5truecm}

\noindent
{\it Proof.} Define a smooth curve $w:I_1\to G$ by $w(t):=\hat c(t)\circ 
o(t)$ ($t\in I_1$).  This curve $w$ is a desired curve.  
\hspace{11.6truecm}q.e.d.

\vspace{0.5truecm}

Furthermore, we can show the following fact from this lemma.  

\vspace{0.5truecm}

\noindent
{\bf Lemma 6.5.} {\sl There exists a smooth curve $w:[0,1]\to G$ with 
$w(t)_{\ast}E_0=E_t$ {\rm(}$t\in[0,1]${\rm)}.}

\vspace{0.5truecm}

\noindent
{\it Proof.} Let $G_{2r}^{AK}(N):=\displaystyle{\mathop{\cup}_{x\in N}
\{\Pi\,\vert\,\Pi:2r{\rm -dimensional}\,\,{\rm anti-Kaehlerian}\,\,
{\rm subspace}\,\,{\rm of}\,\,T_xN\}}$, which is a submanifold of the 
Grassmann bundle of $N$ consisting of $2r$-dimensional subspaces.  The 
group $G$ acts on $G_{2r}^{AK}(N)$ naturally.  
Let $I_2$ be the maximal interval 
such that $t_0\in I_2$ and that $\displaystyle{\mathop{\cup}_{t\in I_2}
E_t\subset G(E_{t_0})}$.  From Lemma 6.4, it follows that $I_2$ is open.  
On the other hand, since $t\to E_t$ ($t\in[0,1]$) is a continuous curve in 
$G_{2r}^{AK}(N)$, $I_2$ is closed.  Therefore we have $I_2=[0,1]$, which 
implies that the above interval $I_1$ is equal to $[0,1]$.  
\begin{flushright}q.e.d.\end{flushright}

\vspace{0.5truecm}

Also we prepare the following lemma.  

\vspace{0.5truecm}

\noindent
{\bf Lemma 6.6.} {\sl Fix $t_0\in[0,1]$.  Let $g:(-\varepsilon,\varepsilon)
\to G$ be a smooth curve such that $g(0)=e$ and that 
$\displaystyle{\frac{d}{dt}\vert_{t=0}g(t)c(t_0)}$ is orthogonal to 
$E_{t_0}$, and $X$ be the vector field along $\exp\,E_{t_0}$ 
defined by $X_x:=\displaystyle{\frac{d}{dt}\vert_{t=0}g(t)x}$ 
{\rm(}$x\in\exp\,E_{t_0}${\rm)}.  
Then $X$ is a normal vector field of $\exp\,E_{t_0}$.} 

\vspace{0.5truecm}

\noindent
{\it Proof.} Denote by $X_{\bf R}^T$ the $T(\exp\,E_{t_0}^{\bf R})$-
component of 
$X\vert_{\exp\,E_{t_0}^{\bf R}}$.  Let $\gamma:{\bf R}\to\exp\,E_{t_0}^{\bf R}$ be a geodesic in $\exp\,E_{t_0}^{\bf R}$ (and hence $N$).  Define a map 
$\delta:(-\varepsilon,\varepsilon)\times{\bf R}\to N$ by 
$\delta(t,s)=g(t)\gamma(s)$.  Since $\delta$ is a geodesic variation, the 
variational vector field $\displaystyle{\frac{\partial \delta}{\partial t}
\vert_{t=0}\,(=X\circ\gamma)}$ is a Jacobi field along $\gamma$.  Hence 
$X^T_{\bf R}\circ\gamma$ is also a Jacobi field.  By using this fact, we have 
$$\frac{d^2}{ds^2}\langle X_{\bf R}^T\circ\gamma,\dot{\gamma}\rangle
=\langle\widetilde{\nabla}_{\dot{\gamma}}\widetilde{\nabla}_{\dot{\gamma}}
(X^T_{\bf R}\circ\gamma),\dot{\gamma}\rangle=-\langle\widetilde R(X^T_{\bf R}
\circ\gamma,\dot{\gamma})\dot{\gamma},\dot{\gamma}\rangle=0.$$
Hence we can express as $\langle X^T_{\bf R}\circ\gamma,\dot{\gamma}
\rangle(s)=as+b\,\,(a,b\in{\bf R})$.  Since $\gamma({\bf R})$ is contained 
in the compact set $\exp\,E_{t_0}^{\bf R}$, we have $\sup\vert\vert\langle 
X^T_{\bf R}\circ\gamma,\dot{\gamma}\rangle\vert\vert<\infty$.  Therefore, we 
see that $\langle X^T_{\bf R}\circ\gamma,\dot{\gamma}\rangle$ is constant.  Hence we have 
$\langle\widetilde{\nabla}_{\dot{\gamma}}(X^T_{\bf R}\circ\gamma),
\dot{\gamma}\rangle=0$.  Since this relation holds for any geodesic $\gamma$ 
in $\exp\,E_{t_0}^{\bf R}$, $X^T_{\bf R}$ is a Killing vector field on a flat 
torus $\exp\,E_{t_0}^{\bf R}$.  This fact together with 
$(X^T_{\bf R})_{c(t_0)}=0$ implies 
that $X^T_{\bf R}\equiv0$.  Denote by $X^T$ the $T(\exp\,E_{t_0})$-component 
of $X$.  We have only to show $X^T\equiv0$.  Since $X^T$ is real holomorphic 
(i.e., $X^T-\sqrt{-1}JX^T\,:\,$ holomorphic) and $X^T_{\bf R}=0$ 
on the totally real submanifold $\exp\,E_{t_0}^{\bf R}$ of 
$\exp\,E_{t_0}$, we see that $X^T=0$ 
along $\exp\,E_{t_0}^{\bf R}$.  Furthermore, it follows from the theorem of 
identity that $X^T=0$ on the whole of $\exp\,E_{t_0}$.  This completes 
the proof.  
\hspace{5.3truecm}q.e.d.

\vspace{0.5truecm}

Let $M$ be an anti-Kaehlerian submanifold with abelian normal bundle in 
$N$.  Assume that $\exp_N(T_x^{\perp}M)$ is properly embedded for each 
$x\in M$.  
By using Lemma 6.6, we can show the following fact.  

\vspace{0.5truecm}

\noindent
{\bf Lemma 6.7.} {\sl Let $x$ be a point of $M$ and $g:{\bf R}\to G$ 
be a $C^{\infty}$-curve such that $g(0)=e,\,\,g(t)x\in M\,(t\in{\bf R})$ and 
that $g(t)_{\ast}T^{\perp}_xM=T^{\perp}_{g(t)x}M$ $(t\in{\bf R})$.  Let 
$c(t):=g(t)x$ $(t\in{\bf R})$.  Then $g(t)_{\ast}:T^{\perp}_xM\to 
T^{\perp}_{c(t)}M$ is the parallel transport along $c\vert_{[0,t]}$ with 
respect to the normal connection $\nabla^{\perp}$ of $M$.}

\vspace{0.5truecm}

\noindent
{\it Proof.} Take an arbitrary $v\in T^{\perp}_xM$.  Let $\gamma_v$ be the 
geodesic in $\exp^{\perp}(T^{\perp}_xM)$ with $\dot{\gamma}_v(0)=v$ and 
define a map $\delta:{\bf R}^2\to N$ by $\delta(t,s):=g(t)(\gamma_v(s))$. 
Since $\delta_{\ast}(\frac{\partial}{\partial t})$ is a normal vector field 
of $\exp(T^{\perp}_{c(t)}M)$ by Lemma 6.6 and $\exp(T^{\perp}_{c(t)}M)$ is 
totally geodesic, we have 
$$\begin{array}{l}
\displaystyle{\widetilde{\nabla}_{\dot c}g(t)_{\ast}v=
\widetilde{\nabla}_{\frac{\partial}{\partial t}\vert_{s=0}}
\delta_{\ast}(\frac{\partial}{\partial s})=
\widetilde{\nabla}_{\frac{\partial}{\partial s}\vert_{s=0}}\delta_{\ast}
(\frac{\partial}{\partial t})}\\
\hspace{1.5truecm}\displaystyle{=\nabla^{\perp_t}_{\frac{\partial}{\partial s}
\vert_{s=0}}\delta_{\ast}(\frac{\partial}{\partial t})\in T^{\perp}_{c(t)}
\exp(T^{\perp}_{c(t)}M)=T_{c(t)}M,}
\end{array}$$
where $\nabla^{\perp_t}$ is the normal connection of $\exp(T^{\perp}_{c(t)}M)$. Hence we have $\nabla^{\perp}_{\dot c}g(t)_{\ast}v=0$.  From the arbtrariness 
of $v$, this implies that $g(t)_{\ast}:T^{\perp}_xM\to T^{\perp}_{c(t)}M$ is 
the parallel transport along $c\vert_{[0,t]}$ with respect to 
$\nabla^{\perp}$.  
\hspace{8.6truecm}q.e.d.

\vspace{0.5truecm}

By using Lemmas 6.5 and 6.7, we can show the following fact.  

\vspace{0.5truecm}

\noindent
{\bf Theorem 6.8.} {\sl Let $M$ be as above.  The normal connection of 
$M$ is flat.}

\vspace{0.5truecm}

\noindent
{\it Proof.} Let $c:I\to M$ be a loop at $x(\in M)$ such that the homotopy 
class $[c]$ of $c$ is the identity element of the fundamental group 
$\pi_1(M,x)$.  
From the assumption, it follows that $t\to T^{\perp}_{c(t)}M$ satisfies the 
same conditions as the above $t\to E_t$.  Hence it follows from Lemma 6.5 
that there exists a smooth curve $w:I\to G$ with $w(t)_{\ast}(T^{\perp}_xM)
=T^{\perp}_{c(t)}M$ ($t\in I$).  Furthermore, it follows from Lemma 6.7 
that $w(1)_{\ast}:T^{\perp}_xM\to T^{\perp}_xM$ is the parallel transport 
along $c$ with respect to $\nabla^{\perp}$.  The element $w(1)$ of $G$ is 
an isometry of the anti-Kaehlerian cylindrical product $\exp_N(T^{\perp}_xM)$ 
having $x$ as a fixed point.  Furthermore, since $[c]$ is the identity 
element of $\pi_1(M,x)$, $w(1)$ preserves the orientation.  Hence, since the 
full orientation-preserving isometry group of an anti-Kaehlerian cylindrical 
product is a free action, $w(1)$ is the identity transformation of 
$\exp_N(T^{\perp}_xM)$ and hence $w(1)_{\ast}$ (i.e., the parallel transport 
along $c$ with respect to $\nabla^{\perp}$) is the identity transformation of 
$T^{\perp}_xM$.  
From the arbitrariness of $c$, it follows that the restricted normal 
holonomy group of $M$ at $x$ is trivial, that is, the normal connection of 
$M$ is flat.  
\hspace{4.6truecm}q.e.d.

\section{Proofs of Theorems A, B and C} 
Let $M$ and $F$ be as in Theorem A.  Fix $x_0\in F$  and $v_o\in 
T^{\perp}_{x_0}F$ with $\exp^{\perp}(v_o)\in M$.  
Without loss of generality, we may assume that $0<\langle v_o,v_o\rangle
<(\varepsilon_F^+)^2$ or $0>\langle v_o,v_o\rangle>
-(\varepsilon_F^-)^2$, 
where $\varepsilon_F^{\pm}$ is as in Section 5.  
Let $L_{x_0},\,\widehat L_{x_0},\,B_{v_o}(F)$ and $\widetilde B_{v_o}(F)$ 
be the quantities as in Section 5 defined for $F$.  
Let $\pi_F:M\to F$ be the focal map onto 
$F$ and $M_{x_0}^0$ be the component containing $v_o$ of 
$(\exp^{\perp_F}\vert_{T^{\perp}_{x_0}F})^{-1}(\pi_F^{-1}(x_0))$, where 
$\exp^{\perp_F}$ is the normal exponential map of $F$.  
Then we can show the following fact.  

\vspace{0.5truecm}

\noindent
{\bf Lemma 7.1.} {\sl The intersection $\widehat L_{x_0}v_o\cap M_{x_0}^0$ 
is open in $\widehat L_{x_0}v_o$.}

\vspace{0.5truecm}

\noindent
{\it Proof.} By imitating the proof of $(11)$ in Page 91 of [B], we can 
show the statement of this lemma.  
\hspace{12.7truecm}q.e.d.

\vspace{0.5truecm}

By using Theorem 5.3, Lemmas 6.5, 6.7 and 7.1, we prove Theorem A.  

\vspace{0.5truecm}

\noindent
{\it Proof of Theorem A.} We suffice to show that $M_{x_0}^0$ is an open 
potion of $\widehat L_{x_0}v_o$.  In fact, 
$M$ is then an open potion of $\widetilde B_{v_o}(F)$ and each fibre of 
$\widetilde B_{v_o}(F)$ are the image by the normal exponential map of 
a principal orbit of a pseudo-orthogonal representation on the normal space 
of $F$ which is equivalent to the direct sum representation of an 
aks-representation and a trivial representation by Theorem 5.3.  
Let $c:[0,1]\to M_{x_0}^0$ be a smooth curve with 
$c(0)=v_o$ and $v_1$ be an element of $T^{\perp}_{\exp^{\perp_F}(v_o)}M$ 
with $\exp^{\perp_M}(v_1)=x_0$.  
Let $\widetilde v_1$ be the $\nabla^{\perp}$-
parallel vector field along $\widetilde c:=\exp^{\perp_F}\circ c$ with 
$\widetilde v_1(0)=v_1$, where $\nabla^{\perp}$ is the normal connection 
of $M$.  Define a vector bundle $E$ along $\widetilde c$ by $E_t:=
T^{\perp}_{\widetilde c(t)}M$ ($t\in[0,1]$).  
For simplicity, set $N:=G/K$.  Since $E_t$ is an 
anti-Kaehlerian and abelian subspace of $T_{\widetilde c(t)}N$ and 
$\exp_N(E_t)$ is properly embedded by the assumption, 
it follows from Lemma 6.5 that there exists a smooth curve $w:[0,1]\to G$ 
with $w(t)(\exp^{\perp_F}(v_o))=\widetilde c(t)$ and 
$w(t)_{\ast}E_0=E_t$ ($t\in [0,1]$).  
Furthermore, it follows from Lemma 6.7 that $w(t)_{\ast}:E_0\to E_t$ is 
the parallel transport along $\widetilde c\vert_{[0,t]}$ with respect to 
$\nabla^{\perp}$.  Hence we have $w(t)_{\ast}v_1=\widetilde v_1(t)$.  
From this fact and $\exp^{\perp_M}(\widetilde v_1(t))=x_0$ ($t\in[0,1]$), 
we have 
$$w(t)(x_0)=w(t)(\exp^{\perp_M}(v_1))=\exp^{\perp_M}(w(t)_{\ast}v_1)=x_0,$$
that is, $w(t)\in K_{x_0}$, where $K_{x_0}$ is the isotropy group of $G$ at 
$x_0$.  Also, we have 
$$\exp_N(c(t))=\exp^{\perp_F}(c(t))=w(t)(\exp^{\perp_F}(v_o))
=\exp_N(w(t)_{\ast}(v_o))$$
and hence $c(t)=w(t)_{\ast}(v_o)\in K_{x_0}v_o$.  From the arbitrariness of 
$c$, it follows that 
$$M_{x_0}^0\subset K_{x_0}v_o. \leqno{(7.1)}$$
Let $\mathfrak H$ be the Lie subalgebra of $\mathfrak{so}_{\rm AK}
(T^{\perp}_{x_0}F)$ generated by the set $\{{\rm pr}_{T^{\perp}_{x_0}F}
\circ\widetilde R(v_1,v_2)\vert_{T^{\perp}_{x_0}F}\,\,\vert\,\,$\newline
$v_1,v_2\in T^{\perp}_{x_0}F\}$ and set 
$H:=\exp_{{\rm SO}_{\rm AK}(T^{\perp}_{x_0}F)}
\mathfrak H$, where $\exp_{{\rm SO}_{\rm AK}(T^{\perp}_{x_0}F)}$ is the 
exponential map of ${\rm SO}_{\rm AK}(T^{\perp}_{x_0}F)$.  Clearly we have 
$H\subset\widehat L_{x_0}$.  Let $v\in T^{\perp}_{v_0}Hv_0\cap T^{\perp}
_{x_0}F$.  Then we have $\langle\widetilde R(v_0,v)v_0,v\rangle=0$ because 
$\widetilde R(v_0,v)v_0\in T_{v_0}Hv_0$.  This implies that ${\rm Span}
\{v_0,v\}$ is an abelian subspace of $T^{\perp}_{x_0}F$.  Hence we see that 
${\rm Span}\{v_0,v\}\subset T^{\perp}_{v_0}(K_{x_0}v_0)$, that is, 
$v\in T^{\perp}_{v_0}(K_{x_0}v_0)$.  
From the arbitrariness of $v$, we have $T^{\perp}_{v_0}Hv_0\cap T^{\perp}
_{x_0}F\subset T^{\perp}_{v_0}(K_{x_0}v_0)$ and hence $T_{v_0}(K_{x_0}v_0)
\cap T^{\perp}_{x_0}F\subset T_{v_0}Hv_0$.  On the other hand, it follows 
from Lemma 7.1 and $(7.1)$ that 
$$T_{v_0}Hv_0\subset T_{v_0}(\widehat L_{x_0}v_0)\subset T_{v_0}M_{x_0}^0
\subset T_{v_0}(K_{x_0}v_0)\cap T^{\perp}_{x_0}F.$$
Therefore, we obtain $T_{v_0}(\widehat L_{x_0}v_0)=T_{v_0}M_{x_0}^0$.  
Similarly, we obtain $T_v(\widehat L_{x_0}v_0)=T_vM_{x_0}^0$ for other 
$v\in M_{x_0}^0$.  Hence we see that $M_{x_0}^0$ is an open potion of 
$\widehat L_{x_0}v_0$.  This completes the proof.  
\hspace{13.1truecm}q.e.d.

\vspace{0.5truecm}

Next we prepare the following lemma to prove Theorem B.  

\vspace{0.5truecm}

\noindent
{\bf Lemma 7.2.} {\sl Let $\pi^{\bf c}:G^{\bf c}\to G^{\bf c}/K^{\bf c}$ be 
the natural projection, $\phi^{\bf c}:H^0([0,1],\mathfrak g^{\bf c})\to 
G^{\bf c}$ be the parallel transport map for $G^{\bf c}$ and $H_u$ be the 
horizontal space of the submersion $\pi^{\bf c}\circ\phi^{\bf c}$ at 
$u\,(\in H^0([0,1],\mathfrak g^{\bf c})$).  Then the restriction 
$(\pi^{\bf c}\circ\phi^{\bf c})\vert_{H_u}$ of $\pi^{\bf c}\circ\phi^{\bf c}$ 
to $H_u$ is regarded as the exponential map of $G^{\bf c}/K^{\bf c}$ at 
$(\pi^{\bf c}\circ\phi^{\bf c})(u)$ under the identification of $H_u$ with 
$T_{(\pi^{\bf c}\circ\phi^{\bf c})(u)}(G^{\bf c}/K^{\bf c})$.}

\vspace{0.5truecm}

\noindent
{\it Proof.} Let $\gamma\,(:{\bf R}\to G^{\bf c}/K^{\bf c})$ be a geodesic in 
$G^{\bf c}/K^{\bf c}$ and $\gamma_u^L$ be the horizontal lift of $\gamma$ to 
$u\,\in(\pi^{\bf c}\circ\phi^{\bf c})^{-1}(\gamma(0))$.  Since 
$\pi^{\bf c}\circ\phi^{\bf c}$ is an anti-Kaehlerian submersion, $\gamma_u^L$ 
is a geodesic in $H^0([0,1],\mathfrak g^{\bf c})$.  
Since $H^0([0,1],\mathfrak g^{\bf c})$ is a flat space, we have 
$\gamma_u^L(t)=u+t\dot{\gamma}_u^L(0)\,(\in H_u)$, where $t\in{\bf R}$.  
From this fact, the statement of this lemma follows.  \hspace{2.6truecm}q.e.d.

\vspace{0.5truecm}

\noindent
{\it Proof of Theorem B.} 
Let $M\hookrightarrow G/K$ be as in the statement of Theorem B and 
$M^{\bf c}\hookrightarrow G^{\bf c}/K^{\bf c}$ be the (extrinsic) 
complexification of $M$, where we note that $G^{\bf c}/K^{\bf c}$ is a 
semi-simple anti-Kaehlerian symmetric space of non-positive curvature.  
Define a distribution $E_0$ on $M^{\bf c}$ by $(E_0)_x:=
\displaystyle{\mathop{\cap}_{v\in T^{\perp}_xM^{\bf c}}\left(
{\rm Ker}\,A^{\bf c}_v\cap{\rm Ker}\,R^{\bf c}(\cdot,v)v\right)}$ ($x\in 
M^{\bf c}$), where $A^{\bf c}$ is the shape tensor of $M^{\bf c}$ and 
$R^{\bf c}$ is the curvature tensor of $G^{\bf c}/K^{\bf c}$.  Then 
$M^{\bf c}$ is an open potion of a product submanifold ${M^{\bf c}}'\times 
G^{\bf c}_0/K^{\bf c}_0\,(\subset {G^{\bf c}}'/{K^{\bf c}}'\times 
G^{\bf c}_0/K^{\bf c}_0=G^{\bf c}/K^{\bf c})$, where 
the decomposition ${G^{\bf c}}'/{K^{\bf c}}'\times G^{\bf c}_0/K^{\bf c}_0$ 
is an anti-Kaehlerian product such that 
the distribution $T(G^{\bf c}_0/K^{\bf c}_0)$ on 
${M^{\bf c}}'\times G^{\bf c}_0/K^{\bf c}_0$ is the extension of $E_0$ and 
${M^{\bf c}}'$ is an anti-Kaehlerian equifocal submanifold in 
${G^{\bf c}}'/{K^{\bf c}}'$.  Denote ${M^{\bf c}}'\times 
G^{\bf c}_0/K^{\bf c}_0$ by $M^{\bf c}$ newly and $T(G^{\bf c}_0/K^{\bf c}_0)$ 
by $E_0$ newly.  
Fix $x\in M^{\bf c}$.  Since $M^{\bf c}$ is proper anti-Kaehlerian 
equifocal, the focal set $F$ of $M^{\bf c}$ at $x$ consists of 
infinitely many complex hyperplanes 
$\{{\it l}_{\lambda}\}_{\lambda\in\Lambda}$ in 
$T^{\perp}_x(M^{\bf c})$.  Take a focal normal vector field $v$ 
such that $v_x\in{\it l}_{\lambda_0}$ for some $\lambda_0\in I$ and that 
$v_x\notin {\it l}_{\lambda}$ ($\lambda\in I\setminus\{\lambda_0\}$).  
Denote by $E$ the focal distribution for $v$.  
Now we shall show that each leaf of $E$ is the image by the normal 
exponential map of an open potion of a complex sphere of a normal space of the 
focal submanifold $F:=f_v(M^{\bf c})$, where $f_v$ is the focal map for $v$.  
Let $L$ be a leaf of $E$.  Denote by $\widetilde E$ the focal distribution on 
$(\pi^{\bf c}\circ\phi^{\bf c})^{-1}(M^{\bf c})$ corresponding to $E$.  
Set $\widetilde F:=(\pi^{\bf c}\circ\phi^{\bf c})^{-1}(F)$, which is a focal 
submanifold corresponding to $\widetilde E$.  It is clear that $L$ is the 
image of some leaf $\widetilde L$ of $\widetilde E$ by $\pi^{\bf c}\circ
\phi^{\bf c}$.  According to Theorem 2 of [K2], 
$\widetilde L$ is an open potion of a complex sphere in the normal space 
$T^{\perp}_{u_0}\widetilde F$ of $\widetilde F$ at some 
$u_0\in \widetilde F$.  
According to Lemma 7.2, 
$(\pi^{\bf c}\circ\phi^{\bf c})\vert_{T^{\perp}_{u_0}\widetilde F}$ is 
regarded as the normal exponential map 
$\exp^{\perp}_{(\pi^{\bf c}\circ\phi^{\bf c})(u_0)}$ of $F$ at 
$(\pi^{\bf c}\circ\phi^{\bf c})(u_0)$ under the identification of 
$T^{\perp}_{u_0}\widetilde F\,(\subset T_{u_0}H^0([0,1],\mathfrak g^{\bf c})
=H^0([0,1],\mathfrak g^{\bf c}))$ with $T^{\perp}_{(\pi^{\bf c}\circ
\phi^{\bf c})(u_0)}F$.  Therefore, we see that $L$ is the image of an open 
potion of a complex sphere in 
$T^{\perp}_{(\pi^{\bf c}\circ\phi^{\bf c})(u_0)}F$ by 
$\exp^{\perp}_{(\pi^{\bf c}\circ\phi^{\bf c})(u_0)}$.  
Let $\mathfrak E:=\{E_i\}_{i\in I}$ be the family of all focal 
distributions on $M^{\bf c}$ whose leaves are the images by the 
normal exponential map of open potions of complex spheres of the normal spaces 
of focal submanifolds.  
Then it follows from the above fact that $E_0\oplus\sum\limits_{i\in I}E_i
=TM^{\bf c}$.  Also, it is clear that $I$ is finite.  
Let $\mathfrak E=\{E_1,\cdots,E_k\}$.  
Take a focal normal vector field $v_1$ with ${\rm Ker}\,f_{v_1\ast}=E_1$ and 
that $F_1:=f_{v_1}(M^{\bf c})$.  Take $w_1\in T^{\perp}F_1$ with 
$\exp^{\perp_{F_1}}(w_1)\in M^{\bf c}$, where 
$\exp^{\perp_{F_1}}$ is the normal exponential map of $F_1$.  According to 
the proof of Theorem A, the partial tube $\widetilde B_{w_1}(F_1)$ includes 
$M^{\bf c}$ as an open potion.  It is clear that $\widetilde B_{w_1}(F_1)$ 
is proper anti-Kaehlerian equifocal.  
Let $\{\widetilde E_1,\cdots,\widetilde E_k\}$ be the family of 
all focal distributions of $\widetilde B_{w_1}(F_1)$ with 
$\widetilde E_i\vert_{M^{\bf c}}=E_i$ ($i=1,\cdots,k$).  
Take a focal normal vector field $v_2$ of $\widetilde B_{w_1}(F_1)$ with 
${\rm Ker}\,f_{v_2\ast}=\widetilde E_2$ and set $F_2:=f_{v_2}
(\widetilde B_{w_1}(F_1))$.  Take $w_2\in T^{\perp}F_2$ with 
$\exp^{\perp_{F_2}}(w_2)\in \widetilde B_{w_1}(F_1)$, where 
$\exp^{\perp_{F_2}}$ is the normal exponetial map of $F_2$.  According to 
the proof of Theorem A, the partial tube $\widetilde B_{w_2}(F_2)$ includes 
$\widetilde B_{w_1}(F_1)$ as an open potion.  It is clear that 
$\widetilde B_{w_2}(F_2)$ is proper anti-Kaehlerian equifocal.  
In the sequel, by repeating $(k-2)$-times the same process, we obtain 
the complete extension $\widehat{M^{\bf c}}$ of $M^{\bf c}$.  From this 
construction of $\widehat{M^{\bf c}}$ and Theorem A, the statements (i) and 
(ii) of Theorem B follow.  \hspace{12.6truecm} q.e.d.

\vspace{0.5truecm}

Next we prove Theorem C.  

\vspace{0.5truecm}

\noindent
{\it Proof of Theorem C.} Let $\{E_0,E_1,\cdots,E_k\}$ be as in the statement 
(ii) of Theorem B.  Fix $x=gK\in M$.  Since $M$ is curvature-adapted, 
each $(E_i)_x$ ($i=1,\cdots,k$) is expressed as 
$\displaystyle{(E_i)_x=\mathop{\oplus}_{(\lambda,\mu)\in S}
\left({\rm Ker}(A_v-\lambda\,{\rm id})\cap 
{\rm Ker}(R(\cdot,v)v-\mu\,{\rm id})\right)^{\bf c}}$ 
for some unit normal vector $v$ of $M$ at $x$, where $A$ is the shape tensor 
of $M$ and $R$ is the curvature tensor of $G/K$, $S$ is a subset of 
$\displaystyle{\left({\rm Spec}\,A_v\times{\rm Spec}\,R(\cdot,v)v\right)
\setminus\{(0,0)\}}$.  Hence we have 
$(E_i)_x\cap T_xM=\displaystyle{\mathop{\oplus}_{(\lambda,\mu)\in S}
\left({\rm Ker}(A_v-\lambda\,{\rm id})\cap 
{\rm Ker}(R(\cdot,v)v-\mu\,{\rm id})\right)}$.  
Also, we have 
$(E_0)_x\cap T_xM=\displaystyle{\mathop{\cap}_{v\in T^{\perp}_xM}
\left({\rm Ker}\,A_v\cap {\rm Ker}\,R(\cdot,v)v\right)}$.  
From these relations, the statement of Theorem C follows.  
\begin{flushright}q.e.d.\end{flushright}

\section{Examples} 
Let $M$ be a principal orbit of a Hermann type action $H\times G/K\to G/K$ 
and $\theta$ be the Cartan involution of $G$ with 
$({\rm Fix}\,\theta)_0\subset K\subset{\rm Fix}\,\theta$ and $\sigma$ be 
an involution of $G$ with $({\rm Fix}\,\sigma)_0\subset H\subset{\rm Fix}\,
\sigma$.  Without loss of generality, we may assume that $\sigma\circ\theta
=\theta\circ\sigma$.  It is shown that $M$ is proper complex equifocal and 
curvature-adapted (see [K3]).  
Denote by $A$ the shape tensor of $M$.  Then $H(eK)$ is 
a totally geodesic orbit (which is a singular orbit except for one case) of 
the $H$-action and $M$ is catched as a partial tube over $H(eK)$.  Let 
$L:={\rm Fix}(\sigma\circ\theta)$.  The submanifold 
$\exp^{\perp}(T^{\perp}_{eK}(H(eK)))$ is totally geodesic and it 
is isometric to the symmetric space $L/H\cap K$, where $\exp^{\perp}$ is the 
normal exponential map of $H(eK)$.  
Let $\mathfrak g,\,\mathfrak f$ and $\mathfrak h$ be the Lie algebras of 
$G,\,K$ and $H$.  Denote by the same symbols the involutions of $\mathfrak g$ 
associated with $\theta$ and $\sigma$.  Set $\mathfrak p:={\rm Ker}
(\theta+{\rm id})\,(\subset\mathfrak g)$ and $\mathfrak q:={\rm Ker}
(\sigma+{\rm id})\,(\subset\mathfrak g)$.  
Take $x:=\exp^{\perp}(\xi)=\exp_G(\xi)K
\in M\cap\exp^{\perp}(T^{\perp}_{eK}(H(eK)))$, where $\xi\in\mathfrak p
={\rm Ker}(\theta+{\rm id})\,(\subset\mathfrak g)$.  For simplicity, set 
$g:=\exp_G(\xi)$.  Let $\Sigma$ be the section of $M$ through $x$, which pass 
through $eK$.  Let $\mathfrak b:=T_{eK}\Sigma$, $\mathfrak a$ be a maximal 
abelian subspace of $\mathfrak p:=T_{eK}(G/K)$ containing $\mathfrak b$, 
$\triangle$ be the root system with respect to $\mathfrak a$ and 
$\mathfrak p=\mathfrak a+\sum\limits_{\alpha\in\triangle_+}
\mathfrak p_{\alpha}$ be the root space decomposition with respect to 
$\mathfrak a$.  
Set $\mathfrak p':=\mathfrak p\cap\mathfrak q (=T^{\perp}_{eK}(H(eK)))$.  
The orthogonal complement ${\mathfrak p'}^{\perp}$ of $\mathfrak p'$ in 
$\mathfrak p$ is equal to $\mathfrak p\cap\mathfrak h$.  
Set $\overline{\triangle}=\{\alpha\vert_{\mathfrak b}\,
\vert\,\alpha\in\triangle\,\,{\rm s.t.}\,\,\alpha\vert_{\mathfrak b}\not=0\}$, 
which is a root system by Theorem B of [K6].  Let $\overline{\triangle}_+$ be 
a positive root system of $\overline{\triangle}$ with respect to some 
lexicographic ordering, 
$\mathfrak p_{\beta}:=\sum\limits_
{\alpha\in\triangle_+\,\,{\rm s.t.}\,\,\alpha\vert_{\mathfrak b}=\pm\beta}
\mathfrak p_{\alpha}$ for $\beta\in\overline{\triangle}_+$, 
$\overline{\triangle}_+^H:=\{\beta\in\overline{\triangle}_+\,\vert\,
{\mathfrak p'}^{\perp}\cap\mathfrak p_{\beta}\not=\{0\}\}$ and 
$\overline{\triangle}_+^V:=\{\beta\in\overline{\triangle}_+\,\vert\,
\mathfrak p'\cap\mathfrak p_{\beta}\not=\{0\}\}$.  Since both $\mathfrak p'$ 
and ${\mathfrak p'}^{\perp}$ are Lie triple systems of $\mathfrak p$ and 
$\mathfrak b$ is contained in $\mathfrak p'$, we have 
${\mathfrak p'}^{\perp}=\mathfrak z_{{{\mathfrak p}'}^{\perp}}(\mathfrak b)
+\sum\limits_{\beta\in\overline{\triangle}_+^H}
({\mathfrak p'}^{\perp}\cap\mathfrak p_{\beta})$ and 
$\mathfrak p'=\mathfrak b+\sum\limits_{\beta\in\overline{\triangle}_+^V}
(\mathfrak p'\cap\mathfrak p_{\beta})$.  
Note that $\overline{\triangle}^V:=\overline{\triangle}^V_+\cup
(-\overline{\triangle}^V_+)$ is the root system of the symmetric space 
$L/H\cap K$.  
Take $\eta\in T^{\perp}_xM$.  For each $X\in{{\mathfrak p}'}^{\perp}\cap
\mathfrak p_{\beta}$ ($\beta\in\overline{\triangle}^H_+$), we have 
$A_{\eta}\widetilde X_{\xi}=-\beta(\bar{\eta})\tanh\beta(\xi)
\widetilde X_{\xi}$ (see the proof of Theorem B of [K3]), where 
$\widetilde X_{\xi}$ is the horizontal lift of $X$ to $\xi$ (see Section 3 
of [K3] about this definition) and $\bar{\eta}$ is the element of 
$\mathfrak b$ with $\exp^{\perp}_{\ast\xi}(\bar{\eta})=\eta$ (where 
$\bar{\eta}$ is regarded as an element of $T_{\xi}\mathfrak p'$ 
under the natural identification of $\mathfrak p'$ with 
$T_{\xi}\mathfrak p'$.  
Also, for each $Y\in T_x(M\cap\exp^{\perp}(\mathfrak p'))
\cap g_{\ast}\mathfrak p_{\beta}$ ($\beta\in\overline{\triangle}^V_+$), 
we have $A_{\eta}Y=-\frac{\beta(\bar{\eta})}{\tanh\beta(\xi)}Y$ 
(see the proof of Theorem B of [K3]).  By using these relations, for the focal 
set $F$ of $\widehat{M^{\bf c}}$ at $x$, we have 
$$\begin{array}{l}
\displaystyle{g_{\ast}^{-1}F=\left(
\mathop{\cup}_{\beta\in\overline{\triangle}_+^V}
\mathop{\cup}_{j\in{\bf Z}}(-\xi+(\beta^{\bf c})^{-1}(j\pi\sqrt{-1}))
\right)}\\
\hspace{1.7truecm}\displaystyle{\cup\left(
\mathop{\cup}_{\beta\in\overline{\triangle}_+^H}
\mathop{\cup}_{j\in{\bf Z}}(-\xi+(\beta^{\bf c})^{-1}((j+\frac12)\pi
\sqrt{-1}))\right),}
\end{array}\leqno{(8.1)}$$
where $\beta^{\bf c}$ is the complexification of $\beta$.  
Let $FD^{cs}:=\{E_i\,\vert\,i=1,\cdots,k\}$ be the family of all focal 
distributions of $\widehat{M^{\bf c}}$ whose leaves are the images by the 
normal exponential map of complex spheres in the normal spaces of focal 
submanifolds and $FD^{cs}_x:=\{(E_i)_x\,\vert\,i=1,\cdots,k\}$.  For each 
$\beta\in\overline{\triangle}$, we set 
$$E^V_{\beta,x}:=g_{\ast}(\mathfrak p_{\beta}\cap\mathfrak p')^{\bf c}\,\,\,\,
(\beta\in\overline{\triangle}^V_+)\,\,\,\,\,\,{\rm and}\,\,\,\,\,\,
E^H_{\beta,x}:=g_{\ast}(\mathfrak p_{\beta}
\cap{\mathfrak p'}^{\perp})^{\bf c}\,\,\,\,
(\beta\in\overline{\triangle}^H_+).$$
Then we have 
$$\mathfrak z_{{\mathfrak p'}^{\perp}}(\mathfrak b)\oplus
\left(\mathop{\oplus}_{\beta\in\overline{\triangle}^V_+}
E^V_{\beta,x}\right)
\oplus
\left(\mathop{\oplus}_{\beta\in\overline{\triangle}^H_+}E^H_{\beta,x}\right)
=T_x\widehat{M^{\bf c}}. \leqno{(8.2)}$$
Also, for each subspace $E$ of 
$T_x\widehat{M^{\bf c}}$, we set $FN(E):=\{v\in T^{\perp}_x\widehat{M^{\bf c}}
\,\vert\,{\rm Ker}(f_{\widetilde v})_{\ast x}=E\}$, where $\widetilde v$ is 
the parallel normal vector field of $\widehat{M^{\bf c}}$ with 
$\widetilde v_x=v$ and $f_{\widetilde v}$ is the focal map for 
$\widetilde v$.  
For $\beta\in\overline{\triangle}^V_+$ with $2\beta,\,\frac12
\beta\notin\overline{\triangle}_+$, $E^V_{\beta,x}$ 
is a member of $FD^{cs}_x$ and, for $\beta'\in\overline{\triangle}^H_+$ with 
$2\beta',\frac12\beta'\notin\overline{\triangle}_+$, $E^H_{\beta',x}$ 
is a member of $FD^{cs}_x$.  In fact, 
$E^V_{\beta,x}$ (resp. $E^H_{\beta',x}$) is the focal distribution for a focal 
normal vector field $v$ with $v_x\in(-\xi+(\beta^{\bf c})^{-1}(0))\setminus
(g_{\ast}^{-1}F\setminus(-\xi+(\beta^{\bf c})^{-1}(0)))$ 
(resp. $v_x\in(-\xi+({\beta'}^{\bf c})^{-1}(\frac{\pi}{2}\sqrt{-1}))\setminus
(g_{\ast}^{-1}F\setminus(-\xi+({\beta'}^{\bf c})^{-1}
(\frac{\pi}{2}\sqrt{-1})))$).  Hence, according to Theorem 2 in [K2], we have 
$E^V_{\beta,x},\,E^H_{\beta',x}\in FD^{cs}_x$.  
However, for $\beta\in\overline{\triangle}^V_+$ with 
$2\beta\in\overline{\triangle}_+$ or $\frac12\beta\in\overline{\triangle}_+$, 
$E^V_{\beta,x}$ is not necessarily a member of $FD^{cs}_x$ but there exists 
$E\in FD^{cs}_x$ with $E\supset E^V_{\beta,x}$.  For example, if 
$\beta\in\overline{\triangle}^V_+,\,\frac12\beta\in\overline{\triangle}^H_+\cap
\overline{\triangle}^V_+$ and $2\beta\notin\overline{\triangle}_+$, then we 
have $E^V_{\beta,x}\notin FD^{cs}_x$ but $E^V_{\beta,x}\oplus 
E^H_{\frac12\beta,x}\in FD^{cs}_x$ and $E^V_{\beta,x}\oplus 
E^V_{\frac12\beta,x}\in FD^{cs}_x$.  
In fact, $E^V_{\beta,x}\oplus E^H_{\frac12\beta,x}$ (resp. 
$E^V_{\beta,x}\oplus E^V_{\frac12\beta,x}$) is the focal distribution for a 
focal normal vector field $v$ with 
$v_x\in(-\xi+(\beta^{\bf c})^{-1}(\pi\sqrt{-1}))\setminus
(g_{\ast}^{-1}F\setminus(-\xi+(\beta^{\bf c})^{-1}(\pi\sqrt{-1})))$ 
(resp. $v_x\in(-\xi+(\beta^{\bf c})^{-1}(0))\setminus
(g_{\ast}^{-1}F\setminus
(-\xi+(\beta^{\bf c})^{-1}(0)))$) but there exists no focal normal vector 
field having $E^V_{\beta,x}$ as a focal distribution.  
Similarly, for $\beta'\in\overline{\triangle}^H_+$ with 
$2\beta'\in\overline{\triangle}_+$ or $\frac12\beta'\in
\overline{\triangle}_+$, 
$E^H_{\beta,x}$ is not necessarily a member of $FD^{cs}_x$ but there exists 
$E'\in FD^{cs}_x$ with $E'\supset E^H_{\beta',x}$.  
Thus, if $\overline{\triangle}$ (which is the root system) is reduced, 
then we have 
$T\widehat{M^{\bf c}}=\oplus_{i=0}^kE_i$ (orthogonal direct sum), where $E_0$ 
is defined by $(E_0)_x:=\displaystyle{\mathop{\cap}_{v\in T^{\perp}_x
\widehat{M^{\bf c}}}\left({\rm Ker}A_v^{\bf c}\cap{\rm Ker}\,R^{\bf c}
(\cdot,v)v\right)}$ ($x\in\widehat{M^{\bf c}}$) and 
$\{E_1,\cdots,E_k\}=FD^{cs}_x$.  
However, if $\overline{\triangle}$ is not reduced, 
then we have $T\widehat{M^{\bf c}}=\sum_{i=0}^kE_i$ but 
the right-hand side is not necessarily an orthogonal direct sum.  
Assume that $\overline{\triangle}$ is reduced.  For each 
$i\in\{1,\cdots,k\}$, we have $(E_i)_x=E^V_{\beta,x}$ or 
$E^H_{\beta,x}$ for some $\beta\in\overline{\triangle}$.  It is easy to show 
that the leaves of $E_i^{\bf R}:=E_i\vert_M\cap TM$ are diffeomorphic to 
a sphere (resp. an affine space) in case of $(E_i)_x=E^V_{\beta,x}$ (resp. 
$E^H_{\beta,x}$).  After all $M$ is orthogonally netted by some foliations 
consisting of (topological) spheres and some foliations consisting of 
leaves which is diffeomorphic to an affine space.  

\vspace{1truecm}

\centerline{{\bf References}}

\vspace{0.5truecm}

{\small
\noindent
[B] M. Br$\ddot u$ck, Equifocal famlies in symmetric spaces of compact type, J.reine angew. Math. {\bf 515} 

(1999), 73--95.

\noindent
[BCO] J. Berndt, S. Console and C. Olmos, Submanifolds and holonomy, Research 
Notes in Mathe-

matics 434, CHAPMAN $\&$ HALL/CRC Press, Boca Raton, London, New York 
Washington, 

2003.

\noindent
[CP] M. Cahen and M. Parker, Pseudo-riemannian symmetric spaces, Memoirs 
of the Amer. Math. 

Soc. {\bf 24} No. 229 (1980).  

\noindent
[HLO] E. Heintze, X. Liu and C. Olmos, Isoparametric submanifolds and 
a Chevalley type rest-

riction theorem, Integrable systems, geometry, and topology, 151-190, 
AMS/IP Stud. Adv. 

Math. 36, Amer. Math. Soc., Providence, RI, 2006.

\noindent
[HOT] E. Heintze, C. Olmos and G. Thorbergsson, Submanifolds with constant 
principal curvatures 

and normal holonomy groups, Intern. J. Math. {\bf 2} (1991), 167--175.

\noindent
[H] S. Helgason, 
Differential geometry, Lie groups and symmetric spaces, Pure Appl. Math. 80, 

Academic Press, New York, 1978.

\noindent
[KN] S. Kobayashi and K. Nomizu, 
Foundations of differential geometry, Interscience Tracts in 

Pure and Applied Mathematics 15, Vol. II, New York, 1969.

\noindent
[K1] N. Koike, 
Submanifold geometries in a symmetric space of non-compact 
type and a pseudo-

Hilbert space, Kyushu J. Math. {\bf 58} (2004), 167--202.

\noindent
[K2] N. Koike, 
Complex equifocal submanifolds and infinite dimensional anti-
Kaehlerian isopara-

metric submanifolds, Tokyo J. Math. {\bf 28} (2005), 
201--247.

\noindent
[K3] N. Koike, 
Actions of Hermann type and proper complex equifocal submanifolds, 
Osaka J. 

Math. {\bf 42} (2005), 599--611.

\noindent
[K4] N. Koike, 
A splitting theorem for proper complex equifocal submanifolds, Tohoku Math. J. 

{\bf 58} (2006) 393-417.

\noindent
[K5] N. Koike, 
A Chevalley type restriction theorem for a proper complex equifocal submani-

fold, Kodai Math. J. {\bf 30} (2007) 280-296.

\noindent
[K6] N. Koike, 
On curvature-adapted and proper complex equifocal submanifolds, 
submitted 

for publication.  

\noindent
[O] C. Olmos, 
The normal holonomy group, Proc. Amer. Math. Soc. {\bf 110} 
(1990), 813--818.  

\noindent
[OS] T. Oshima and J. Sekiguchi, The restricted root system of a semisimple 
symmetric pair, 

Advanced Studies in Pure Math. {\bf 4} (1984), 433--497. 

\noindent
[PT] R. S. Palais and C. L. Terng, Critical point theory and submanifold 
geometry, Lecture Notes 

in Math. {\bf 1353}, Springer-Verlag, Berlin, 1988.

\noindent
[R] W. Rossmann, 
The structure of semisimple symmetric spaces, Can. J. Math. {\bf 1} 
(1979), 157--180.

\noindent
[Si] J. Simons, 
On the transitivity of holonomy systems, Ann. of Math. 
{\bf 76} (1962), 213--234.

\noindent
[Sz1] R. Sz$\ddot{{{\rm o}}}$ke, Complex structures on tangent 
bundles of Riemannian manifolds, Math. Ann. {\bf 291} 

(1991) 409-428.

\noindent
[Sz2] R. Sz$\ddot{{{\rm o}}}$ke, Automorphisms of certain Stein 
manifolds, Math. Z. {\bf 219} (1995) 357-385.

\noindent
[Sz3] R. Sz$\ddot{{{\rm o}}}$ke, Adapted complex structures and 
geometric quantization, Nagoya Math. J. {\bf 154} 

(1999) 171-183.

\noindent
[Sz4] R. Sz$\ddot{{{\rm o}}}$ke, Involutive structures on the 
tangent bundle of symmetric spaces, 
Math. Ann. {\bf 319} 

(2001), 319--348.

\noindent
[TT1] C. L. Terng and G. Thorbergsson, 
Submanifold geometry in symmetric spaces, J. Differential 

Geom. {\bf 42} (1995), 665--718.

\noindent
[TT2] C. L. Terng and G. Thorbergsson, 
Taut immersions into complete Riemannian manifolds, 

Tight and Taut submanifolds (Berkeley, Cal., 1994), 181--228, Math. 
Sci. Res. Inst. Publ. 

{\bf 32}, Cambridge Univ. Press, Cambridge, 1997.

\noindent
[W1] H. Wu, Holonomy groups of infinite metrics, Pacific J. Math. {\bf 20} 
(1967), 351--392.

\noindent
[W2] B. Wu, Isoparametric submanifolds of hyperbolic spaces, 
Trans. Amer. Math. Soc. {\bf 331} 

(1992), 609--626.
}

\vspace{1truecm}


\end{document}